\documentclass{amsart}
\usepackage{verbatim} 
\usepackage[utf8]{inputenc}
\usepackage{mathbbol,mathrsfs,amsfonts,amssymb,amsthm,stmaryrd}
\usepackage{amsmath}
\usepackage{enumitem}
\usepackage[all]{xypic}
\usepackage[mathcal]{euscript}
\usepackage{graphicx,supertabular}
\usepackage[dvipsnames]{xcolor}
\usepackage{extarrows}
\usepackage[colorlinks,hyperindex]{hyperref}
\hypersetup{linkcolor=blue, citecolor=blue}

\newcommand\GL{\operatorname{GL}}
\newcommand{\Prim}{{\operatorname{Prim}}}
\newcommand{\Ad}{\operatorname{Ad}}

\newcommand{\om}{\omega}

\newcommand{\ZZ}{\mathbb{Z}}
\newcommand{\TT}{\mathbb{T}}
\newcommand{\RR}{\mathbb{R}}
\newcommand{\C}{\mathbb{C}}
\newcommand{\R}{\mathbb{R}}
\newcommand{\T}{\mathbb{T}}
\newcommand{\Z}{\mathbb{Z}}
\newcommand{\HH}{{\mathcal H}}
\newcommand{\K}{{\mathcal K}}
\newcommand{\U}{{\rm U}}
\newcommand{\PU}{{\rm PU}}
\newcommand{\UU}{{\rm U(1)}}
\newcommand{\Aut}{{\rm Aut}}
\newcommand{\id}{{\rm id}}

\newcommand{\Ind}{{\rm Ind}}
\newcommand{\diag}{{\rm diag}}

\newcommand{\tra}{{\ \pitchfork\ }}
\newcommand{\tri}{{\ \!\!\pitchfork}}

\newcommand{\lk}{{\langle}}
\newcommand{\rk}{{\rangle}}

\newcommand{\fat}{{\mathbf 1\!\!\!\mathbf 1\!\!\!\mathbf 1}} 
\newcommand{\fatn}{{\mathbf 0\!\!\!\mathbf 0\!\!\!\mathbf 0}} 

\theoremstyle{theorem}

\newtheorem{thm}{Theorem}[section]

\newtheorem{prop}[thm]{Proposition}
\newtheorem{lem}[thm]{Lemma}
\newtheorem{cor}[thm]{Corollary}

\theoremstyle{definition}
\newtheorem{defi}[thm]{Definition}
\newtheorem{rem}[thm]{Remark}
\newtheorem{exa}[thm]{Example}

\begin{document}

\title[Non-Commutative T-Duality]{Non-Commutative T-Duality\\
{\footnotesize{The Dynamical Duality Theory and 2-Dimensional Examples}}}

\author[Echterhoff]{Siegfried 
  Echterhoff} 
\author[Schneider]{Ansgar Schneider}

\address{Westf\"alische Wilhelms-Universit\"at M\"unster\\ 
  Mathematisches Institut\\
  Einsteinstr. 62\\ 
  D-48149 M\"unster\\ 
  Germany} 
  \address{Georgstra{\ss}e 19\\
D-53111 Bonn\\
Germany}

\email{echters@math.uni-muenster.de, 
ansgar.schneider@uni-muenster.de}

\thanks{This work was partially supported by the Deutsche Forschungsgemeinschaft
(SFB 878)}

\begin{abstract}
A duality theory
of bundles of C$^*$-algebras
whose fibres are twisted transformation group algebras
is established.
Classical T-duality is obtained as a special case, 
where all fibres are commutative tori, 
i.e. untwisted group algebras for $\Z^n$. Our theory also includes 
the bundles considered by Mathai and Rosenberg in their work on non-commutative
T-duals, in which they allow twisted group algebras on one side of the duality.
\end{abstract}

\maketitle

\begingroup
\let \boldsymbol=\relax
\tableofcontents
\endgroup

\section{Introduction}
In classical (or topological) $T$-duality one starts with a circle bundle $p:E\to B$ over a locally compact  base space $B$ together with 
a class $\delta\in H^3(E,\ZZ)$  (often called $H$-flux). Then $T$-duality provides an involution $(E,\delta)\mapsto (E^{\#}, \delta^{\#})$ 
of circle bundles over $B$ with $H$-flux, satisfying a number of interesting properties. 
We refer to Rosenberg's CBMS memoir \cite{Ros} for a detailed axiomatic definition of this classical notion of $T$-duality (see \cite[Axiomatics 6.1.2]{Ros})
and an explanation  of how this relates to string theory. There is a completely topological construction of  $T$-duality due to Bunke and Schick \cite{BS}.
Another approach due to Raeburn and Rosenberg in \cite{RR} (see also \cite[Chapter 7]{Ros}) is completely C*-algebraic in nature: 
given the data $(E,\delta)$ as above, there is a unique stable continuous-trace
algebra $CT(E, \delta)$ (i.e., an algebra of continuous $C_0$-sections of a locally trivial field of compact operators over $E$) 
with Dixmier-Douady invariant 
$\delta\in H^3(E,\ZZ)$, and the results in \cite{RR} show that there is an essentially  unique $\RR$-action on $CT(E, \delta)$ which covers the 
action  of $\RR$ on $E$ given by inflating the given circle action.  Moreover,
they show that there is a unique circle bundle $E^{\#}$ over $B$ together with a class $\delta^{\#}\in H^3(E^{\#},\ZZ)$ 
such that 
$$CT(E,\delta)\rtimes \RR\cong CT(E^{\#}, \delta^{\#}).$$
Note that this C*-algebraic description of $T$-duality has the advantage that it gives a direct connection of the  (twisted) $K$-theories of $(E,\delta)$ and 
$(E^{\#},\delta^\#)$: by the Connes-Thom-isomorphism there is a natural isomorphism
$$K^*(E, \delta)=K_*(CT(E,\delta))\cong K_{*+1}(CT(E,\delta)\rtimes \RR)=K^{*+1}(E^{\#}, \delta^{\#}).$$
Note that Takai-duality provides an isomorphism 
$$CT(E^{\#}, \delta^{\#})\rtimes\widehat{\RR}\cong CT(E,\delta)\rtimes \RR\rtimes \widehat{\RR}\cong CT(E,\delta)\otimes \K(L^2(\RR))\cong CT(E,\delta).$$
which shows that this construction really provides an involution on the category of circle bundles over $B$ with $H$-flux.

Unfortunately, for higher dimensional torus bundles the theory becomes much more involved since in general there does not exist a (classical) dual 
pair $(E^\#, \delta^\#)$ for a given pair $(E,\delta)$ for a principal   $\TT^n$-bundle $p:E\to B$ if $n\geq 2$. 
However, in \cite{MR1} Mathai and Rosenberg show that in case $n=2$ there always exist (non-unique) non-commutative 
$T$-duals which are section algebras of bundles of stable non-commutative tori. Moreover, the possible duality pairs can be classified with the help of 
a certain Mackey-obstruction map. The construction extends the ideas explained above: 
recall that the non-commutative $2$-tori are just the 
twisted group algebras $C^*(\ZZ^2,\om)$ for $[\om]\in H^2(\ZZ^2, \U(1))$ (and similarly for higher dimensional tori).
Now given a pair $(E,\delta)$ as above in which $p:E\to B$ is a principal $\TT^2$-bundle,
there exists (in general a non-unique) action
of $\RR^2$ on $CT(E,\delta)$ which covers the given $\TT^2$-action on $E$ and then the $T$-dual is constructed as the 
crossed product $CT(E,\delta)\rtimes \RR^2$.  Indeed, writing $E$ locally as $U_i\times \TT^2$, the action of $\RR^2$ on the
invariant  ideal $CT(U_i\times\TT^2,\delta)\subseteq CT(E,\delta)$ is induced from an action of the lattice 
$\ZZ^2\subseteq \RR^2$ on $CT(U_i,\delta|_{U_i})$.  Hence the ideal 
 $CT(U_i\times \TT^2,\delta)\rtimes \RR^2$ is stably isomorphic to $CT(U_i,\delta|_{U_i})\rtimes \ZZ^2$, which is a section algebra of 
a  continuous C*-bundle with fibres $\mathcal K\rtimes_{\alpha_x}\ZZ^2\cong \mathcal K\otimes C^*(\ZZ^2, \om_x)$ at $x\in U_i$, 
where  $[\om_x]\in H^2(\ZZ^2,\U(1))$ is the Mackey obstruction 
 for implementing the action $\alpha_x$ on the fibre $\K$ at $x$ as the adjoint action of a homomorphism $U:\ZZ^2\to  \U(\mathcal H)$ (if $\K=\K(\mathcal H)$). 
 Later, in \cite{MR2} the authors extended these results to higher dimensional tori, where the obstruction for the 
existence of a possibly non-commutative $T$-dual is given by the requirement that the  class $\delta\in H^3(E,\ZZ)$ vanishes 
on the fibres of $p:E\to B$ (which is automatic if $n\leq 2$).  Again, the Connes-Thom-isomorphism always 
provides an isomorphism (up to dimension shift by $n$) of the relevant $K$-theories.

It is obvious that the resulting duality theory lacks symmetry: while on one side  
 we have a classical pair $(E,\delta)$, there is often a non-commutative torus bundle on the dual side. 
In this  paper we start to investigate a symmetric version of non-commutative T-duality in which we allow 
{\bf non-commutative torus bundles on both sides} of the duality. 

The framework we develop in this article is general enough to work for 
general locally compact abelian groups $G$ with discrete and 
co-compact subgroups $N\subset G$ and is not necessarily 
coupled to the motivating example $G=\R^n, N=\Z^n.$
Let us summarise some of its main content:
we introduce some notation and review some basic knowledge about C$^*$-dynamical 
systems in section \ref{SecPrelims}.
Then section
 \ref{SecNCToverPoint} starts with a recapitulation of 
Mathai's and Rosenberg's T-duality over the one-point space. It serves 
as a motivation for our  framework 
for general non-commutative T-duality over the one point space as it is
given in \ref{SecItCrPrAndTrans}.
This means to identify a subcategory $NCT(N;\widehat G)^\tri$ 
of the  category of all C$^*$-dynamical
systems whose objects are stable, twisted transformation 
group algebras of $N$, 
equipped with transverse $\widehat G$-actions (Definition \ref{DefiTransvers}).
This subcategory is dual to $NCT(N^\perp;G)^\tri$ 
by the duality functor $\_\rtimes G$
(Theorem \ref{ThmDualityOverThePoint}).
In particular, for $G=\R^n, N=\Z^n$, we obtain a 
self-duality of the category of stable, non-commutative tori
$NCT(\Z^n;\R^n)^\tri$ with transverse $\R^n$-actions.
In section \ref{SecClassifiRemarks} we construct a cohomological invariant 
\begin{eqnarray*}
[NCT(N;\widehat G)^\tri]\to H^2(N,\UU)
\end{eqnarray*}
on the set of isomorphism classes of these objects. 
This allows us to re-obtain the classical (i.e. commutative) subtheory inside 
our theory: it is the kernel of this map (Theorem \ref{ThmCommutativeSubTheory}).
Section \ref{Sec2DimNCT} describes the complete  picture of 
non-commutative T-duality in two and three dimensions, i.e. for $G=\R^n, N=\Z^n$, $n=2,3$,
over the one-point space.

In section \ref{SecGlobalNCTDuality} 
we turn to the global situation of non-commutative T-duality. 
The objects which we consider therein are bundles of
C$^*$-algebras whose fibres are stable, twisted group algebras
which satisfy a certain local triviality property which we call
$\omega$-triviality (Definition \ref{DefiOmegaTrivaial}).
For $G=\R^n, N=\Z^n$ the theory over the one-point space
directly generalises to the bundle case. 
In particular,
we establish a duality for these bundles, 
and we re-obtain the classical T-duality as a subcategory
characterised point-wise by
a trivial cohomological invariant.
For more general groups
some technical assumptions have to be made.

Section \ref{SecExamplesOverTheCircle} presents 
in detail  some examples of non-commutative T-duality over 
the circle.
In particular, we give an example of a locally 
$\omega$-trivial bundle which does not arise as 
a dual of a commutative bundle, thereby 
showing that in our framework the class of 
objects is bigger than in the approaches made so far.

\newpage
\section{Notation and Preliminaries}
\label{SecPrelims}

\subsection{Abelian Groups}
\label{SubSectAbGr}
Throughout this paper $G$ will always denote an abelian,
Hausdorff, second-countable, locally compact group.
Its dual group, the goup of characters
$\widehat G:={\rm Hom}(G,\UU)$, is equipped with the compact-open topology.
It is again an abelian, Hausdorff, second-countable and locally compact group \cite{Ru}. 
The bidual $\widehat{\!\widehat G}$
is canonically isomorphic to $G$, and we use both of the 
notations $\lk g,\chi\rk$ and $\lk\chi,g\rk$ to denote $\chi(g)\in \UU$, for $g\in G, \chi\in \widehat G$.
We assume that we have given a discrete and co-compact subgroup $N\subset G$. Then its annihilator 
$N^\perp:={\{\chi\in\widehat G:\chi|_N=1\}}$ is discrete and co-compact in $\widehat G$,
and there are cannonical identifications 
$\widehat N=\widehat G/N^\perp$ and $\widehat{N^\perp}=G/N.$

The most prominent and guiding example of groups 
that fit into this situation is for $n\in\mathbb N$ the self-dual example
$$\Z^n\hookrightarrow\R^n\twoheadrightarrow \T^n,$$
where $\T^n$ is the $n$-fold torus $\R^n/\Z^n$.
Another self dual example for $n\in \mathbb N$ is 
$$\mathbb Q^n\hookrightarrow\mathbb A^n\twoheadrightarrow \mathbb S^n,$$
where the (discrete) rational numbers $\mathbb Q$ sit inside the adeles $\mathbb A$, and 
the quotient is the solenoid $\mathbb S$, the dual group of the rationals (see \cite{HR}). 

\subsection{C*-Dynamical Systems}
By the term C*-algebra we will typically mean a separable C*-algebra.
The C*-auto\-morphism group $\Aut(A)$ of a C*-algebra $A$ is equipped with the
topology of point-wise convergence. An action $\alpha:G\to \Aut(A)$ is
just a continuous group homomorphism (usually called strongly continuous), 
and such a triple $(A,G,\alpha)$ is called a  {\em C*-dynamical system}.

If $(A,G,\alpha)$ is a C*-dynamical system (with $G$ abelian), 
then its {\em dual}
is  the system $(A\rtimes_\alpha G,\widehat G,\hat \alpha)$,
where the crossed product $A\rtimes_\alpha G$ is
the enveloping C*-algebra of the Banach *-algebra
$L^1(G,\alpha,^\star)$ which is {$L^1(G,A)$} equipped with {the} product
$$(f\ast f')(g):=\int_G f(h) \alpha_{h}(f'(g-h))\ dh$$
and with involution $f^\star(g):=\alpha_g(f(-g))^*$, for $f,f'\in {L^1(G,A)}$.
The dual action $\hat\alpha$ is given on the dense subspace
${L^1(G,A)}$ just by point-wise multiplication:
$\hat\alpha_\chi(f):=\lk\chi,\_\rk f(\_)$, for $\chi\in\widehat G,f\in {L^1(G,A)}$.

Recall the famous Takai duality theorem (e.g., see \cite[Theorem 7.1]{Wbook}):
\begin{thm}
\label{ThmTakai}
There is a $G$-equivariant isomorphism
\begin{eqnarray*}
\Big(\big(A\rtimes_\alpha G\big)\rtimes_{\hat\alpha}\widehat G,G,\hat{\hat{\alpha}}\Big)
\cong \Big(A\otimes\K(L^2(G)),G,\alpha\otimes ({\rm Ad}\circ \rho)\Big),
\end{eqnarray*}
where $\K(L^2(G))$ is the algebra of compact operators on $L^2(G)$, $\rho$ is the right
regular representation of $G$ on $L^2(G)$, and ${\rm Ad}$ denotes the conjugation
action of the unitary operators on the compacts. 
\end{thm}

\subsection{Morita Equivalent Actions}
Assume that $(A,G,\alpha)$ and $(B, G, \beta)$ are two C*-dynamical systems. 
Recall that a {\em Morita equivalence} 
$(E,\gamma)$ between $(A,G,\alpha)$ and $(B,G,\beta)$ consists of an $A$-$B$-equivalence bimodule 
$E$ together with an action $\gamma:G\to \Aut(E)$ which is compatible with the given actions $\alpha$ and $\beta$ in the sense that 
$$\alpha_g(_A\lk\xi,\eta\rk)\gamma_g(\zeta)=\gamma_g(_A\lk \xi,\eta\rk \zeta)=
\gamma_g( \xi\lk\eta, \zeta\rk_B)=\gamma_g(\xi)\beta_g(\lk \eta,\zeta\rk_B)$$
holds for all $\xi,\eta,\zeta\in E$ and $g\in G$. 
If $(E,\gamma)$ is such a Morita equivalence, 
then $C_c(G,E)$ becomes a $C_c(G,A)$-$C_c(G,B)$-bimodule by defining the left and right actions and the 
inner products by the convolution formulas
\begin{align*}
f\cdot \xi(g)&=\int_Gf(h)\gamma_g(\xi(g-h)\,dh\\
\xi\cdot f'(g)&=\int_G\xi(h)\beta_h(f'(g-h))\,dh\\
_{C_c(G,A)}\lk \xi,\eta\rk(g)&= \int_G  {}_A\lk \xi(h), \gamma_g(\eta(g-h)\rk\,dh\\
\lk \xi,\eta\rk_{C_c(G,B)}(g)&= \int_G \beta_{-g}(\lk \xi(h),\eta(g+h)\rk_B\,dh
\end{align*}
for all $\xi,\eta\in C_c(G,E), f\in C_c(G,A)$ and $f'\in C_c(G,B)$.  With these operations, $C_c(G,E)$ completes to a $(A\rtimes_\alpha G)$-$(B\rtimes_\beta G)$-equivalence bimodule $E\rtimes_\gamma G$ which 
equipped with the dual action $\widehat{\gamma}:G\to \Aut(E\rtimes_\gamma G)$ given by point-wise
multiplication,
$$(\widehat{\gamma}_\chi\xi)(g)=\lk g,\chi\rk\xi(g),
\quad \text{for all } \chi\in \widehat{G}, \xi\in C_c(G,E),$$
gives an equivariant Morita equivalence $(E\rtimes_\gamma G,\widehat{\gamma})$ for the dual systems
$(A\rtimes_{\alpha}G, \widehat{G},\widehat{\alpha})$ and $(B\rtimes_{\beta}G,\widehat{G},\widehat{\beta})$  (see \cite{Co, ech-dual} for further details).

\begin{exa}\label{dual-mor}
If there is an $\alpha$-$\beta$-equivariant 
isomorphism $\Phi:A\to B$, 
we consider $A$ as an $A$-$B$-bimodule with inner products given by 
$$
_A\lk a,b\rk=ab^*\quad\text{and}\quad \lk a,b\rk_B=\Phi(a^*b),
$$
then $(A,\alpha)$ gives a Morita equivalence between $(A,G,\alpha)$ and $(B,G,\beta)$.

Moreover, for any system $(A,G,\alpha)$, the system $\Big(A\otimes\K(L^2(G)),G,\alpha\otimes ({\rm Ad}\circ \rho)\Big)$
is Morita equivalent to $(A,G,\alpha)$ with bimodule $A\otimes L^2(G)$ and action 
$\alpha\otimes \rho:G\to \Aut(A\otimes L^2(G))$. Thus it follows from Takai's duality theorem that 
$(A,G,\alpha)$ is Morita equivalent to the double dual system $\Big(\big(A\rtimes_\alpha G\big)\rtimes_{\hat\alpha}\widehat G,G,\hat{\hat{\alpha}}\Big)$.
\end{exa}
The Morita equivalence between $\alpha$ and $\hat{\hat\alpha}$ is a primal case of  an important  Morita equivalence which will play a fundamental r\^ole 
in this 
paper. Let us review some other special cases of 
Morita equivalences which will appear in this work. 
(For reference, see \cite[Section 8.11]{P79}.)

\begin{exa}\label{exterior equivalence}
\mbox{}
(1) \em{Exterior Equivalence: }
If $\alpha$ and $\beta$ are actions on the same C*-algebra $A$, then they 
are called {\em exterior equivalent} if there exists 
a strictly continuous map to the unitary group of the multiplier algebra of $A$, $v:G\to \U{\rm M}(A); g\mapsto v_g$ such that 
$$
\alpha_g(a)=v_g\beta_g(a)v_g^*,\quad v_{g+h}=v_g\beta_g(v_h),\quad
 \text{for all } g,h\in G, a\in A.
$$
Let $E=A$ be the canonical $A$-$A$-equivalence and define $\gamma: G\to \Aut(E)$ 
by $\gamma_g(a)=v_g\beta_g(a)$ for $g\in G$ and $a\in E$. Then
it is easily checked that this implements a 
Morita equivalence between $(A,G,\alpha)$ and $(A,G,\beta)$.
We  say that {\em $v$ implements the exterior equivalence 
between $\alpha$ and $\beta$}.
In this case there is a canonical isomorphism
$$
\Phi_v: A\rtimes_{\beta}G\to A\rtimes_{\alpha}G, 
\quad  \Phi_v(f)(g)=f(g)v_g^*,\quad f\in C_c(G,A),
$$
which is $\widehat{\alpha}$-$\widehat{\beta}$-equivariant, 
hence induces an isomorphism for the dual actions.
\medskip
\\
(2) {\em Outer Conjugacy:} Two systems 
$(A,G,\alpha)$ and $(B,G,\beta)$ are {\em outer conjugate}, if there exists an isomorphism 
$\Phi:A\to B$ such that $\beta$ is exterior equivalent to the action $\alpha'=\Phi\circ \alpha\circ \Phi^{-1}$.
\medskip
\\
(3) {\em Stable Outer Conjugacy:} This is rather the most general form of Morita 
equivalence, in fact, it is an equivalent notion for C*-algebras which have
a countable approximate identity (or, equivalently, which contain strictly positive elements)
 \cite{Co}.
  Two systems 
$(A,G,\alpha)$ and $(B,G,\beta)$ are {\em stably outer conjugate} if
the two systems 
$(A\otimes\K,G,\alpha\otimes{\rm id})$ and $(B\otimes\K,G,\beta\otimes{\rm id})$ 
are outer conjugate, where $\K$ is the algebra of compact operators
on some separable Hilbert space.

\end{exa}

If $\alpha:G\to \Aut(A)$ is an action, we can restrict it to the subgroup $N\subset G$,
and then $C_c(G,A)$ completes to give an $(A\rtimes_{\alpha}N)$-Hilbert module, 
$E_N^G(A)$,  if the right action of $A\rtimes_\alpha N$ and the 
$A\rtimes_\alpha N$-valued inner products are defined on the level of $C_c(N,A)$ as 
\begin{align*}
\lk \xi,\eta\rk_{A\rtimes N}(n)&=\int_G \alpha_h(\xi(h))^*\eta(m-h)\,dh,\\
\xi\cdot f(g)&=\int_N \xi(g+n)\alpha_{g+n}(f(-n))\,dn
\end{align*}
for $\xi,\eta\in C_c(G,A), f\in C_c(N,A)$. 

There is a canonical left action of the dual system
$(A\rtimes_{\alpha}G, N^\perp,\widehat{\alpha})$ on $E_N^G(A)$ 
which is given by the covariant representation $(\Phi, U)$ 
in which $\Phi(f)\xi=f*\xi$ for $f,\xi\in C_c(G,A)$ and $U_\chi\xi=\chi\cdot\xi$ (point-wise multiplication). The integrated form $\Phi\times U: A\rtimes_{\alpha}G\rtimes_{\widehat{\alpha}}N^\perp
\to \mathcal{L}(E_N^G(A))$ defines a left action of 
$A\rtimes_{\alpha}G\rtimes_{\widehat{\alpha}}N^\perp$ on $E_N^G(A)$, which by 
\cite[Proposition 2.1]{ech-dual}  implements an isomorphism onto $\K(E_N^G(A))$ ($\mathcal L(\cdot)$ and $\K(\cdot)$ denote the adjointable and compact operators of a module, respectively). 
Of course, this is just a reformulation of Green's famous imprimitivity theorem \cite[Theorem 17]{Green}.
In \cite[Proposition 3.4 and Lemma 3.6]{ech-dual} it is shown that the resulting Morita equivalence
between $A\rtimes_\alpha G\rtimes_{\widehat{\alpha}} N^{\perp}$ and $A\rtimes_{\alpha}N$ {is} equivariant with respect to certain actions by $G$ and $\widehat{G}$. For notation, 
given an action $\alpha:G\to \Aut(A)$ for the abelian group $G$ and $N\subseteq G$ is a closed subgroup of $G$ we let 
 $\alpha^{dec}:G\to \Aut(A\rtimes_\alpha N)$ denote the action given by 
 $$\alpha^{dec}_h(f)(g)=\alpha_h(f(g))\quad g,h\in G, f\in C_c(G,A)$$
 and if $\beta:G/N\to \Aut(B)$ is an action of the quotient group, we denote by 
 $\inf\beta:G\to \Aut(B)$ 
 the inflation of $\beta$ to $G$. Combining the results of \cite[Proposition 3.4 and Lemma 3.6]{ech-dual} we then get
 
\begin{prop}\label{prop-dualmor}
In the above situation $E_N^G(A)$ becomes an $(A\rtimes_{\alpha}G\rtimes_{\widehat{\alpha}}N^\perp)$-$(A\rtimes_{\alpha}N)$-imprimitivity bimdule. Moreover, if we define actions $\gamma$ and $\widehat\gamma$
of $G$ and $\widehat{G}$ 
on $E_N^G(A)$ by
$$(\gamma_g\xi)(h):=\xi(g+h)\quad \text{and}\quad
(\widehat{\gamma}_\chi\xi)(g):=\langle g,\chi\rangle\xi(g)$$
for $\xi\in C_c(G,A)$, then $(E_N^G(A), \gamma)$ is a Morita equivalence 
between  
{$$(A\rtimes_{\alpha}G\rtimes_{\widehat{\alpha}}N^{\perp}, G, 
\inf\widehat{\widehat{\alpha}|_{N^\perp}})\quad \text{and}\quad (A\rtimes_\alpha N, G, \alpha^{dec})$$}
and 
$(E_N^G(A),\widehat{\gamma})$ is a Morita equivalence between
{$$(A\rtimes_{\alpha}G\rtimes_{\widehat{\alpha}}N^\perp, \widehat G, 
\widehat{\alpha}^{dec})\quad\text{and}\quad (A\rtimes_\alpha N, \widehat G, \inf\widehat{\alpha|_N}).$$}
\end{prop}

\subsection{Actions on $\K$, Twisted Group Algebras, and 2-Cocycles}
\label{SecZeugsUeberK}
Consider the short exact sequence 
$$
1\to\UU\to\U(\HH)\stackrel{\rm Ad}{\to} \PU(\HH)\to 1,
$$
where $\U(\HH)$ is the unitary group of some separable 
Hilbert space $\HH$, and $\PU(\HH)$
is the projective unitary group, the quotient of the unitaries by its center.
It induces a (not very long) exact sequence in Borel cohomology 
\begin{eqnarray}
\label{EqTheMaSequence}
\dots\to H^1(G,\U(\HH))\to H^1(G,\PU(\HH))\stackrel{\rm Ma}{\to}H^2(G,\UU) 
\end{eqnarray}
which terminates at $H^2(G,\UU)$ due to the non-commutativity of 
the involved coefficient groups. The (negative of the usual) 
connecting homomorphism\footnote{
\label{FootTheRightMackey}
We choose the convention ${\rm Ma}(\alpha):=-[\partial V]$, for a 
Borel lift $V:G\to\U(\HH)$ of $\alpha:G\to\PU(\HH)$.}
 ${\rm Ma}$
is called {\em Mackey obstruction}.
Now, because $\U(\HH)={\rm UM}(\K)$, the unitary group of the multiplier algebra
of the compacts $\K=\K(\HH)$, and because 
all automorphisms of $\K$ are inner, conjugation
defines a canonical isomorphism $\PU(\HH)={\rm Aut}(\K)$. 
So if  $\alpha:G\to \Aut(\K)$ is an action on 
the compacts, i.e. $\alpha\in H^1(G,\PU(\HH))$,
it defines a class ${\rm Ma}(\alpha)\in H^2(G,\UU)$, 
the { Mackey obstruction of $\alpha$}.
The action $\alpha$ is said to be {\em unitary}
if its Mackey obstruction vanishes. 
Note that any class $[\omega]\in H^2(G,\UU)$ arises as a Mackey obstruction 
of some action $\alpha:G\to \Aut(\K)$: 
Just put $\alpha_\omega(g):={\rm Ad}(L_\omega(g))$, where 
$L_\omega:G\to {\rm U}(L^2(G))$ is the left regular $\omega$-representation of $G$, i.e., 
\begin{eqnarray}
\label{EqTheLeftRegRep}
(L_\omega(g)\xi)(h)=
\omega(g, h-g)\xi(h-g),\quad \xi\in L^2(G), g,h\in G.
\end{eqnarray}
Then ${\rm Ma}(\alpha_\omega)=[\omega^{-1}]=-[\omega]$.

The following statement shows that actions on 
$\K$ are classified up to Morita (or exterior) equivalence
by the Borel 
cohomology group $H^2(G,\UU)$. 
We refer to \cite[{Section} 6.3]{CKRW} for a more general result.

\begin{prop}\label{prop-Mackey}
Suppose $\alpha,\beta:G\to \Aut(\K)$ are two actions of $G$ on $\K$. 
Then the following are equivalent:
\begin{enumerate}
\item $\alpha$ and $\beta$ are exterior equivalent.
\item $\alpha$ and $\beta$ are Morita equivalent.
\item ${\rm Ma}(\alpha)={\rm Ma}(\beta)\in H^2(G,\UU)$.
\end{enumerate}
\end{prop}

If ${\rm Ma}(\alpha)=[\omega]$, then the crossed product 
$\K\rtimes_\alpha G$ is isomorphic to\footnote{
Here we see that our sign convention of ${\rm Ma}$ is the right one, for 
otherwise we would obtain 
$$\K\otimes (\C\rtimes_{\omega^{-1}} G)\cong \K\rtimes_\alpha G.$$
} 
$\K\otimes(\C\rtimes_\omega G)$, where $\C\rtimes_\omega G$ denotes the twisted 
group C*-algebra of $G$ with respect to the cocycle $\omega$. This is the 
enveloping C*-algebra 
of the Banach *-algebra $L^1(G,\omega,^\star)$ 
given by the Banach space $L^1(G)$ with convolution and
involution given by
\begin{eqnarray}
\label{EqOpForTwistedAlg}
(f*f')(g)&:=&\int_G f(h)f'(g-h)\omega(h, g-h)\ dh\quad\text{and}\\
f^\star(g)&:=&\overline{\omega(g, -g)f(-g)}.\nonumber
\end{eqnarray}
The isomorphism
$\Psi:\K\otimes (\C\rtimes_\omega G)\to \K\rtimes_\alpha G$ is given on the level of 
$L^1$-functions by
\begin{eqnarray}
\label{EqIsoOfTwistedAndCrossed}
\Psi(k\otimes f)(g)=f(g)kV(g), \quad k\in \K, f\in L^1(G),
\end{eqnarray}
where $V:G\to\U(\HH)$ is a  Borel lift  of $\alpha$
(i.e. 1-cochain) such that its boundary is $\omega^{-1}$ (see \cite[Theorem 1.4.15]{Echt96}).
Note that this isomorphism is equivariant with respect to 
the canonical (dual) actions of $\widehat{G}$ on 
$\K\rtimes_{\alpha}G$ and  on $\K\otimes (\C\rtimes_\omega G)$ 
given by point-wise multiplication with characters.
It is an immediate consequence of Proposition \ref{prop-Mackey}
and of Takai duality that, up to equivariant Morita equivalence, 
the twisted group algebra equipped with the dual $\widehat{G}$-action 
is also classified  by the  cohomology class $[\omega]\in H^2(G,\U(1))$.

One should regard the twisted group algebra $\C\rtimes_\omega G$ 
as a deformation of $C_0(\widehat{G})\cong \C\rtimes_1 G$, 
where $1$ denotes the trivial cocycle on $G$. 
In particular, the twisted group algebras $\C\rtimes_\omega\Z^n$ 
are deformations 
of $C(\T^n)$, and they are called non-commutative $n$-tori. 
In this picture, the dual action 
of $\T^n$ on $\C\rtimes_\omega\Z^n$ 
is the analogue of the translation action of $\T^n$ on $C(\T^n)$
in the commutative case.
We have a natural isomorphism between 
the additive group $M^u(n,\R)$ of strictly upper triangular real matrices 
and $H^2(\R^n,\UU)$
which is given by sending a matrix {$A$}  to the class of the cocycle 
$\omega_A$ given by
$$
\omega_A(x,y):=\exp({2\pi i (A x)^{{t}} y}).
$$
Under this identification, the restriction map 
\begin{eqnarray*}
H^2(\R^n,\UU)\to H^2(\Z^n,\UU);\quad 
[\omega]\mapsto[\omega|_{\Z^n\times \Z^n}],
\end{eqnarray*}
 which is surjective, has kernel  given by the collection 
 of all classes corresponding to the set 
$M^u(n,\Z)$ of strictly upper triangular matrices with integer coefficients.
So the non-commutative $n$-tori are classified (up to $\T^n$-equivariant 
Morita equivalence)
 by 
$$
H^2(\Z^n,\UU)\cong M^u(n,\R)/M^u(n,\Z)\cong \T^{n(n-1)/2}.
$$
We refer to \cite{BK} for further details.
Moreover, for every action $\alpha$ of $\Z^n$ on $\K$ we can find an action 
$\beta$ of $\R^n$ 
on $\K$ such that the restriction $\beta|_{\Z^n}$ is Morita equivalent to 
$\alpha$ (choose an action 
$\beta$ corresponding to any class $[\eta]\in H^2(\R^n,\UU)$ 
which restricts to ${\rm Ma}(\alpha)$).

If $\omega\in H^2(G,\UU)$ for the abelian group $G$, then $\omega$ 
determines a continuous homomorphism
$h_\omega:G\to \widehat{G}$ given by 
\begin{equation}\label{eq-hom}
\langle h_\omega(g),h\rangle:=\omega(g,h)\omega(h,g)^{-1},\quad \text{for all } g,h\in G.
\end{equation}
$h_\omega$ only depends 
on the class $[\omega]\in H^2(G,\UU)$ and 
  $h_\omega= 0$ if  and only if $[\omega]=0$.
The kernel $S\subset G$ of $h_\omega$
is called the {\em symmetry group} of $\omega$. 
A cocycle $\omega$ is said to be {\em totally skew} 
if $S=\{0\}$, 
and $\omega$ is said to be {\em type I} if the image 
$h_\omega(G)$ is closed in $\widehat{G}$. 
Recall that for any C*-algebra $A$, $\Prim(A)$ denotes the space of primitive ideals 
of $A$ equipped with the Jacobson topology.
With these notations, the following results have been 
shown by Baggett and Kleppner in \cite[Section 3]{BK}. 

\begin{thm}\label{thm-hom}
For $\omega\in Z^2(G,\UU)$ the following are true:
\begin{enumerate}
\item There is a canonical bijection between 
$\widehat{S}$ and {${\rm Prim}(\C\rtimes_\omega G)$} given 
by induction of representations. 
In particular, $\C\rtimes_\omega G$ is simple if and only if $\omega$ is totally skew.
\item The image $h_\omega(G)$ is always a dense subgroup of 
$S^\perp\subset \widehat{G}$,
thus $\omega$ is type I if and only if $h_\omega(G)=S^\perp$.
\item The C*-algebra $\C\rtimes_\omega G$ is type I if and only if $\omega$ is type I.
\end{enumerate}
Moreover, the map 
\begin{eqnarray}
h_{\_}:H^2(G,\UU)\hookrightarrow {\rm Hom}(G,\widehat{G});\quad [\omega]\mapsto h_\omega
\label{EqTheInjectionOfH2}
\end{eqnarray}
is injective.
\end{thm}

Combining (1) and (3) of the above theorem we  see that $
\C\rtimes_\omega G$ 
is simple and type I if and only 
if $\omega$ is type I and totally skew,
i.e.
{$
h_\omega:G{\stackrel{\cong}{\to}}\widehat G
$
is an isomorphism.
 Since every separable, simple, 
type I C*-algebra is (isomorphic to) an algebra of compact operators on some 
separable Hilbert space,  we see that $\C\rtimes_\omega G\cong \K(\HH)$ 
for some separable Hilbert space $\HH$ if (and only if) $\omega$ 
is totally skew and type I. 
Then the dual action 
$
\widehat{G}\to \Aut(\C\rtimes_\omega G)=\Aut(\K(\HH))$ 
is again classified by a class
$[\widehat{\omega}]\in H^2(\widehat{G},\UU)$. 
This class has been computed by one of the authors in 
\cite[Lemma 3.3.5]{Echt96}:

\begin{prop}\label{prop-dualcocycle}
Suppose that $\omega\in Z^2(G,\UU)$ is a type I and totally skew 
2-cocycle. Then
the Mackey-obstruction for the dual action 
$\widehat{G}\to \Aut(\C\rtimes_\omega G)$ 
is given by the class of the cocycle 
$(h_\omega)_*\omega^{-1}\in Z^2(\widehat{G},\UU)$, 
where the push-forward is pullback along the inverse $h_\omega^{-1}$, 
i.e.
$$
(h_\omega)_*\omega^{-1}(\chi,\psi)={\omega(h_\omega^{-1}(\chi), h_\omega^{-1}(\psi))}^{-1},
\quad\chi,\psi\in \widehat{G}.
$$

\end{prop}
We call $(h_\omega)_*\omega^{-1}$ the {\it dual} cocycle of $\omega$,
and we leave the following lemma as an exercise for the reader.

\begin{lem}\label{lem-doubledual}
Suppose that $\omega\in Z^2(G,\UU)$ 
is type I and totally skew. 
Then $h_{(h_\omega)_*\omega^{-1}}= h_\omega^{-1}$ which implies that the
dual cocycle $(h_\omega)_*\omega^{-1}$ is also totally skew and type I. 
This also implies that the double dual cocycle  agrees with 
the original one:
$$
\big(h_{(h_\omega)_*\omega^{-1}}\big)_*\big((h_\omega)_*\omega^{-1}\big)^{-1}=\omega.
$$
\end{lem}

\begin{exa}\label{ex-totskew}
Let $G=\R^n$ and let us identify $\widehat{\R^n}$ with $\R^n$ via the canonical isomorphism 
$x\mapsto\chi_x$ with $\langle\chi_x,y\rangle=\exp(2\pi  i\   x^{t} y)$. 
Let  $\omega_A(x,y):=\exp({2\pi i (Ax)^{t} y})$ for some strictly 
upper diagonal matrix $A\in M(n,\R)$. 
Then 
$$\langle h_{\omega_A}(x),y\rangle=\exp({2\pi i( (Ax)^{t}y- (A y)^{t} x)})=
\exp({2\pi i (\Sigma_Ax)^{t} y})
=\langle \chi_{\Sigma_A x},y\rangle$$
with $\Sigma_A:=A-A^t$  the skew symmetric matrix corresponding to $A$. 
Thus we see that up to the identification $\widehat{\R^n}\cong \R^n$ the homomorphism 
$h_{\omega_A}$ is given by the linear map $x\mapsto \Sigma_A x$. 
It follows that $\omega_A$ is always type I and 
$\omega_A$ is totally skew if and only if $\Sigma_A$ is invertible. 
The dual cocycle is then given 
by $\omega_B$ with $B=\Sigma_A^{-1} A \Sigma_A^{-1}$.
\end{exa}

\section{Stable NC Tori -- NC T-Duality over the One-Point
Space}
\label{SecNCToverPoint}

\subsection{Introduction}\label{SecNCToverPoint1}
We start our discussion of T-duality with bundles over a point.
The trivial principal $\T^n$-bundle over the point is just the
$n$-torus $\T^n$ equipped with the translation action of $\T^n$ on itself. 
Suppose that  $\delta\in H^3(\T^n,\Z)$ allows an  action $\alpha$ of  $\R^n$ 
on the corresponding stable continuous-trace C*-algebra $CT(\T^n,\delta)$ 
which covers the inflated action  of $\R^n$ on $\T^n={\rm Prim}(CT(\T^n,\delta))$. 
Since $\T^n=\R^n/\Z^n$,  it follows from \cite[Theorem]{ech-ind} that
$CT(\T^n,\delta)$ is equivariantly isomorphic to the induced algebra 
${\rm Ind}_{\Z^n}^{\R^n}(\K,\tilde\alpha)$, where 
$\tilde\alpha$ denotes the action of the stabiliser $\Z^n$ 
on the fibre $\K=CT(\T^n,\delta)|_z$ of $CT(\T^n,\delta)$
over some {chosen  point} $z\in \T^n$.  
Recall that for any action $\beta:N\to \Aut(B)$ on a C*-algebra $B$, 
the induced system $(\Ind_N^G(B,\beta), G, \Ind(\beta))$ is given by 
\begin{align*}
\Ind_N^G(B,\beta):=\{F\in C_b(G,B): F(g+n)=&\beta_{-n}(F(g)), \quad\text{for all }g\in G, n\in N\\
& \;\text{and}\;\big(gN\mapsto\|f(g)\|\big)\in C_0(G/N)\},
\end{align*}
equipped with the point-wise operations, and with $G$-action $\Ind(\beta)$
which is just given by left (sign!) translation.
By the  discussion in the previous section, we may assume, up to equivariant Morita equivalence, that 
$\tilde\alpha=\mu|_{\Z^n}$ for some action $\mu:\R^n\to\Aut(\K)$.  But then 
we obtain an isomorphism
$$\Phi: \Ind_{\Z^n}^{\R^n}(\K,\tilde\alpha)\to C(\T^n, \K)\cong\K\otimes C(\T^n),\quad  \Phi(f)(\dot{g})={\mu}_g(f(g)),$$
which transforms the induced action $\Ind(\tilde\alpha)$ 
to the diagonal action\footnote{
To keep notation down, whenever there is a canonical 
action $\gamma$ of a quotient
$G/N$ (or $\widehat G/N^\perp$) 
we will simply denote by $\inf$ (rather than by $\inf\gamma$) 
the inflated action of $G$ (or $\widehat G$).}  
$\mu\otimes{\rm inf}$.
Thus we learn the following facts (which have been observed before by Mathai and Rosenberg in \cite{MR1, MR2}): 
Firstly, the class $\delta$ is trivial, i.e., $CT(\T^n,\delta)\cong \K\otimes C(\T^n)$.
Secondly, up to equivariant Morita equivalence the action $\alpha$ is given by a diagonal
action $\mu\otimes{\rm inf}$.

Recall from \cite{MR1, MR2} 
that in the above setting the (possibly non-commu\-tative) dual torus is given (again up to 
equivariant Morita equivalence)
by the crossed product $\K\otimes C(\T^n)\rtimes_{\mu\otimes{\rm inf}}\R^n$, 
equipped with the dual action of 
$\widehat{\R^n}\cong{\R^n}$.
 By the equivariant version of Green's imprimitivity theorem (this is a special case of
 \cite[Theorem 4.11]{EKQR}}), this system is Morita equivalent to 
$\K\rtimes_{\mu}\Z^n$ equipped with
the  action of  ${\R^n}$ 
which is inflated from the dual action of
$\T^n$ on $\K\rtimes_{\tilde\alpha}\Z^n$. 
 By the discussion in the previous section we know that 
 $\K\rtimes_{\mu}\Z^n$ is equivariantly  
 isomorphic to $\K\otimes (\C\rtimes_\omega \Z^n)$ equipped with the 
 inflated action
$\id \otimes {\inf}$ if ${\rm Ma}(\mu)=[\omega].$

To summarise the point-wise  duality picture of Mathai and Rosenberg,
there are stabilised commutative tori $\K\otimes C(\T^n)$ with a diagonal action 
$\mu\otimes\inf$ of $\R^n$ on one side, and there are stabilised non-commutative 
tori $\K\otimes (\C\rtimes_\omega\Z^n)$ 
with action $\id\otimes \inf$ on the other side. 
This motivates the content of this section which is the 
investigation of C*-dynamical systems
\begin{eqnarray}
\label{EqUrknall}
\Big(\K\otimes (\C\rtimes_{\omega} N), \widehat G,\hat \mu\otimes {\inf}\Big),
\end{eqnarray}
where we will typically assume that the 2-cocycle $\omega$ on $N$ 
has an extension to $G$, 
and $\hat\mu:\widehat G\to\Aut(\K)$ is an action which is not necessarily trivial.
By Takai duality the study of (\ref{EqUrknall}) is equivalent to the study of its
dual system
\begin{eqnarray}
\label{EqDualUrknall}
\Big(\big(\K\otimes (\C\rtimes_{\omega} N)\big)\rtimes_{\hat\mu\otimes\inf} 
\widehat G, G, {\widehat{\hat{\mu}\otimes\inf}}\Big),
\end{eqnarray}
and we need to understand under which circumstances there is a 
Morita equivalence
\begin{eqnarray*}
(\ref{EqDualUrknall})
\sim
\Big(\K\otimes (\C\rtimes_{\hat\omega} N^\perp), G,\mu\otimes {\inf}\Big).
\end{eqnarray*}

\subsection{Iterated Crossed Products and Transversality}
\label{SecItCrPrAndTrans}
For an action $\hat\mu:\widehat G\to \Aut(\K)$, 
let us analyse the iterated crossed product
\begin{eqnarray}
\label{EqDualUrknall2}
\big(\K\otimes (\C\rtimes_{\omega} N)\big)\rtimes_{\hat\mu\otimes\inf} \widehat G,
\end{eqnarray}
where we assume that the 2-cocycle $\omega$ has an extension to $G$ 
which we again denote by $\omega$. We need

{\begin{defi}\label{def-wedge0}
Suppose that $G$ is a locally compact abelian group.
\begin{enumerate}
\item 
 By the {\em Heisenberg cocycle} $\vee$ on 
${G}\times \widehat{G}$ we understand the 2-cocycle  given by
$$(h,\psi)\vee(g,\chi):=\lk \psi,g\rk.$$

\item
The Heisenberg cocycle on $\widehat G\times G$ 
is denoted by $\wedge$, i.e.
$$
(\psi,h)\wedge(\chi, g):=\lk h,\chi\rk.
$$

\item For $\omega\in Z^2(G,\UU)$ and $\hat\omega\in Z^2(\widehat{G},\UU)$ we denote by 
$\omega\vee\hat{\omega}$ the product $\omega\cdot\vee\cdot\hat{\omega}$ in which we regard
$\omega$ and $\hat\omega$ as cocycles on ${G}\times \widehat{G}$ by pullback along the projections to $G$ and $\widehat{G}$, respectively. 
Similarly, we define $\hat\omega\wedge\omega\in Z^2(\widehat G\times G,\UU)$.
\end{enumerate}

\end{defi}}

\begin{lem}\label{lem-wedge1} 
If $[\hat\omega]$ is the Mackey obstruction of $\hat\mu$, then 
(\ref{EqDualUrknall2}) is $G$-equivariantly isomorphic to 
\begin{eqnarray}
\label{EqDualUrknall3}
\K\otimes (\C\rtimes_{\omega\vee\hat\omega} ( N\times \widehat G)),
\end{eqnarray}
where $G$ acts dually on the second factor of $N\times \widehat G$.
\end{lem}
\begin{proof}
The product and the involution of  the crossed product (\ref{EqDualUrknall2})
 on the basis of the inflation action are given in terms of the pairing 
$\lk\_,\_\rk:N\times\widehat G\to \UU$ which on the level 
of $L^1$-functions can be expressed in terms of the cocycle
$\vee$.
In fact, similar to (\ref{EqIsoOfTwistedAndCrossed}), we
can define an isomorphism 
$$
\xymatrix{
\K\otimes L^1(N\times G,\omega\vee\hat\omega)
\ar[r]^-\cong_-\Psi
& 
L^1(N\times G,\K)
\cong
L^1(G, \K\otimes L^1(N,\omega))\\
\ar@{}[u]|{\resizebox{!}{0.4cm}{$\cap$}}
\K\otimes{\big(\C\rtimes}_{\omega\vee\hat\omega}(N\times\widehat G)\big)
&
\ar@{}[u]|{\resizebox{!}{0.4cm}{$\cap$}}
L^1(G, \K\otimes (\C\rtimes_\omega N) )
\\
&
\ar@{}[u]|{\resizebox{!}{0.4cm}{$\cap$}}
\big(\K\otimes (\C\rtimes_\omega N)\big)\rtimes_{\hat\mu\otimes\inf}\widehat G 
}
$$
by 
$\Psi(k\otimes f)(n,\chi)= f(n,\chi)k\hat V(\chi)$
where $\hat V:\widehat{G}\to \U(\HH)$ is a Borel map  with  $\hat\mu=\Ad \circ \hat V$ and such that $\hat\omega^{-1}=\partial \hat V$.
It is  straightforward to check that this isomorphism is $G$-equivariant with respect to the dual actions.
\end{proof}

Applying Proposition \ref{prop-dualmor} to (\ref{EqDualUrknall3}) 
for the subgroup $N\times \widehat G\subset G\times \widehat G$
we obtain the Morita equivalent system 
\begin{eqnarray}
\label{EqDualUrknall4}
\K\otimes (\C\rtimes_{\omega\vee\hat\omega} 
(  G\times \widehat G))\rtimes_{\rm dual} (N^\perp\times0),
\end{eqnarray}
whereon $G$ acts by the dual action on the second group factor of the inner crossed product.
We want to apply Proposition \ref{prop-dualcocycle} to this inner 
crossed product, so we are interested in the properties of 
$$
h_{\omega\vee\hat\omega}:G\times \widehat G\to 
\widehat G\times G
$$
defined in (\ref{eq-hom}). 

\begin{lem}\label{lem-homega-Hesisenberg}
The homomorphism $
h_{\omega\vee\hat\omega}:G\times \widehat G\to \widehat G\times G
$ is given by 
\begin{eqnarray*}
h_{\omega\vee\hat\omega}=\left( 
\begin{matrix}
h_{\omega}&{\rm id}_{\widehat G} \\
-{\rm id}_{ G}& h_{\hat\omega}
\end{matrix}
\right): (g,\chi)\mapsto (h_\omega(g)+\chi,h_{\hat\omega}(\chi)-g).
\end{eqnarray*}
\end{lem}

\begin{proof} 
We have
\begin{eqnarray*}
\lk h_{\omega\vee\hat\omega}(g,\chi),(h,\psi)\rk
&=&
\omega\vee\hat\omega((g,\chi),(h,\psi))
\ 
\omega\vee\hat\omega((h,\psi),(g,\chi))^{-1}\\
&=&
\lk h_\omega(g),h\rk
\lk h_{\hat\omega}(\chi),\psi\rk
\lk \chi,h\rk \lk\psi,g\rk^{-1}\\
&=&
\lk h_\omega(g)+\chi,h\rk
\lk h_{\hat\omega}(\chi)-g,\psi\rk.
\end{eqnarray*}
The canonical identification $G\cong\widehat{\!\widehat G}$ yields the result.
\end{proof}

\begin{lem}\label{LemThreeCondOnPhi}
The following three conditions are equivalent:
\begin{enumerate}
\item $\phi:= {\rm id}_G+h_{\hat\omega} \circ h_\omega :G\to G$
 is an isomorphism.
\item $\hat\phi:= 
{\rm id}_{\widehat G}+h_\omega \circ h_{\hat\omega} :
\widehat G\to \widehat G$ is an isomorphism.
\item $h_{\omega\vee\hat\omega}$ is an isomorphism.
\end{enumerate}
\end{lem}

\begin{proof}
The equivalence of the first two statements follows from the observation that
$\hat\phi$  is the dual of $\phi$.
Alternatively,
if $\phi^{-1}$ exists, 
then a one-line calculation shows that
 ${\rm id}_{\widehat G}- h_\omega\circ\phi^{-1} \circ h_{\hat \omega}$ is an 
inverse for $\hat\phi$.

Now, if  $h_{\omega\vee\hat\omega}$ is an isomorphism, then
$$
h_{\omega\vee\hat\omega}^\dag:=
{\rm flip}\circ h_{\omega\vee\hat\omega}\circ{\rm flip}
=
\left( 
\begin{matrix}
h_{\hat\omega}&-{\rm id}_{ G} \\
{\rm id}_{\widehat G}& h_{\omega}
\end{matrix}
\right) 
$$
is an isomorphism, 
where ${\rm flip}:\widehat G\times  G\to  G\times \widehat G$
is transposition. Then
the composed ismorphism is
$$
h_{\omega\vee\hat\omega}\circ h_{\omega\vee\hat\omega}^\dag 
= \left( 
\begin{matrix}
{\rm id}_{\widehat G}+h_{\omega} \circ h_{\hat\omega}  & 0 \\
0 & {\rm id}_{ G}+h_{\hat \omega} \circ h_{\omega}
\end{matrix}
\right)
= \left( 
\begin{matrix}
\hat \phi  & 0 \\
0 & \phi
\end{matrix}
\right).
$$
Conversely, if $\phi, \hat\phi$ 
are isomorphisms then
$$
h_{\omega\vee\hat\omega}^{-1}:=
h_{\omega\vee\hat\omega}^{\dag}
\circ
\left( 
\begin{matrix}
\hat\phi^{-1}  & 0 \\
0 &\phi^{-1}
\end{matrix}
\right)
$$
exists and is obviously a right inverse.  The property of being
 a left inverse 
requires a small calculation. First we have
\begin{eqnarray*}
h^{-1}_{\omega\vee\hat\omega}\circ h_{\omega\vee\hat\omega}
&=&
\left( 
\begin{matrix}
h_{\hat\omega}\circ\hat\phi^{-1}\circ h_\omega+\phi^{-1}  
&  
h_\omega\circ \phi^{-1}-\hat\phi^{-1}\circ h_{\hat\omega}
\\
\hat\phi^{-1}\circ h_{\omega} -h_{\omega}\circ \phi^{-1} 
&
\hat\phi^{-1}+h_{\omega}\circ\phi^{-1}\circ h_{\hat\omega}
\end{matrix}
\right).
\end{eqnarray*}
Inserting now $\hat\phi^{-1}=
{\rm id}_{\widehat G} - h_\omega\circ\phi^{-1}\circ h_{\hat\omega}$
into this, we derive at
\begin{eqnarray*}
h^{-1}_{\omega\vee\hat\omega}\circ h_{\omega\vee\hat\omega}
&=&
\left( 
\begin{matrix}
\id_G&0
\\
0& \id_{\widehat G}
\end{matrix}
\right).
\end{eqnarray*} 
\end{proof}

Let's assume one of the equivalent conditions of Lemma \ref{LemThreeCondOnPhi}.
Then the inner crossed product in (\ref{EqDualUrknall4}) is isomorphic to
the compacts. 
Having determined $h_{\omega\vee\hat\omega}^{-1}$ we can compute
the dual cocycle
$ (h_{\omega\vee\hat\omega})_*({\omega\vee\hat\omega})^{-1}
 $ according to  
Proposition \ref{prop-dualcocycle}. The result helps us to
understand the remaining outer crossed product  with $N^\perp$
in (\ref{EqDualUrknall4})
and also the remaining $G$-action on it.
To state the result,
some more notation is useful.
\begin{defi}
\label{DefiOfTheDualMoped}
For a 2-cocycle $\omega$ on $G$
and a 2-cocycle 
$\hat\omega$ on $\widehat G$, we define new 
cocycles on $G$ and $\widehat G$ by
$$
\omega{\rtimes}\hat\omega :=\omega\cdot h_{\omega}^*\hat\omega^{-1},
\qquad
\hat\omega{\rtimes}\omega :=\hat\omega\cdot h_{\hat\omega}^*\omega^{-1}.
$$
Moreover, if $\phi=\id_G+h_{\hat\omega}\circ h_\omega$ is an isomorphism, 
we define
$$
\omega \bar\rtimes \hat\omega := \phi_*(\omega{\rtimes}\hat\omega),
\qquad
\hat \omega \bar\rtimes \omega := \hat\phi_*(\hat\omega{\rtimes}\omega).
$$
\end{defi}

\begin{lem}\label{LemComputationOfWidehatOmega}
In $H^2(\widehat G\times G,\UU)$ the following 
equality holds\footnote{
We denote the group operation on $U(1)$-valued 
2-cocycles multiplicatively, whereas we denote 
the group operation on any cohomology group additively.}:
\begin{eqnarray}
\label{EqTheDualCocycleSpelledOut}
[ (h_{\omega\vee\hat\omega})_*({\omega\vee\hat\omega})^{-1} ]
&=&
[\hat\omega\bar\rtimes\omega] 
{+}
[\omega\bar\rtimes\hat \omega]
{-}
(\hat\phi\times\id_G)_* [\wedge],
\end{eqnarray}
where the classes on $\widehat G$ and $G$
are understood as classes on $\widehat G\times G$
by pullback along the projections. 
\end{lem}

\begin{proof}
Before we start with the actual computation 
we need to be aware of some general cocycle properties.
Let $\nu$ be a 2-cocycle on any abelian group.
Twofold application of the cocycle identity gives
\begin{eqnarray*}
\nu(x+y,-(x+y)) &=& 
\nu(y,-x-y)\ \nu(x,-x)^{-1}\ \nu(x,y)^{-1}\\
&=&
\nu(0,x)\ \nu(-x,-y)^{-1}\ \nu(y,-y)\ \nu(x,-x)^{-1}\ \nu(x,y)^{-1}
\end{eqnarray*}
which means that $c(x):=\nu(x,-x)\nu(0,0)$
is a cochain that implements
$\nu(x,y)\sim \nu(-y,-x)^{-1}$.

Furthermore, a fourfold application of the 
cocycle identity gives 
\begin{eqnarray*}
\nu(a+x,b+y)&=&
\nu(a,b)\ \nu(x,y)\nonumber\\
&&\cdot\, \nu(x,b)\ \nu(b,x)^{-1}\\
&&
\cdot\,\nu(a,x)^{-1}\nu(b,y)^{-1}\nu(a+b,x+y).
\end{eqnarray*}
Or, if we define a 2-cocycle $\tilde\nu$ on the product of the group with itself
by  $\tilde\nu\big( (a,x),(b,y)\big):=\nu(a+x,b+y)^{-1},$
then
\begin{eqnarray*}
\tilde\nu\big( (a,x),(b,y)\big)&=&
\nu(a,b)^{-1}\
\nu(x,y)^{-1}\
\langle h_\nu(x),b\rangle^{-1}\
(d\nu)\big( (a,x),(b,y)\big)
\\
&\sim&
\nu(a,b)^{-1}\
\nu(x,y)^{-1}\
\langle x,h_\nu(b)\rangle,
\end{eqnarray*}
where  $d$ is the boundary operator on the 
product of the group with itself.

Let us now turn to the actual computation.
By definition we have to compute
$
(h_{\omega\vee\hat\omega}^{-1})^{*}({\omega\vee\hat\omega})^{-1},
$ 
for 
$$
h_{\omega\vee\hat\omega}^{-1}=
\left( 
\begin{matrix}
h_{\hat\omega}&-{\rm id}_{ G} \\
{\rm id}_{\widehat G}& h_{\omega}
\end{matrix}
\right)
\circ
\left( 
\begin{matrix}
\hat\phi^{-1}&0 \\
0&\phi^{-1}
\end{matrix}
\right).
$$
Let us use the shorthands $\chi':=\hat\phi^{-1}\chi$ 
and $g':=\phi^{-1}g$. Then we have
 \begin{eqnarray}
&&(h_{\omega\vee\hat\omega})_*({\omega\vee\hat\omega})^{-1}
\big( (\chi,g),(\psi,h)\big)
=\nonumber\\
&&\hspace{2cm}\omega(h_{\hat\omega}(\chi')-g',h_{\hat\omega}(\psi')-h')^{-1}
\label{OmegaMitVielDrin}\\
&&\hspace{3cm}\cdot\,\langle \chi'+h_{\omega}(g'),
h_{\hat\omega}(\psi')-h'\rangle^{-1}
\label{MitVielDrin}\\
&&\hspace{4cm}\cdot\
\hat\omega(\chi'+h_\omega(g'),\psi'+h_\omega(h'))^{-1}
\label{HatOmegaMitVielDrin}
\end{eqnarray}
The middle term (\ref{MitVielDrin}) decomposes
to
\begin{eqnarray*}
(\ref{MitVielDrin})
&=&
\langle\chi', h'\rangle\
\langle h_\omega(g'), h_{\hat\omega}(\psi')\rangle^{-1}\
\hat\omega(\chi', \psi')\
\hat\omega(\psi',\chi')^{-1}\
\omega(g',h')\
\omega(h',g')^{-1}\\
&\sim&
\langle\chi', h'\rangle\
\langle h_\omega(g'), h_{\hat\omega}(\psi')\rangle^{-1}\
\hat\omega(\chi', \psi')\
\hat\omega(-\chi',-\psi')\
\omega(g',h')\
\omega(h',g')^{-1}
\end{eqnarray*}
But as $((\chi,g),(\psi,h))\mapsto \lk \chi,h\rk$ 
is cohomologous to 
$((\chi,g),(\psi,h))\mapsto \lk \psi,g\rk^{-1}$ 
(just by the cochain $(\chi,g)\mapsto \lk g,\chi\rk)$
this can be transformed to
\begin{eqnarray*}
(\ref{MitVielDrin})
&\sim&
\langle g', \psi'\rangle^{-1}\
\langle h_\omega(g'), h_{\hat\omega}(\psi')\rangle^{-1}\
\hat\omega(\chi', \psi')\
\hat\omega(\psi',\chi')^{-1}\
\omega(g',h')\
\omega(h',g')^{-1}.
\end{eqnarray*}
To transform 
(\ref{OmegaMitVielDrin})
and (\ref{HatOmegaMitVielDrin}), 
we apply the above identity for $\tilde\nu$.
We find 
\begin{eqnarray*}
(\ref{HatOmegaMitVielDrin})
&\sim&
\hat\omega(\chi',\psi')^{-1}\
\hat\omega(h_\omega(g'),h_\omega(h'))^{-1}\
\langle  (h_{\omega}(g'), h_{\hat\omega}(\psi')\rangle
\end{eqnarray*}
and
\begin{eqnarray*}
(\ref{OmegaMitVielDrin})
&\sim&
\omega(h_{\hat\omega}(\chi'),h_{\hat\omega}(\psi'))^{-1}\
\omega(-g',-h')^{-1}\
\langle -g', h_\omega(h_{\hat\omega}(\psi'))\rangle\\
&\sim&
\omega(h_{\hat\omega}(\chi'),h_{\hat\omega}(\psi'))^{-1}\
\omega(h',g')\
\langle -g', h_\omega(h_{\hat\omega}(\psi'))\rangle.
\end{eqnarray*}
Multiplying these partial results
we get
\begin{eqnarray}
\label{EqTheLastUnicorn}
(\ref{OmegaMitVielDrin})\cdot
(\ref{MitVielDrin})\cdot
(\ref{HatOmegaMitVielDrin})
&\sim&
\lk g',\psi' +h_{\omega}(h_{\hat\omega}(\psi'))\rk^{-1}\\
&&\cdot
\omega(g',h')\
h_\omega^*\hat\omega(g',h')^{-1}
\nonumber\\
\nonumber
&&\cdot
\hat\omega(-\chi',-\psi')\
h_{\hat\omega}^*\omega(\chi',\psi')^{-1}
\end{eqnarray}
which is the claimed formula up to the minus sign inside the argument 
of $\hat\omega$.
However, the injectivity of the map 
$[\hat\omega]\mapsto h_{\hat\omega}\in {\rm Hom}(\widehat G,G)$
implies that $\hat\omega(\chi,\psi)\sim\hat\omega(-\chi,-\psi)$.
So the lemma  is proven.
\end{proof}

\begin{rem}
\label{RemImportantNotice}
It is very important for later purposes (Theorem \ref{THEMAINTHEOREM}) to observe 
at this point that the computation done in the proof of Lemma
\ref{LemComputationOfWidehatOmega}
is based on explicit cochains 
(given in terms of 
$\hat\omega,\omega$ or $\lk\_,\_\rk$)
up to and including
(\ref{EqTheLastUnicorn}).
Only the very last step 
$\hat\omega(\chi,\psi)\sim\hat\omega(-\chi,-\psi)$
required an abstract argument.
If $\hat\omega$ is cohomologous to a  
bicharacter then this last relation 
can also be made explicit:
Let 
$\hat\omega=dc\ \hat\eta$ for some 
bicharacter $\hat\eta$, then
\begin{eqnarray*}
\hat\omega(\chi,\psi)
&=&dc(\chi,\psi)\ \hat\eta(\chi,\psi)\\
&=&dc(\chi,\psi)\ \hat\eta(-\chi,-\psi)\\
&=&dc(\chi,\psi)\ dc(-\chi,-\psi)^{-1}\ dc(-\chi,-\psi)\hat\eta(-\chi,-\psi)\\
&=&d\tilde c(\chi,\psi)\ \hat\omega(-\chi,-\psi),
\end{eqnarray*}
for $\tilde c(\chi):= c(\chi)c(-\chi)^{-1}$.
These explicit cohomology relations are abstract 
cocycle identities and do not depend on $\UU$ as a module.
In fact, Lemma \ref{LemComputationOfWidehatOmega}
remains valid 
if (1) the involved cocycles are 
 not just $\UU$-valued but 
$C(B,\UU)$-valued (trivial module structure)
for some space $B$, if (2) they are 
cohomologous to bihomomorphisms
(rather than bicharacters),
and if (3) we have control over the 
continuity properties of the 
quantities $h_\omega,h_{\hat\omega}$
which then should be regarded
as bundle maps 
$$
\xymatrix{
B\times G\ar[r]^{h_\omega} \ar[d]&
B\times \widehat G \ar[r]^{h_{\hat\omega}}\ar[d]&
B\times G\ar[d]\\
B\ar@{=}[r]&B\ar@{=}[r]&B.
}
$$
This will ensure that if a cocycle with values in $C(B,\UU))$ is point-wise pulled back by 
$h_\omega,h_{\hat\omega}$ or 
pushed foreward by $\phi$ (which involves inversion in ${\rm Aut}(G)$!), 
then the resulting point-wise defined object is again mapping to $C(B,\UU)$.
Equality (\ref{EqTheDualCocycleSpelledOut})
then  holds in 
$H^2(\widehat G\times G, C(B,\UU))$.
\end{rem}

The structure of 
$[
 (h_{\omega\vee\hat\omega})_*({\omega\vee\hat\omega})^{-1}
]$
given by its three summands now immediately yields the following 
corollary. 

\begin{cor}
\label{CorTheDualOverThePoint}
If $\phi$ is an isomorphism, the crossed product (\ref{EqDualUrknall4}) 
is $G$-equivariantly isomorphic to 
\begin{eqnarray}
\label{EqDualUrknall5}
\K\otimes\K\otimes
(\C\rtimes_{\hat\omega\bar\rtimes\omega}N^\perp),
\end{eqnarray}
where $G$ acts by $\id\otimes\mu\otimes {\rm in^\phi}$ for an action $\mu$ with 
Mackey obstruction ${\rm Ma}(\mu)=[\omega\bar\rtimes\omega]$
in 
$H^2(G,\UU)$, and 
$({\rm in^\phi}_g(f))(n):=\lk g,\hat\phi^{-1}(n)\rk^{-1} f(n)$.
\end{cor}

So except from the part of the action given by $\rm in^\phi$
we have found a structure rather similar  to the one with which we have started.
To manipulate it a little further we need an extra assumption.

\begin{defi}
\label{DefiTransvers}

\begin{enumerate}
\item
A homomorphism $G\to G$ is said to be an {\em automorphism of} 
$(N,G)$ if it is an automorphism of $G$ and if it 
maps $N$ bijectively to itself.
We denote by $\Aut(N,G)$ the set of all of those.

\item
A pair of cocycles $\omega:G\times G\to\UU$, 
$\hat\omega:\widehat G\times\widehat G\to\UU$
(or their cohomology classes)
is called {\em transverse} if 
$\phi=\id_G+h_{\hat\omega}\circ h_\omega{\in\Aut(N,G)}$.
(Equivalently, one might require that $\hat\phi=\id_{\widehat G}+h_{\omega}\circ h_{\hat\omega}$ is an 
automorphism of $(N^\perp,\widehat G)$.)
\item
The binary relation defined by transversality is denoted by 
$$\tra\subset H^2(G,\UU)\times H^2(\widehat G,\UU),
$$
i.e. $\omega\tra\hat\omega$ if and only if $\omega$ and $\hat\omega$ are transverse.
\item
Let  $\hat\mu\otimes \inf$ be a 
$\widehat G$-action
on $\K\otimes (\C\rtimes_\eta N)$.
The actions $\hat\mu$ or $\hat\mu\otimes\inf$ or the dynamical 
system $(\K\otimes (\C\rtimes_\eta N),\widehat G,\hat\mu\otimes\inf)$
are called {\em transverse} if the 
cocycle $\eta:N\times N\to\UU$
has an extension $\omega$ to $G$
such that $\omega\tra{\rm Ma}(\hat\mu)$.
\end{enumerate}
\end{defi}

\begin{cor}
\label{CorSuffCond}
If the dynamical system $\Big(\K\otimes (\C\rtimes_{\omega} N), \widehat G,\hat \mu\otimes {\inf}\Big)$ 
is transverse, then  its dual system
$$
\Big(\big(\K\otimes (\C\rtimes_{\omega} N)\big)\rtimes_{\hat\mu\otimes\inf} 
\widehat G, G, \widehat{\hat\mu\otimes\inf}\Big) \cong 
\Big(\K\otimes\K\otimes
(\C\rtimes_{\hat\omega\bar\rtimes\omega}N^\perp), G, \id\otimes\mu\otimes {\rm in^\phi}\Big)
$$
is $G$-equivariantly isomorphic 
to 
$$
\K\otimes (\C\rtimes_{\hat\omega\rtimes\omega}N^\perp),
$$
where $G$ acts by $\mu\otimes \inf$ with
Mackey obstruction ${\rm Ma}(\mu)=[\omega\bar\rtimes\hat\omega]$ 
$(=[\phi_*({\omega\rtimes\hat\omega})])$.
\end{cor}

\begin{proof}
Using $\K\otimes\K\cong \K$ and Corollary \ref{CorTheDualOverThePoint} it suffices to show that 
$$\Big(\K\otimes
(\C\rtimes_{\hat\omega\bar\rtimes\omega}N^\perp), G, \mu\otimes {\rm in^\phi}\Big)
\cong
 \Big(\K\otimes (\C\rtimes_{\hat\omega\rtimes\omega}N^\perp), G, \mu\otimes \inf\Big).$$
By transversality, $\hat\phi$ induces an isomorphism of $N^\perp$,
so 
it induces an isomorphism 
\begin{eqnarray}
\label{EqTheSimplifiedDual}
\hat\phi^\star:\C\rtimes_{\phi_*(\hat\omega\rtimes\omega)} N^\perp
\cong
\C\rtimes_{\hat\omega\rtimes\omega} N^\perp
\end{eqnarray}
given by pullback: $f\mapsto f\circ\phi$. 
Similarly, the inversion on the group
 $\ominus: N^\perp\to N^\perp$
induces an automorphism by pullback
\begin{eqnarray}
\ominus^\star:\C\rtimes_{\hat\omega\rtimes\omega} N^\perp
\cong
\C\rtimes_{\ominus^*(\hat\omega\rtimes\omega)} N^\perp.
\end{eqnarray}
Note 
that 
$$
({\rm in^\phi}_g(f))(\phi(-n))=\lk g,-n\rk^{-1} f(\phi(-n))=
\lk g,n\rk f(\phi(-n)),
$$
so the composition $\ominus^\star\circ\phi^\star$ 
turns $\rm in^\phi$ into the ordinary inflation action.
However, pullback of a 2-cocycle along the
inversion gives a cocycle that is similar to the original one, i.e. they
have the same cohomology class
(see Remark \ref{RemImportantNotice}). Then  their twisted group algebras
are equivariantly isomorphic.
\end{proof}

By Takai duality, we know that the dual of the constructed system 
$$(\K\otimes (\C\rtimes_{\hat\omega\rtimes\omega}N^\perp),\mu\otimes\inf)$$
(i.e. the bidual of the original system) is Morita equivalent to
the original (transverse) system. 
The following lemma shows that the dual system of a transverse system 
is  transverse again.

\begin{lem}
\label{LemDualityOnTransversalitySet}
The assignment $(\omega,\hat\omega)\mapsto (\omega\bar\rtimes\hat\omega,\hat\omega\rtimes\omega)$
defines a bijection
$\tra\to \tra$
such that
$$
\xymatrix{
\tra\ar[d]\ar[r]^\cong&\tra\ar[d]\\
\Aut(N,G)\ar[r]^{\id}&\Aut(N,G)
}
$$
commutes, where the vertical  arrows are given by the tautological map.
\end{lem}

\begin{proof}
Firstly, $\phi':=\id+h_{\hat\omega\rtimes\omega}\circ h_{\omega\bar\rtimes\hat\omega}$ is an isomorphism 
of $(N,G)$:
A one-line computation gives
$
h_{h_{\hat\omega}^*\omega^{-1}}=h_{\hat\omega}\circ h_\omega\circ h_{\hat\omega}: \widehat G\to\widehat G,
$
and so
$
h_{\hat\omega\rtimes\omega} = h_{\hat\omega\cdot h_{\hat\omega}^*\omega^{-1}}
=
h_{\hat\omega}+h_{\hat\omega}\circ h_\omega\circ h_{\hat\omega}
= h_{\hat\omega}\circ \hat \phi, 
$
wherein as before $\hat\phi = \id+ h_{\omega}\circ h_{\hat\omega}$.
The same algebra gives
$h_{\omega\bar\rtimes\hat\omega}=h_{\phi_*(\omega\rtimes\hat\omega)}
=\hat\phi^{-1}\circ h_{\omega\rtimes\hat\omega}\circ \phi^{-1}
=\hat\phi^{-1}\circ h_{\omega}$. 
Therefore
$$
\phi'=\id+h_{\hat\omega\rtimes\omega}\circ h_{\omega\bar\rtimes\hat\omega}=\id+(h_{\hat\omega}\circ\hat\phi)\circ
(\hat\phi^{-1}\circ h_{\omega})=\phi
$$
which is an automorphism of $(N,G)$ by assumption.
This also shows that the diagram of the lemma commutes.

Secondly, another straight forward calculation shows that the inverse of 
$(\omega,\hat\omega)\mapsto (\omega\bar\rtimes\hat\omega,\hat\omega\rtimes\omega)$
is given  by
$$
(\omega,\hat\omega)\mapsto
(\omega\rtimes \hat\omega,\hat\omega\bar\rtimes\omega).
$$
\end{proof}

Corollary \ref{CorSuffCond} tells that transversality 
gives a sufficient condition to answer the
question raised in (\ref{EqDualUrknall}).
Combining it together with Lemma \ref{LemDualityOnTransversalitySet}
we have found a class of C$^*$-dynamical systems
which is closed under taking crossed products:
Let us
denote by $NCT(N;\widehat G)$ 
the 2-category of systems
$(\K\otimes(\C\rtimes_\omega N),\widehat G,\hat\mu\otimes\inf),$ 
which has $\widehat G$-equivariant 
Morita equivalences as 1-morphisms
and equivariant isomorphisms between them as 2-morphisms.
There is a proper subcategory  
$NCT(N;\widehat G)^\tri\subset NCT(N;\widehat G)$
which consists of systems which are 2-isomorphic  
(i.e. Morita equivalent) to a transverse representative.
This whole section is summarised in

\begin{thm}
\label{ThmDualityOverThePoint}
The duality functor ${\_\rtimes G}$
defined on all C$^*$-dynamical systems with group $G$
restricts to a duality of transverse dynamical systems:
$$
\xymatrix{
{\begin{array}{c}
\textrm{C$^*$-Dynamical Systems}\\
\textrm{with {G}roup $G$}
\end{array}}
\ar[r]^-\sim
&
{\begin{array}{c}
\text{C$^*$-Dynamical Systems}\\
\text{with {G}roup $\widehat G$}
\end{array}}
\\
\ar@{}[u]|{{\resizebox{0.4cm}{!}{$\cup$} }}
NCT(N^\perp;G) &
\ar@{}[u]|{{\resizebox{0.4cm}{!}{$\cup$} }}
NCT(N;\widehat G)
\\
\ar@{}[u]|{{\resizebox{0.4cm}{!}{$\cup$} }}
NCT(N^\perp;G)^\tri \ar[r]^-\sim&
\ar@{}[u]|{{\resizebox{0.4cm}{!}{$\cup$} }}
NCT(N;\widehat G)^\tri
}.
$$
\end{thm}

\subsection{Classification Remarks}
\label{SecClassifiRemarks}
Recall from section \ref{SecZeugsUeberK}
that the C*-dynamical systems $\K\rtimes_\alpha N$ or $\C\rtimes_\omega N$
equipped with their canonical $\widehat N$-actions are classified up to equivariant 
Morita equivalence by second Borel cohomology
$$
{\rm Ma}(\alpha),[\omega]\in H^2(N,\UU).
$$
We start with an example that illustrates 
that the objects with which we are dealing are more involved.

\begin{exa}
\label{ExaTwoDimAbiguity} 
Let $G=\R^2$ and $N=\Z^2$, and choose 
${\frac{1}{3}},{\frac{2}{3}}\in\R/\Z=\T\cong H^2(N,\UU)$. 
Denote by $\hat\mu_3$ an action of  $\widehat G=\R^2$ on $\K$ with
Mackey obstruction $3\in \R\cong H^2(\R^2,\UU)$.
Then there is a $\widehat G$-equivariant Morita equivalence
$$
\Big(\K\otimes (\C\rtimes_{\frac{2}{3}}N), \hat\mu_3\otimes {\rm inf} \Big)
\sim
\Big(\K\otimes (\C\rtimes_{\frac{1}{3}}N), {\rm id}\otimes {\rm inf} \Big).
$$
\end{exa}

\begin{proof}
A lengthy but direct proof is given in Appendix \ref{AppSiggisExample}.
Using our theory of transversality one can significantly shorten the proof. 
This is done in section \ref{Sec2DimNCT} below.
\end{proof}

To understand why this example could possibly be true
let us try to understand what can be said in general about 
an equality of classes 
$[\omega_1],[\omega_2]\in H^2(N,\UU)$ if there is a
$\widehat G$-equivariant Morita equivalence 
\begin{equation}\label{EqMoritaEq}
\K\otimes (\C\rtimes_{\omega_1} N)\sim \K\otimes (\C\rtimes_{\omega_2} N), 
\end{equation}
where $\widehat G$ acts on 
both sides diagonally by actions on the compacts 
$\hat\mu_i:\widehat G\to{\rm Aut}(\K)$, $i=1,2$, tensor
the inflated actions 
${\rm inf}_i:\widehat G\to {\rm Aut}(\C\rtimes_{\omega_i} N),i=1,2,$ 
on the respective twisted group C*-algebras.
Denote by $S\subset N$ the symmetry group of $\omega_2$ that 
is the kernel of the map $h_{\omega_2}:N\to \widehat N$.
The dual group $\widehat S$ of $S$ is homeomorphic to
the primitive spectrum of $\K\otimes (\C\rtimes_{\omega_2} N)$, 
and, by the Daums-Hofmann Theorem (see \cite{Wbook}),
 the $\UU$-valued functions thereon are isomorphic
to the center of the unitary group of its multiplier algebra, 
i.e. there is a short exact sequence 
\begin{eqnarray*}
1\to C(\widehat S,\UU)\to {\rm UM}(\K\otimes(\C\rtimes_{\omega_2} N))
\stackrel{\rm Ad}{\to}
{\rm Inn}(\K\otimes (\C\rtimes_{\omega_2} N))\to 1,
\end{eqnarray*}
and this sequence is $\widehat N$-equivariant for the actions that
are induced by the dual action on $\K\otimes (\C\rtimes_{\omega_2} N)$.
On $C(\widehat S,\UU)$ this action is just translation in the argument
after restriction $\widehat N\to \widehat S$. Using the polish topology 
on ${\rm Inn}(A)$ induced 
from the polish strict topology of ${\rm UM}(A)$ together with 
\cite[Corollary 0.2]{RR}, this short exact sequence 
induces a (not very long) exact sequence in Borel cohomology
\begin{eqnarray*}
\xymatrix{
\cdots\ar[r] &
H^1(\widehat N,C(\widehat S,\UU))
\ar[r]
&
H^1(\widehat N,  {\rm UM}(\K\otimes (\C\rtimes_{\omega_2} N)))
\ar 
`[r] 
`[l] 
`[lld]
`[r]
[dl]
&\\
&
H^1(\widehat N,{\rm Inn}(\K\otimes (\C\rtimes_{\omega_2} N)))
\ar[r]^-\delta
&
H^2(\widehat N, C(\widehat S,\UU))
}
\end{eqnarray*}
which terminates at $H^2(\widehat N, C(\widehat S,\UU))$ due to the 
non-commutativity of the involved coefficient groups.

Now, let us tensor both sides of (\ref{EqMoritaEq}) with another copy 
of the compacts which we then equip with an action 
$\hat\mu_1^{op}$ which has the inverse Mackey obstruction 
of $\hat\mu_1$:  ${\rm Ma}(\hat\mu_1^{op})=-{\rm Ma}(\hat\mu_1)$. 
As
$\hat\mu_1^{op}\otimes\hat \mu_1$ is 
Morita equivalent to the trivial action,
we have to deal with the Morita equivalent 
actions
$\id\otimes\id\otimes \inf_1$  and 
$\hat\mu^{op}_1\otimes\hat\mu_2\otimes\inf_2$.
Because being Morita equivalent is the same as being (stably) outer conjugate,
we conclude that the action on the left hand side 
$(\id\otimes {\rm id}\otimes {\rm inf}_1)$ 
is conjugate to 
\begin{eqnarray*}
{\rm Ad}(u)\circ(\hat\mu_1^{op}\otimes\hat\mu_2\otimes {\rm inf_2})&=&
\underbrace{{\rm Ad}(u)\circ(\hat\mu_1^{op}\otimes\hat\mu_2\otimes{\rm id})}
\circ ({\rm id}\otimes\id\otimes {\rm inf_2})\\
&&\qquad\qquad\quad =:\gamma,
\end{eqnarray*}
for some continuous 1-cocycle 
$u:\widehat G\to {\rm UM}(\K\otimes\K\otimes(\C\rtimes_{\omega_2} N))$.
It follows right from the definition of $\gamma$ that it vanishes 
on $N^\perp$. Indeed, we have the following:

\begin{lem}
The above defined map 
$\gamma$ factors over the quotient map
$$
\xymatrix{
\widehat G\ar[r]^-\gamma\ar@{->>}[d]& 
 {\rm Inn}(\K\otimes\K\otimes(\C\rtimes_{\omega_2} N))\\ 
\widehat N\ar[ur]_-{\dot\gamma}
}
$$
and $\dot\gamma$ satisfies the cocycle relation
$$
\dot\gamma(\hat n+\hat m)=
\dot\gamma(\hat n)\circ \big((\hat n\cdot\dot\gamma(\hat m)\big),\quad \hat n,\hat m\in \widehat N,
$$
where $\cdot$ is precisely the action of 
$\widehat N$ on the inner automorphisms 
that occurred in the short exact sequence above.
\end{lem}
\begin{proof}
Let's denote the isomorphism which implements
the stable outer conjugacy between the two actions by 
$\Phi$, i.e.
$$
\Phi\circ ({\rm id}\otimes\id\otimes {\rm inf}_1)\circ\Phi^{-1} 
=
\gamma \circ ({\rm id}\otimes\id\otimes {\rm inf}_2),
$$ 
and let
$\lk\eta,\_\rk_i:= ({\rm id}\otimes\id\otimes {\rm inf}_i)(\eta)$.
We just compute:
\begin{eqnarray*}
\gamma(\chi+\psi)&=&
\Phi\circ \lk\chi+\psi,\_\rk_1\circ\Phi^{-1} 
\circ \lk\chi+\psi,\_\rk_2^{-1}\\
&=&
\Phi\circ \lk\chi,\_\rk_1\circ\Phi^{-1} 
\circ
\Phi\circ \lk\psi,\_\rk_1\circ\Phi^{-1} 
\circ \lk\chi+\psi,\_\rk_2^{-1}\\
&=&
\Big(\Phi\circ \lk\chi,\_\rk_1\circ\Phi^{-1} 
\circ
\lk\chi,\_\rk_2^{-1}\Big)
\circ
\lk\chi,\_\rk_2\circ\\
&&\Big(\Phi\circ \lk\psi,\_\rk_1\circ\Phi^{-1} 
\circ \lk\psi,\_\rk_2^{-1}\Big)
\circ \lk\chi,\_\rk_2^{-1}\\
&=&
\gamma(\chi)\circ
\lk\chi,\_\rk_2\circ\gamma(\psi)
\circ \lk\chi,\_\rk_2^{-1}.
\end{eqnarray*}
The multiplication action on the unitary group of the
multiplier algebra turns into conjugation when passsing
to the inner automorphisms by ${\rm Ad}$,
so we obtain the cocycle identity 
$\gamma(\chi+\psi)=\gamma(\chi)\circ \big(\chi|_N\cdot\gamma(\psi)\big)$
for $\gamma$.
In particular, it follows that 
$\gamma(n^\perp+\psi)=\gamma(n^\perp)\circ \gamma(\psi)=\gamma(\psi)$,
for $n^\perp\in N^\perp$, i.e. 
$\gamma$ is constant on the cosets.
\end{proof}

The cocycle  $\dot\gamma$ determines a cohomology class
which -- by abuse of notation -- we again denote by 
$$
\gamma\in H^1\Big(\widehat N,
{\rm Inn}\big(\K\otimes\K\otimes(\C\rtimes_{\omega_2} N)\big)\Big).
$$
By the (not very long) exact sequence it defines an obstruction class
$$
\delta(\gamma)\in H^2(\widehat N, C(\widehat S,\UU)).
$$
If this obstruction class vanishes, $\gamma$ has a unitary lift, 
and this lift then implements a Morita equivalence between
the two actions ${\inf}_i,i=1,2$.
The next Lemma shows that in this case the classes 
$[\omega_1],[\omega_2]\in H^2(N,\UU)$ agree if
we assume that these classes extend to $G$:

\begin{lem}
Let
$[\omega_1],[\omega_2]\in H^2(N,\UU)$ be two classes in the image of 
the restriction map $H^2(G,\UU)\to H^2(N,\UU)$. If 
$$
\C\rtimes_{\omega_1} N \sim \C\rtimes_{\omega_2} N
$$
is a $\widehat G$-equivariant Morita equivalence, 
where $\widehat G$ acts on both sides by the inflated actions, 
then
$$
[\omega_1]=[\omega_2]\in H^2(N,\UU).
$$
\end{lem}

\begin{proof}
If $\mu_i$ is a $G$-action on $\K$ with Mackey obstruction $[\omega_i]$,
then $\K\otimes C(G/N)$ with $\mu_i\otimes{\text{(left translation)}}$
is the (pre-)dual of  $\C\rtimes_{\omega_i} N$, $i=1,2$.
But these systems are induced systems
$
\Ind^{G}_N(\K,\mu_i|_N)
$
by \cite[Theorem]{ech-ind}
which are classified up to Morita equivalence 
by their actions $\mu_i|_N:N\to{\Aut}(\K)$.
\end{proof}

Under certain circumstances one might 
indeed conclude that the obstruction class
$\delta(\gamma)$ vanishes:

\begin{cor} 
Assume that $N$ is torsion-free, and assume 
$\omega_2$ (or $\omega_1$) is totally skew.
If there is a $\widehat G$-equivariant Morita equivalence
$$
\K\otimes (\C\rtimes_{\omega_1}N)\sim \K\otimes (\C\rtimes_{\omega_2}N),
$$
with diagonal $\widehat G$-actions on both algebras 
as assumed in (\ref{EqMoritaEq}), then
$$
[\omega_1]=[\omega_2]\in H^2(N,\UU).
$$
\end{cor}
\begin{proof}
A totally skew cocycle has trivial symmetry group $S=\{0\}$, 
and if $N$ is torsion-free, then $\widehat N$ is connected,
and so the whole obstruction group vanishes:
$$
H^2(\widehat N, C(\widehat S,\UU))=H^2(\widehat N,\UU)\hookrightarrow{\rm Hom}(\widehat N,N)=0.
$$
\end{proof}

In Example \ref{ExaTwoDimAbiguity} both of the involved
non-commutative tori were rational. However, if one of
the involved classes is irrational, then its symmetry group is
trivial, so one can apply the above corollary.

\begin{exa}
Let $G=\R^2, N=\Z^2$, and let $r\in \R/\Z\cong H^2(\Z^2,\UU)$ 
be an irrational number, and let
$\omega_r$ be a corresponding cocycle. 
If there is a $\R^2$-equivariant Morita equivalence
$$
\K\otimes (\C\rtimes_{\omega} \Z^2)\sim \K\otimes (\C\rtimes_{\omega_r} \Z^2),
$$
where $\R^2$ acts on both algebras diagonally as in (\ref{EqMoritaEq}), 
then 
$$
[\omega]=[\omega_r]\in H^2(\Z^2,\UU).
$$
\end{exa}

Let us continue with some further analysis of the class $\delta(\gamma)$ constructed
above. Consider the following diagram of induced maps
$$
\xymatrix{
H^2(\widehat N, C(\widehat S,\UU))\ar[d]^{q^*}&\hspace{-3.4cm}\ni \delta(\gamma)&\\
H^2(\widehat G, C(\widehat S,\UU))\ar[d]^{i^*}&
H^2(\widehat G,\UU)\ar[d]^{i^\#}\ar[l]_-{c_*}&\hspace{-1cm}\ni{\rm Ma}(\hat\mu_1^{op}\otimes\hat\mu_2)\ar@<-5ex>@{|->}[d]^?\\
H^2(N^\perp, C(\widehat S,\UU))\ar@<-0.6ex>[r]_-{{\rm ev}_\#}&
H^2(N^\perp,\UU)\ar@<-0.6ex>@{_(->}[l]_-{c_\#}\ar@{}[r]|{\resizebox{0.27cm}{!}{$\ni$}}
&\hspace{-1.4cm} 0,\\
}
$$
where $c:\UU\to C(\widehat S,\UU)$ is the obvious inclusion, and 
${\rm ev}:C(\widehat S,\UU)\to \UU$ is the evaluation at $0\in \widehat S$
which is  $N^\perp$-equivariant, because $N^\perp$ acts trivially on both sides.
As ${\rm ev}_\#\circ c_\#={\rm id}$, $c_\#$ is injective.
It is easily seen that the Mackey obstruction class 
${\rm Ma}(\hat\mu_1^{op}\otimes\hat\mu_2)$ of the action 
$\hat\mu_1^{op}\otimes \hat\mu_2:\widehat G\to{\rm Aut}(\K\otimes\K)$, 
is connected to $\delta(\gamma)$ by\footnote{The sign reflects
the fact that the Mackey obstruction is defined to be the negative of the 
actual connecting homomorphism in (\ref{EqTheMaSequence}).
}
$
q^*(\delta(\gamma))=-c_*({\rm Ma}(\hat\mu_1^{op}\otimes\hat\mu_2)).
$
But the composition $i^*\circ q^*$ is zero, so by commutativity of the square
and by injectivity of $c_\#$, we see that 
$i^\#({\rm Ma}(\hat\mu_1^{op}\otimes\hat\mu_2))=0$, and we have just proven
the following statement:

\begin{lem}
\label{LemTheDualMackeyObstructions}
Let
$
\K\otimes (\C\rtimes_{\omega_1} N)\sim \K\otimes (\C\rtimes_{\omega_2} N)
$
be  a $\widehat G$-equivariant Morita equivalence, where $\widehat G$ acts on 
both sides diagonally by actions on the compacts 
$\hat\mu_i:\widehat G\to{\rm Aut}(\K)$, $i=1,2$, tensor
the inflated actions 
${\rm inf}_i:\widehat G\to {\rm Aut}(\C\rtimes_{\omega_i} N),i=1,2,$ 
on the respective crossed products.
Then the Mackey obstructions of the restricted actions
$\hat\mu_i|_{N^\perp}$ agree:
$${\rm Ma}(\hat\mu_1|_{N^\perp})={\rm Ma}(\hat\mu_2|_{N^\perp})\in H^2(N^\perp,\UU).$$ 
\end{lem}

An instance of Lemma \ref{LemTheDualMackeyObstructions} 
already appeared in Example \ref{ExaTwoDimAbiguity},
wherein the $\R^2$-action $\hat\mu_3$ 
has trivial Mackey obstruction when restricted to $\Z^2$.
In general, Lemma \ref{LemTheDualMackeyObstructions}  gives 
an obstruction map
\begin{eqnarray}
\label{EqTheClassiMapForNCT}
{\rm Ma}_{(N;\widehat G)}:{[}NCT(N;\widehat G){]}\to H^2(N^\perp,\UU)
\end{eqnarray}
{(here the brackets $[.]$ denote 2-isomorphism classes, 
i.e. Morita equivalence classes).
This map will enable us to identify the ``commutative theory" 
inside our theory as we explain next.

If  $\hat\mu\otimes\inf$ is a transverse 
$\widehat G$-action
on 
$\K\otimes (\C\rtimes_\omega N)$,
then the previous section 
shows that
its crossed product 
$\K\otimes (\C\rtimes_\omega N)\rtimes_{\hat\mu\otimes\inf} \widehat G$
together with the dual $G$-action is Morita equivalent to 
an algebra 
$\K\otimes(\C\rtimes_{\hat\eta} N^\perp)$ equipped with a
transverse action $\mu\otimes\inf$.
Corollary \ref{CorSuffCond} determines 
the class $[\hat\eta]\in H^2(N^\perp,\UU)$ to a certain
extend. Namely
\begin{eqnarray}
\label{EqTheDualAlgClass}
[\hat\eta]\in ({\rm Ma}(\hat\mu)|_{N^\perp})^\tri\subset H^2(N^\perp,\UU),
 \end{eqnarray}
where, for $\hat\theta\in H^2(N^\perp,\UU)$, we have used the notation
\begin{eqnarray}
\label{EqDefiOfThetaPerp}
\hat\theta^\tri :=
\Big\{ 
[\hat\omega\rtimes\omega]|_{N^\perp}
\ :\ [\hat\omega]|_{N^{\perp}}=\hat\theta,\omega\tra\hat\omega
\Big\}.
\end{eqnarray}
For $\hat\theta$ in the image of the restriction map 
$\widehat{\rm res}:H^2(\widehat G,\UU)\to 
H^2(N^\perp,\UU)\big)$ 
this set is never empty.
In fact, the trivial cocycle $1$ is transverse 
to any cocycle $\hat\omega$, so for $\hat\theta=[\hat\omega]|_{N^\perp}$
we have
$$
\hat\theta= [\hat\omega]|_{N^\perp}=
[\hat\omega\rtimes1]|_{N^\perp}
\in\hat\theta^\tri.
$$
The element relation (\ref{EqTheDualAlgClass}) restricts 
the possible classes for building the dual algebra 
to the set $({\rm Ma}(\hat\mu)|_{N^\perp})^\tri$.
This set is a Morita invariant of both the original system
$(\K\otimes(\C\rtimes_{\omega}N),\hat\mu\otimes\inf)$
and of its dual system. (Just because of Lemma \ref{LemTheDualMackeyObstructions} 
which states that the class ${\rm Ma}(\hat\mu)|_{N^\perp}$ 
is a Morita invariant of the original system.)
However, the element $[\hat\eta]$ itself
is not a specified element in 
$({\rm Ma}(\hat\mu)|_{N^\perp})^\tri$.
In fact, it varies with the choices 
of the extension of $\omega$ from $N$ to $G$, and 
these choices can make a difference
(see section \ref{Sec2DimNCT} for a detailed example).
Nevertheless, 
the following lemma shows that 
in a very important case the set 
$({\rm Ma}(\hat\mu)|_{N^\perp})^\tri$
reduces to a singleton.

\begin{lem}
\label{LemAboutThetaTransverse}
For $\hat\theta\in H^2(N^\perp,\UU)$ 
the following three statements are equivalent:
\begin{enumerate}
\item
$\hat\theta=0$,
\item
$\hat\theta^\tri =\{0\},
$
\item
$0\in\hat\theta^\tri.
$
\end{enumerate}
\end{lem}

\begin{proof}
Let us compute the image
of ${\hat\theta^\tri}$ under
$
h_{\_}:H^2(N^\perp,\UU)\hookrightarrow {\rm Hom}(N^\perp,G/N).
$
Let $x:=[\hat\omega\rtimes\omega]|_{N^\perp} \in {\hat\theta^\tri}$.
In the proof of Lemma \ref{LemDualityOnTransversalitySet}
we have already seen that
\begin{eqnarray*}
h_{\hat\omega\rtimes\omega}
&=&\phi\circ h_{\hat\omega},
\end{eqnarray*}
for $\phi = \id+ h_{\hat\omega}\circ h_\omega$.
By transversality, $\phi:G {\stackrel{\cong}{\to}} G$ induces an isomorphism
$\dot\phi:G/N {\stackrel{\cong}{\to}} G/N$. Now, the commutativities of the outer square,  of the bottom triangle
and of the two trapezoids in
$$
\xymatrix{
N^\perp\ar[rr]^{h_x}\ar@{^(->}[ddd]\ar[dr]_{h_\theta}&& G/N\\
&G/N\ar[ur]_{\dot\phi}&\\
& G\ar@{->>}[u]\ar[dr]^\phi&\\
\widehat G\ar[rr]_{h_{\hat\omega\rtimes\omega}}
\ar[ur]^-{h_{\hat\omega}}
&&G\ar@{->>}[uuu]
}
$$
imply that the upper triangle commutes, i.e. $h_x=\dot\phi\circ h_\theta$.
But as $\dot\phi$ is an isomorphism, $h_x$ vanishes if and only if $h_\theta$
vanishes. The lemma is then obvious.
\end{proof}

The important thing about this last lemma is that we 
have found an invariant that
can distinguish the commutative systems,
i.e. those which are equivariantly Morita equivalent 
to a system $(\K\otimes C(G/N),\mu\otimes\inf)$,
from those who are genuinely non-commutative.
Let us denote by 
$CT(N;\widehat G)\subset NCT(N;\widehat G)$ the
corresponding subcategory. 
Together with the results of the previous section 
we have found:
\begin{thm}
\label{ThmCommutativeSubTheory}

The set of 2-isomorphism classes 
$[CT(N;\widehat G)]$ 
is the kernel of the composition
$[NCT(N;\widehat G)^\tri]\to [NCT(N^\perp;G)^\tri]\to H^2(N,\UU)$
in
$$
\resizebox{!}{0.52cm}{
$
\xymatrix{
H^2(N^\perp,\UU)&&
\ar[ll]_-{{\rm Ma}_{(N;\widehat G)}}
[NCT(N;\widehat G)]& [NCT(N^\perp;G)]
\ar[rr]^-{{\rm Ma}_{(N^\perp; G)}}
&&H^2(N,\UU)
\\
&&
\ar@{}[u]|{{\resizebox{0.4cm}{!}{$\cup$} }}
[NCT(N;\widehat G)^\tri]\ar@<0.31ex>[r]&\ar@<0.31ex>[l] 
\ar@{}[u]|{{\resizebox{0.4cm}{!}{$\cup$} }}
[NCT(N^\perp;G)^\tri]\\
&&
\ar@{}[u]|{{\resizebox{0.4cm}{!}{$\cup$}}}
[CT(N;\widehat G)]&
\ar@{}[u]|{{\resizebox{0.4cm}{!}{$\cup$} }}
[CT(N^\perp; G)].
}
$
}$$
\end{thm}

\subsection{Example: Duality for the NC-Torus in Dimension 2}
\label{Sec2DimNCT}
Let $G:=\R^2=\widehat G$, and $N:=\Z^2=N^\perp$.
Recall the isomorphisms $H^2(\Z^2,\UU)\cong \T$ and 
$H^2(\R^2,\UU)\cong \R$.
Let us identify the transversality relation $\tra\subset \R\times\R$:
A cocycle corresponding to $\theta\in \R$ is given by
$\omega_\theta(x,y):=\exp(2\pi i \theta x_2y_1)$, so one obtains
$$
h_\theta=
\left( 
\begin{matrix}
0  & \theta \\
-\theta & 0
\end{matrix}
\right):\R^2\to \R^2.
$$
Then for some $\hat\theta\in \R$ 
we have 
$$
\phi=\id_{\R^2}+h_{\hat\theta}\circ h_\theta=(1-\theta\hat\theta)\cdot\id_{\R^2}
$$
which is invertible as long as $\theta\hat\theta\not=1$. 
Moreover,
it restricts to an isomorphism of $\Z^2$ if and only if 
$\theta\hat\theta=0$ or  
$\theta\hat\theta=2$.
Therefore $\tra \subset\R\times \R$
consists of the union of the coordinate axis 
and of the graph of $x\mapsto \frac{2}{x}$ 
(s. Figure  \ref{VeryTransverse}).
\begin{figure}
\includegraphics[scale=0.3]{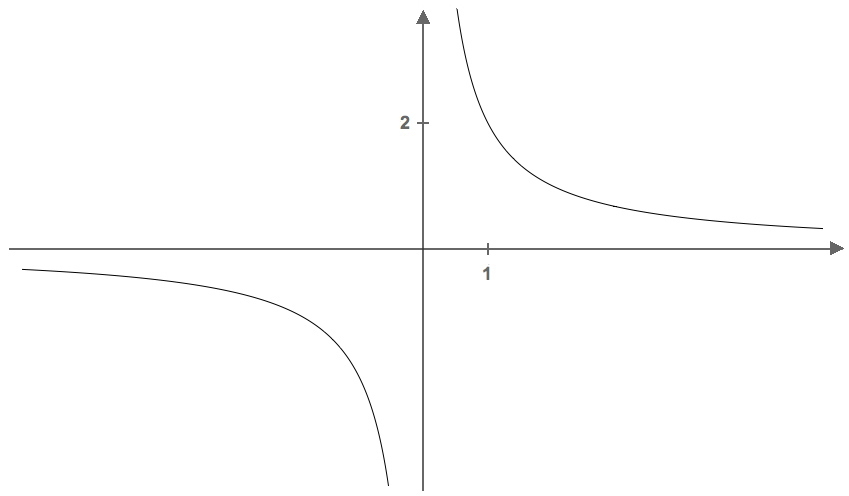}
\caption{$\tra \subset\R\times \R$}
\label{VeryTransverse}
\end{figure}

We can now answer the following question about transverse actions completely:
\begin{itemize}[leftmargin=*]
\item[Q:]
If $\dot \theta\in \T$, what are the transverse 
actions $\hat\mu\otimes\inf$ on 
$\K\otimes(\C\rtimes_{\dot\theta} N)$?
\item[A:]
If $\dot\theta=\dot 0\in \T$, then every action $\hat\mu\otimes\inf$ is 
transverse.
If $\dot\theta\not=\dot 0\in \T$, and $\theta\in\R$ is some lift of $\dot\theta$,
then every action $\hat\mu\otimes\inf$ is 
transverse which has a Mackey obstruction 
$$
{\rm Ma}(\hat\mu)\in\left\{0,\frac{2}{\theta+n}\ \Big|\  n\in\mathbb Z\right\}.
$$ 
\end{itemize}

Let us compute the dual system of
$(\K\otimes(\C\rtimes_{\dot\theta} \Z^2),\hat\mu\otimes\inf)$, 
for a transverse action $\hat\mu$ with ${\rm Ma}(\hat\mu)=\hat\theta$:
\begin{itemize}[leftmargin=*]
\item[(i)]
In the simple case of $\dot\theta=0$ we can choose $\theta=0$ as
a transverse lift to $\hat\theta$. Then we have
$$
\K\otimes C(\T)\rtimes_{\hat\mu\otimes\inf}\R^2
\sim
\K\otimes(\C\rtimes_{\dot{\hat\theta}} \Z^2),
$$
on which $\R^2$ acts by $\id\otimes\inf$, and $\dot{\hat\theta}$ is 
the restriction of $\hat\theta$ to $\Z^2$.
\item[(ii)]
In case of $\dot\theta\not=\dot 0$, there are two subcases.
First, if $\hat\theta=0$, then let $\theta$ be any lift of 
$\dot\theta$, and we have
$$
\K\otimes(\C\rtimes_{\dot\theta} \Z^2)\rtimes_{\hat\mu\otimes\inf}\R^2
\sim
\K\otimes C(\T^2),
$$
on which $\R^2$ acts by $\mu\otimes\inf$ with ${\rm Ma}(\mu)=\theta$.
Second, if $\hat\theta\not=0$, then transversality means that there 
is a lift $\theta$ of $\dot\theta$
such  that $\hat\theta=2/\theta$.
We now only need to compute the cocycles 
occurring in Corollary \ref{CorTheDualOverThePoint}: 
\begin{eqnarray*}
h_{\hat\theta}^*\omega_{\theta}^{-1}(x,y)
&=&\exp(2\pi i \theta (-\hat\theta x_1)(\hat\theta y_2) )^{-1}\\
&=&\exp(2\pi i 2 \hat\theta x_1 y_2 )\\
&\sim&\exp(2\pi i  2 \hat\theta x_2 y_1 )^{-1},
\end{eqnarray*}
so
$[\hat\omega_{\hat\theta}\cdot h_{\hat\theta}^*\omega_{\theta}^{-1}]=-\hat\theta$, and similarly
$[\omega_\theta\cdot h_\theta^*\hat\omega_{\hat\theta}^{-1}]=-\theta$.
Pushing these classes forward along $\phi=\hat\phi=-\id_{\R^2}$ doesn't
change the class.
So if we denote by $\dot{\hat\theta}\in\T$ the restricted class of $\hat\theta$
and if $\mu_{-\theta}$ is an action with Mackey obstruction $-\theta$, 
then the dual is given by
\begin{eqnarray}
\label{EqSomeDualOverHere}
(\K\otimes(\C\rtimes_{-\dot{\hat\theta}} \Z^2),\mu_{-\theta}\otimes\inf).
\end{eqnarray}
\end{itemize}
Note that we re-obtain Example \ref{ExaTwoDimAbiguity} at this stage:
In fact, for $\theta=\frac{2}{3}$ and $\hat\theta=3$ the dual
system of $(\K\otimes(\C\rtimes_{\dot\theta} \Z^2),\hat\mu\otimes\inf)$
is according to (\ref{EqSomeDualOverHere})
\begin{eqnarray*}
(\K\otimes C(\T^2),\mu_{-\frac{2}{3}}\otimes\inf).
\end{eqnarray*}
Then we might take the bidual according to (i)
which is
\begin{eqnarray*}
(\K\otimes (\C\rtimes_{{-\frac{2}{3}}}\Z^2),\id\otimes\inf).
\end{eqnarray*}
However, 
$\frac{1}{3}=-\frac{2}{3}\ {\rm mod}\ \Z$,
and this is exactly  
what we observed in Example \ref{ExaTwoDimAbiguity}.

\subsection{Duality for 3-dimensional NC-tori}
We consider now the standard lattice $\Z^n\subseteq \R^n$. Recall from \cite{BK}  that
for every  $2$-cocycle $\omega$ on $\R^n$ there is a unique  strictly upper triangular 
real
$n\times n$ matrix $A$ such that 
$\omega$ is cohomologous to $\omega_A:\R^n\times\R^n\to\T$ given by
$\omega_A(x,y)=e^{2\pi i\lk A x,y\rk}$. 
Moreover, every $2$-cocycle on $\Z^n$ is cohomologous 
to the restriction of some $\omega_A$ to $\Z^n\times \Z^n\subseteq \R^n\times \R^n$ and two 
such restrictions lie in the same class in $H^2(\Z^n,\T)$  
  if and only if the difference $A-A'$  has  integer entries. If we identify $\R^n$ with $\widehat{\R^n}$ 
via $x\mapsto\chi_x$ with $\chi_x(y)=e^{2\pi i \lk x,y\rk}$, then a straight-forward computation shows 
that $h_{\omega_A}:\R^n\to \R^n$ 
is given by matrix multiplication with the skew symmetric matrix $A-A^t$. 
Now, given another strictly upper triangular matrix $B$ defining a $2$-cocycle $\widehat{\omega}=\omega_B$ on 
$\widehat{\R^n}\cong \R^n$, then the pair $(\omega_A, \omega_B)$ is transverse if and only if the matrix
$$\Phi:= I_n+ (B-B^t)(A-A^t)$$
lies in $\GL(n,\Z)$. As a consequence, given an action  $\widehat\mu\otimes\inf$ of $\widehat{\R^n}\cong \R^n$ 
on \mbox{$\K\otimes (\C\rtimes_\eta\Z^n)$} such that ${\rm Ma}(\widehat\mu)=[\omega_B]$, then the system 
$\big( \K\otimes (\C\rtimes_\eta\Z^n), \widehat{\R^n}, \widehat\mu\otimes \inf\big)$ is transversal if and only 
if there exists a strictly upper triangular matrix $A$ such that $\eta\sim \omega_A|_{\Z^n\times Z^n}$ and such that 
$I_n+ (B-B^t)(A-A^t)\in \GL(n,\Z)$. 

 Of course, this will always be the case if $B=0$, which results to give Mathai-Rosenberg duality between commutative and 
 noncommutative tori.  
 For the case $n=2$ we completely solved this question in the previous example. In particular we saw that 
 every non-commutative $2$-torus has non-commutative duals. We shall now show that this is not always the case 
 in higher dimensions.
 
 For this we have a look 
at the case $n=3$: Suppose that 
$$A=\left(\begin{matrix} 0 & a_1 & a_2\\ 0&0& a_3\\ 0&0&0\end{matrix}\right)\quad\text{and}\quad 
B=\left(\begin{matrix} 0 & b_1 & b_2\\ 0&0& b_3\\ 0&0&0\end{matrix}\right).$$
 Then a short computation
shows that 
\begin{align*}&\Phi:=I_3+(B-B^t)\cdot (A-A^t)\\
&=\left(\begin{matrix} 1-(b_{1}a_{1}+b_{2}a_{2})  & -b_{2}a_{3} &b_{1}a_{3}\\
-b_{3}a_{2} & 1-(b_{1}a_{1}+b_{3}a_{3}) & -b_{1}a_{2}\\
b_{3}a_{1} &- b_{2}a_{1} &1-(b_{2}a_{2}+b_{3}a_{3})\end{matrix}\right).
\end{align*}
A lengthy but straightforward computation gives
\begin{equation}\label{eq-det}
\det\Phi= 1-2(c_1+ c_2+c_3)+ 2(c_1c_2+c_1c_3+c_2c_3)+ c_1^2+c_2^2+c_3^2,
\end{equation}
with $c_i:=a_{i}b_{i}$, $i \in \{1,2,3\}$.  In order that 
$\Phi\in \GL(3,\Z)$ we need $\det\Phi=\pm 1$. Since all entries of the matrix $\Phi$ need to be integers,
it follows that all mixed terms $c_ic_j$ with $i\neq j$ are integers as well. Also, all sums $c_i+c_j$ with $i\neq j$ 
need to be integers. Using this, it easily follows that $2c_i\in \Z$ for all $i\in\{1,2,3\}$ (e.g., we get
$2c_1= (c_1+c_2)+(c_1+c_3)-(c_2+c_3)\in \Z$). But then it follows that $c_1^2+c_2^2+c_3^2\in \Z$.
It is an easy exercise to check that the sum of three squares of half-integers is an integer 
if and only if all three numbers are integers themselves. Thus we get $c_1, c_2, c_3\in \Z$ and 
 we need to search for integer solutions $c_1, c_2, c_3$ for the equation $\det \Phi=\pm 1$.
Using standard forms for quadratic equations one can check that there are precisely the following two sets of solutions
 $$ c_3= -c_1-c_2 \quad\text{or}\quad c_3=2-c_1-c_2 \quad \forall c_1, c_2\in \Z,$$
 in which cases the determinant will be $1$.
 Moreover, the matrix $\Phi$ will have integer entries if and only if 
 $n_{ij}:=a_ib_j\in \Z$ for all $1\leq i,j\leq 3$. If  $B=0$, the 
 conditions will always be fulfilled.  On the other hand, if just one entry $b_i\neq 0$, this will force
 any of the pairs $a_j, a_k$ to be rationally dependent, since if both $a_j, a_k\neq 0$,
 then the equations $0\neq a_j b_i=n_{ji}$ and $0\neq \alpha_k b_i=n_{ki}$ implies 
 that $a_j=\frac{n_{ji}}{n_{ki}} a_k$. 
 
Thus, a general non-commutative 3-torus which admits a non-commutative dual will be attached to a matrix 
of the form $A=\theta \left(\begin{matrix} 0& \alpha_1&\alpha_2\\ 0&0&\alpha_3\\0&0&0\end{matrix}\right)$
with $0\neq \theta \in \R$ and $\alpha_1,\alpha_2,\alpha_3\in \mathbb Q$. Then 
 $B=\frac{1}{\theta}\left(\begin{matrix} 0& \beta_1&\beta_2\\ 0&0&\beta_3\\0&0&0\end{matrix}\right)$
 with $\beta_1,\beta_2,\beta_3\in \mathbb Q$ such that $\alpha_i\cdot \beta_j\in \Z$ for all $1\leq i,j\leq 3$ and 
 $\alpha_1\beta_1+\alpha_2\beta_2+\alpha_3\beta_3\in \{0,2\}$ will match up with $A$ to give a  dual pair of non-commutative three dimensional tori. 
 As an example: if all $\alpha_i=1$, then $(\beta_1,\beta_2,\beta_3)=(n,k,-n-k)$ as well as 
$(\beta_1,\beta_2,\beta_3)=(n,k,2-n-k)$, with $n,k\in \Z$, will give compatible parameters for $B$.
 
 \subsection{Example: Tensor products of duality pairs}
 
 Duality pairs are closed under taking tensor products: Suppose that $(\omega_A, \omega_B)$ is a duality pair for 
 $\Z^n\subseteq \R^n$ and $(\omega_C,\omega_D)$ is a duality pair for $\Z^m\subseteq \R^m$. 
 Then $(\omega_A\cdot\omega_C,\omega_B\cdot\omega_D)$ is a duality pair for $\Z^{n+m}\subseteq \R^{n+m}$,
 where $\omega_A\cdot\omega_C$ is the cocycle on $\R^{n+m}$ given by the product
 $$\omega_A\cdot\omega_C((x_1,y_1), (x_2, y_2))=\omega_A(x_1, x_2)\omega_C(y_1,y_2)\quad\forall\; (x_1,y_1),(x_2,y_2)\in \R^{n+m}$$
(and similarly for $\omega_B\cdot\omega_D$). Note that $\omega_A\cdot\omega_C=\omega_{\diag(A,C)}$ if we denote by 
$\diag(A,C)$ the matrix $\left(\begin{matrix} A&0\\0&B\end{matrix}\right)\in M_{n+m}(\R)$. 
Since
\begin{align*}
&I_{n+m}  + \big(\diag(B,D)-\diag(B,D)^t\big)\big(\diag(A,C)-\diag(A,C)^t\big)\\
&=\diag\Big(I_n+(B-B^t)(A-A^t), I_m+(D-D^t)(C-C^t)\Big)
\end{align*}
it follows that  $(\omega_A\cdot\omega_C,\omega_B\cdot\omega_D)$ is a dual pair if and only if 
$(\omega_A, \omega_B)$ and $(\omega_C, \omega_D)$ are both dual pairs.

\section{Non-Commutative C*-Dynamical  T-Duality}
\label{SecGlobalNCTDuality} 

In classical (commutative) C*-dynamical T-duality the  objects are
principal torus bundles $E\to B$ equipped with 
a locally trivial bundle of compact operators $F\to E$
that is trivialisable over the fibres $E_b\to b$.
The C*-algebra of sections 
vanishing at infinity $A:= \Gamma_0(E,F)$
is a bundle of C*-algebras
whose fibres are stable commutative tori $\K\otimes C(\T^n)$.
We consider more general bundles whoses fibres 
are twisted group algebras.

In this section the word \emph{space} will always 
mean a second countable, locally compact Hausdorff space.

\subsection{{$\boldsymbol{C_0(B)}$}-Algebras and 
Continuous Bundles of C*-Algebras
}
\mbox{}
{Recall that} a C*-algebra $A$ is called a {\em $C_0(B)$-algebra} 
for a space $B$, if 
$A$ is equipped with a fixed non-degenerate $*$-homomorphism 
$\Phi:C_0(B)\to {\rm ZM}(A)$, 
the center of the multiplier algebra ${\rm M}(A)$ of $A$. 
If $A$ is a $C_0(B)$-algebra, then for 
any closed subset  $X\subset B$ we let $I_X$ denote the closed ideal 
$\Phi(C_0(B\setminus X))A$ of $A$ 
and the quotient $A|_X:=A/I_X$ is called the {\em restriction of $A$ to $X$}.
For a single point $\{b\}=X$, $A|_b$ is called the {\em fibre of $A$ over $b$}. 
The elements of $A$ can be viewed as sections of a fibre-bundle over {$B$} with fibres 
{$A|_b$} 
by writing $a(b):=a+I_b\in A|_b$, for  $a\in A$ and $b\in B$. 
A $C_0(B)$-algebra $A$ is called
{\em a continuous bundle of C*-algebras over $B$}
if these sections are continuous in the sense 
that $b\mapsto \|a(b)\|$ is continuous.  
We refer to \cite[Appendix C]{Wbook} for a detailed treatment of $C_0(B)$-algebras.

An action $\alpha:G\to \Aut(A)$ on a 
$C_0(B)$-algebra $A$, is called {\em fibre-preserving}
if it is $C_0(B)$-linear, i.e. if $\alpha_g(f\cdot a)=f\cdot \alpha_g(a)$
 for all $a\in A$, $g\in G$ and $f\in C_0(B)$, where we write $f\cdot a$ for $\Phi(f)a$.
If $\alpha$ is $C_0(B)$-linear, it induces actions $\alpha|_b$ on each fibre $A|_b$ 
 by setting $(\alpha|_b(g))(a(b)):=(\alpha_g(a))(b)$. 
 Then $\alpha$ is completely determined by 
 the actions on the fibres.
 
 If $A_1$ and $A_2$ are two $C_0(B)$-algebras, then we say that an 
 $A_1$-$A_2$-equivalence bimodule $E$
is $C_0(B)$-linear, if $ f\cdot \xi=\xi\cdot f$ 
for all $f\in C_0(B)$, $\xi \in E$, where the left and right 
 actions of $C_0(B)$ on $E$ are given via  
 extending the left and right actions of {$A_1$} and $A_2$ 
 to their multiplier algebras. 
 We then say that $A_1$ and $A_2$ are $C_0(B)$-linearly Morita equivalent.

 We will typically deal with a situation in which the Morita equivalences 
 are assumed to be both, $G$-equivariant and $C_0(B)$-linear.

\subsection{Families of Twisted Group Algebras 
{and $\boldsymbol{\omega}$-Triviality}}

Let $B$ be a locally compact space, and
consider a $C_0(B)$-linear
action on $\K\otimes C_0(B)$, i.e.
a continuous homomorphism 
$\mu:G\to C(B,\PU(\HH))$, where  $C(B,\PU(\HH))$
is equipped with the compact open topology
(this is the $C_0(B)$-linear automorphism group
of $\K\otimes C_0(B)$ with the topology of point-wise convergence).
We can construct such a homomorphism $\mu$ 
literally by formula (\ref{EqTheLeftRegRep})
from a cocycle $\omega\in Z^2(G,C(B,\UU))$
by {the left regular $\omega$-representation}.
I.e. one just has to apply 
(\ref{EqTheLeftRegRep})
point-wise for ${\omega|_b}\in Z^2(G,\UU)$ (the evaluation of
$\omega$ at $b\in B$). The resulting function indeed is 
continuous as a map  $G\to C(B,\PU(\HH))$ 
\cite[Proof of Prop. 3.1]{HORR}.
This gives rise to a map
\begin{eqnarray}
\label{EqDefiningXi}
{\xi_{B,G}}:H^2(G,C(B,\UU))\to \mathcal E_G(B),
\end{eqnarray}
where $\mathcal E_G(B)$  denotes the Morita equivalence 
classes of systems $(\K\otimes C_0(B),\mu,G)$. 
It is shown in   \cite{CKRW} that
$\mathcal E_G(B)$ is a group by forming the balanced 
tensor product over $B$, and the above map is a 
homomorphism.
For the one-point space $B={\rm pt}$, Proposition \ref{prop-Mackey} 
tells that the above map is an
an isomorphism, namely the inverse of the Mackey obstruction.
Yet, for general $B$ this fails:
\begin{prop}[{\cite[Sec. 6.3]{CKRW}}] 
If the second $\check{\rm  C}$ech cohomology 
$\check H^2(B,\Z)$ is countable, 
then
there is an exact sequence
\begin{eqnarray*}
0\to H^2(G,C(B,\UU))\stackrel{\xi_{B,G}}{\longrightarrow} 
\mathcal E_G(B)\longrightarrow
 {\rm Hom}(G,\check H^2(B,\Z)).
\end{eqnarray*}
\end{prop}

In the remainder of this article
we will exclusively deal with actions $\mu$ on 
$\K\otimes C_0(B)$
which are in the image
of $\xi_{B,G}$.
Instead of analysing 
the crossed product $(\K\otimes C_0(B))\rtimes_\mu G$ 
we will therefore mostly consider
 twisted transformation group algebras
 defined as follows:
If $\omega\in Z^2(G,C(B,\UU))$,
then the Banach space $L^1(G,C_0(B))$ is
turned into a Banach *-algebra by the same formulas as in 
(\ref{EqOpForTwistedAlg}).
Its enveloping C*-algebra is denoted by 
$$C_0(B)\rtimes_\omega G.$$
It has a canonical action of the dual group $\widehat G$
by point-wise multiplication of characters.
Furthermore, it is a C*-algebra over $B$ and its fibres are 
$\big(C_0(B)\rtimes_\omega G\big)|_b\cong \C\rtimes_{{\omega|_b}} G$,
where again ${\omega|_b}$ denotes 
the restriction of $\omega$ to $b\in B$. (We refer to \cite{EW02} for a more detailed 
discussion of twisted transformation group algebras.)
If now $\mu:G\to C(B,\PU)$ is a given action
in the image of $\xi_{B,G}$,
then
there is a
$\widehat G$-equivariant and 
$C_0(B)$-linear
isomorphism
$$
\Psi :\K\otimes (C_0(B)\rtimes_\omega G)\cong 
(\K\otimes C_0(B))\rtimes_\mu G
$$
given literally by formula (\ref{EqIsoOfTwistedAndCrossed}).
In other words, we can go back and forth between families 
of twisted group algebras and their corresponding 
crossed product algebras without losing information.

We are now coming back to our general 
situation in which
$N\subseteq G$ is a cocompact discrete 
subgroup of the abelian group $G$.

\begin{defi}
\label{DefiOmegaTrivaial}
Let $B$ be a space and let $A$ be a (separable and stable) $C_0(B)$-algebra.
\begin{enumerate}
 \item
A is called {$\omega$\it-trivial} if there exists  
$\omega\in Z^2(G,C(B,\UU))$ together with a $C_0(B)$-linear
Morita equivalence $A\sim C_0(B)\rtimes_{\omega}N$.
The pair of such a cocycle together with such a Morita 
equivalence is called 
an {$\omega$\it-trivialisation}.

\item
$A$ is called {\it locally $\omega$-trivial} if
there exists an open covering $(U_i)_{i\in I}$ of $B$ such that for all
$i\in I$ the restricted algebras $A|_{\overline U_i}$ over the closure 
$\overline U_i\supset U_i$  are $\omega$-trivial.
\item 
The open sets together with their 
$\omega$-trivialisations are called 
{\it charts}, and a collection of charts covering 
all of $B$ is called an {\it atlas}.
\end{enumerate}
(We use the term {{\it (local) $\hat \omega$\it-triviality}} if 
we consider the groups $N^\perp,\widehat G$ instead of
$N,G$.)
\end{defi}
Note that although the Morita equivalence in this definition only refers 
to the cocycle on $N$, we require the cocycle to extend to $G$.

\begin{exa}
\label{ExaTwoLocTrivAlgs}
Let $G:=\R^2,N:=\Z^2$. Let $B:=[0,1]$  be 
the interval and 
$S^1:=B/(0\sim 1)$ be the circle.

\begin{enumerate}
\item
Let $\omega^{(1)}\in Z^2(\R^2,C(B,\UU))$ 
be  given by the class function\footnote{
I.e. $
\omega^{(1)}_t(x,y):=\exp(2\pi i [\omega^{(1)}_t] x_2y_1)$.
}
$t\mapsto [\omega^{(1)}_t]\in\R\cong H^2(\R^2,\UU)$ 
which is defined by Figure \ref{FancyOmega}.
\begin{figure}
\includegraphics[scale=0.45]{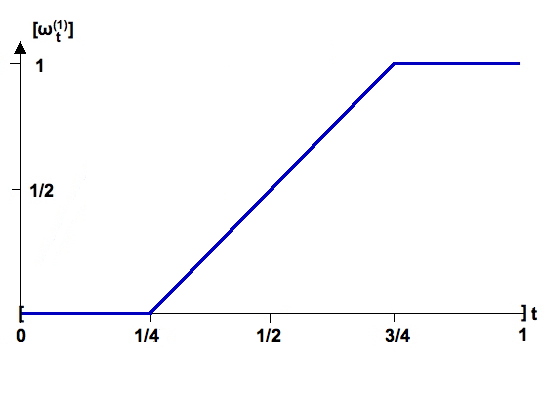}
\caption{A function for defining a bundle of algebras.}
\label{FancyOmega}
\end{figure}

Let $X:=C(B)\rtimes_{\omega^{(1)}} \Z^2$, which is a globally
$\omega$-trivial (non-stable) algebra over $B$.
Note that $\omega^{(1)}_{0}|_{\Z^2\times\Z^2}=\omega^{(1)}_1|_{\Z^2\times\Z^2}$, so
the fibres $X|_{0}$ and $X|_1$ over the two endpoints of $B$
are canonically isomorphic
(they are equal after identifying $X|_t\cong \C\rtimes_{\omega_t}\Z^2$
by the obvious map).
Let us denote by $A^{1}$ the algebra which is obtained 
by gluing along this isomorphism, i.e. the pullback in 
$$
\xymatrix{
A^{1}\ar@{.>}[d]\ar@{.>}[r] & X\ar[d]^{i_{0}^*\times i_1^*}\\
 X|_{0}\ar[r]^-{\id\times{\rm can}}&X|_{0}\oplus  X|_1
},
$$
and denote by $A^{(1)}$ its stabilisation $\K\otimes A^{1}$.
If we let $\eta\in Z^2(N,C(S^1,\UU))$ be the cocycle
given point-wise by $\omega^{(1)}|_{N\times N}$,
then
$A^{(1)}$ is canonically isomorphic to 
$\K\otimes C(S^1)\rtimes_{\eta}\Z^2$.
This is not an  $\omega$-trivial algebra any more, 
yet it is still a locally $\omega$-trivial algebra
over the circle $S^1$.
\medskip

\item
Let us do a slightly more involved construction
with the cocycle
$\omega^{(2)}\in Z^2(\R^2,C(B,\UU))$ whose class function is 
given by
Figure \ref{MoreFancyOmega}.
\begin{figure}
\includegraphics[scale=0.45]{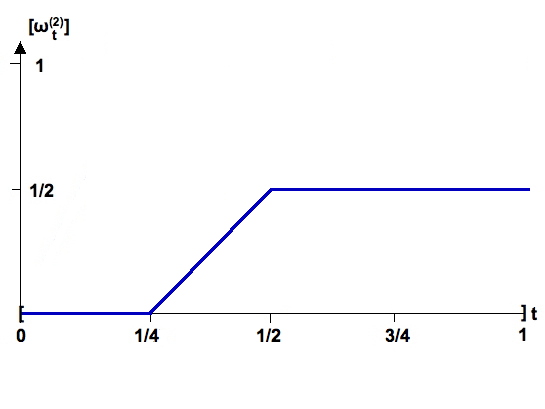}
\caption{Another function for defining a bundle of algebras.}
\label{MoreFancyOmega}
\end{figure}

Let $Y:= C(B)\rtimes_{\omega^{(2)}} \Z^2$.
In this case, the fibres $Y|_{0}$, $Y|_1$ over the two endpoints
of $B$ are two non-isomorphic but Morita equivalent 
non-commutative tori with classes 0 and $\frac{1}{2}$. 
So there is a stable isomorphism
$\varphi:\K\otimes Y|_{0}\cong \K\otimes Y|_1$. We use this 
isomorphism to glue the stabilisation $\K\otimes Y$ 
over the endpoints to itself:
Let $A^{(2)}$ be the pullback in
$$
\xymatrix{
A^{(2)}\ar@{.>}[d]\ar@{.>}[r] &\K\otimes Y\ar[d]^{i_{0}^*\times i_1^*}\\
\K\otimes Y|_{0}\ar[r]^-{\id\times\varphi}&(\K\otimes Y|_{0})\oplus (\K\otimes Y|_1)
}.
$$
We claim that $A^{(2)}$ is a locally $\omega$-trivial 
algebra over the circle $S^1$. It is clear that $Y$ itself gives a chart for, say, 
$V:=(\frac{1}{8},\frac{7}{8})\subset S^1$. The critical issue is to find a local
$\omega$-trivialisation around the gluing point: For, say,  
$U:=[0,\frac{1}{4})\sqcup (\frac{3}{4},1]/(0\sim 1)\subset S^1$
we have
$$
A^{(2)}|_{\overline U}\cong\Big\{(f,g) \in (\K\otimes Y|_{[0,-\frac{1}{4}]})\oplus
(\K\otimes Y|_{[\frac{3}{4},1]})\ \Big|\ \varphi(f|_{0})=g|_1\Big\}.
$$
But the isomorphism $\varphi$ extends to a fibre-wise isomorphism
$$\tilde\varphi:\K\otimes Y|_{[0,\frac{1}{4}]}\stackrel{\cong}{\to} \K\otimes Y|_{[\frac{3}{4},1]}$$
just because all fibres are the same. So we obtain
\begin{eqnarray*}
A^{(2)}|_{\overline U}
&\cong&\Big\{(\tilde\varphi f,g) \in (\K\otimes Y|_{[\frac{3}{4},1]]})\oplus
(\K\otimes Y|_{[\frac{3}{4},1]})\ \Big|\ (\tilde\varphi f)|_{1}=g|_1\Big\}\\
&\sim& C(\overline U)\rtimes_{\omega_{\frac{1}{2}}}\Z^2,
\end{eqnarray*}
and this means that $V,U$ give an atlas for $A^{(2)}$.
\end{enumerate}
\end{exa}
\medskip

\subsection{Duality for Polarisable Pairs}

Recall that for a {fibre-preserving} 
action $\alpha:G\to {\rm Aut}(A)$, 
i.e. a $C_0(B)$-linear action,
the crossed product C*-algebra $A\rtimes_\alpha G$ is 
again a $C_0(B)$-algebra.

\begin{defi}
\label{DefiNCPair}
Let $B$ be a space.
\begin{enumerate}
\item
An action $\alpha:G\to \Aut(A)$ on a locally 
$\hat\omega$-trivial  $C_0(B)$-algebra $A$ is called
\emph{locally transverse} if there exists 
an atlas $(U_i,\hat\omega_i)_{i\in I}$ for $A$ and 
cocycles $\omega_i\in Z^2(G,C(\overline U_i,\UU))$
which are point-wise transverse to  $\hat\omega_i$
together with $\widehat G$-equivariant and $C_0(\overline U_i)$-linear
Morita equivalences
\begin{eqnarray*}
(A\rtimes_\alpha G)|_{\overline U_i}\sim C_0(\overline U_i)\rtimes_{\hat\omega_i\wedge\omega_i}(N^\perp\times G).
\end{eqnarray*}
(The $\widehat G$-equivariance is required for the dual actions on both sides.)
These local data $U_i,\hat\omega_i,\omega_i$ and the Morita equivalences 
are called  {\it transverse charts} which altogether constitute a
{\it transverse atlas}.

\item
A {\it NC pair $(A,\alpha)$ over $B$}
is a locally $\hat\omega$-trivial $C_0(B)$-algebra $A$ together with a 
locally transverse action.
\item We use the term {\it NC dual pair} for a locally 
$\omega$-trivial algebra with  a locally transverse 
$\widehat G$-action.
\end{enumerate}
\end{defi}

\begin{rem}
Local transversality of an action $\alpha$ can be rephrased
in terms of the local actions $\alpha|_{\overline U_i}$ on 
$A|_{\overline U_i}\cong \K\otimes C_0(\overline U_i)\rtimes_{\hat\omega_i} N^\perp$.
It is equivalent to require, firstly, that
$\alpha|_{\overline U_i}$ factorises as
$\mu_i\otimes\inf$, and, secondly, that 
the element in $\mathcal E_G(\overline U_i)$ given
by $\mu_i:\overline U_i\to C(G,\PU(\HH))$ is in the image of 
$H^2(G,C(\overline U_i,\UU))\hookrightarrow \mathcal E_G(\overline U_i)$
such that its pre-image in 
$H^2(G,C(\overline U_i,\UU))$
is point-wise transverse
to $\hat\omega_i$.
\end{rem}

If $(A,\alpha)$ is a NC pair, then the crossed 
product $A\rtimes_\alpha G$ is a $C_0(B)$-algebra  and 
its fibres can be computed point-wise according to 
section \ref{SecItCrPrAndTrans} which 
shows that the dual action of $\widehat G$ on 
$A\rtimes_\alpha G$ is transverse in each fibre.
However, the crossed product algebra $A\rtimes_\alpha G$
need not to be locally $\omega$-trivial, and even if it is,
then the point-wise computation of 
section \ref{SecItCrPrAndTrans}
determines the dual action  only up to Morita
equivalence in each fibre which is in general not enough to determine 
a unique action with these properties. So we cannot 
simply conclude that the dual action is locally transverse.

\begin{defi}
\label{DefiPolarkreis}
Let $B$ be a space.
\begin{enumerate}
\item
A {\it local  polarisation} over $ U\subset B$ is
just a continuous family of group automorphims
$\varphi:\overline U\to \Aut(\widehat G\times G)$,
where $\Aut(\widehat G\times G)$ is equipped with 
the Bracconier topology. 
This is the topology generated by the compact open
topology and the pre-images of the open sets by 
inversion in ${\rm Aut}(\widehat G\times G)$. By the exponetial 
law (\cite{B64}) this means that $\varphi$ is continuos if and only if

\begin{eqnarray*}
\overline U\times \widehat G\times G 
&\to& \overline U\times \widehat G\times G\\
(u,\chi,g)&\mapsto&(u,\varphi(u)(\chi,g))
\end{eqnarray*}
is a homoeomorphism.

\item 
Let $(A,\alpha)$ be a NC pair,
and let $(U_i,\hat\omega_i,\omega_i)_{i\in I}$
be a transverse atlas. The atlas
is called {\it polarisable}
if there exist local polarisations $\varphi_i$ over $U_i$ such that
\begin{eqnarray*}
[\hat\omega_i\wedge\omega_i] = \varphi_i^*[\wedge] \in 
H^2(\widehat G\times G,C(\overline U_i,\UU)).
\end{eqnarray*}
\item A NC pair $(A,\alpha)$ is called polarisable 
if it permits a polarisable atlas.
\end{enumerate}
\end{defi}

Polarisability  is a slight restriction on the class of objects we consider.
In fact, if an atlas is polarisable, then it follows that the involved cocycles
$\hat\omega_i,\omega_i$ are cohomologous to bihomomorphisms, 
just because 
\begin{eqnarray}
\label{EqBihomomor}
\omega_i=(\hat\omega_i\wedge\omega_i)|_{(0\times G)\times (0\times G)}\sim
(\varphi^*_i\wedge)|_{(0\times G)\times (0\times G)},
\end{eqnarray}
and the latter cocycle is a bihomomorphism
$G\times G\to C(\overline U_i,\UU)$.
This means that the bundle theory does not cover all the cases of the
point-wise theory in the previous sections.
However, for $G=\R^n$ and $N=\Z^n$ all cocycles are
cohomologous to bicharacters, and, as we see next, polarisability is no restriction 
at all:

\begin{prop}
For $G=\R^n, N=\Z^n$ every 
NC pair is polarisable.
\end{prop}

\begin{proof}
Let $U_i,\hat\omega_i,\omega_i$ be a chart.
The first thing to note is that for $\R^{2n}$ we always have  
an isomorphism\footnote{
$Z^2(G,\UU)$ has the topology of 
almost everywhere point-wise convergence,
and $H^2(G,\UU)$ has the corresponding 
quotient topology.
Then $H^2(\R^{n},\UU)\cong \R^{{n(n-1)}/{2}}$ as
topological groups.
}
\cite[Sec. 5]{EW93}
$$
H^2(\R^{2n},C(\overline U_i,\UU))\cong C(\overline U_i,H^2(\R^{2n},\UU)).
$$
which is given by point-wise evaluation. 
Secondly, 
the Braconnier topology coincides with the usual 
topology on ${\rm Aut}(\R^n\times \R^n)=Gl(2n,\R)$.

Now, if $\hat\eta\wedge\eta$ is a  type I and totally skew cocycle
on $\R^{2n}$,
then $h_{\hat\eta\wedge\eta}:\R^{2n}\to\R^{2n}$
is an invertible, anti-symmetric 
matrix. 
But if $As(2n,\R)$ denotes the 
set of all anti-symmetric $2n\times2n$-matrices, then 
$$
\begin{tabular}{cc}
\begin{minipage}{6cm}
\begin{eqnarray*}
c:Gl(2n,\R)&\twoheadrightarrow& {Gl(2n,\R)\cap As(2n,\R)},\\
\varphi&\mapsto& \hat\varphi\circ h_\wedge \circ \varphi
\end{eqnarray*}
\end{minipage}
&
\begin{minipage}{4cm}
\begin{eqnarray*}
\text{  }
h_\wedge=\left(
\begin{matrix}
0&\Eins\\
-\Eins&0
\end{matrix}
\right),
\end{eqnarray*}
\end{minipage}
\end{tabular}
$$
is surjective and admits local sections.
(Here the dual map $\hat{\varphi}$ coincides with the transpose map $\varphi^t$.)
In fact, surjectivity is a standard fact about
anti-symmetric matrices. To conclude the existence 
of local sections, it is sufficient to observe that $c$ is a 
submersion, i.e. we claim that its derivative
\begin{eqnarray*}
dc(\varphi):Mat(2n,\R)&\to& As(2n,\R)\\
M&\mapsto& \hat M\circ \hat\varphi\circ h_\wedge+h_\wedge \circ \varphi\circ M
\end{eqnarray*}
is surjective for all $\varphi\in Gl(2n,\R)$.
But this is immediate: for any $X\in As(2n,\R)$ we can choose a $Y\in Mat(2n,\R)$
such that $X=Y-\hat Y$, then $M_Y:=\varphi^{-1}\circ h_\wedge^{-1}\circ Y$
satisfies 
\begin{eqnarray*}
dc(\varphi)(M_Y)&= &\hat Y\circ \hat h_\wedge^{-1}\circ h_\wedge+ Y
=\hat Y\circ(-\Eins) +Y= X.
\end{eqnarray*}

So (after passing to a possibly finer covering,
again called $U_i$)
we find 
lifts in
$$
\xymatrix{
&&Gl(2n,\R)\ar[d]\\
\overline U_i\ar[rr]_-{h_{\hat\omega_i\wedge\omega_i}}
\ar@{..>}[rru]^{\exists\varphi_i}
&&{Gl(2n,\R)\cap As(2n,\R)}
}
$$
and so 
$
h_{(\varphi_i)^*\wedge}=\hat\varphi_i\circ h_\wedge\circ\varphi_i 
= h_{\hat\omega_i\wedge\omega_i}
$
which means that $[(\varphi_i)^*\wedge]=[\hat\omega_i\wedge\omega_i]$.
\end{proof}

If $(A,\alpha)$ is a C$^*$-dynamical system over 
$B$ (i.e. $A$ is a $C_0(B)$-algebra and $\alpha$ is fibre-preserving)  and $(A',\alpha')$ is a C$^*$-dynamical system over $B'$, 
then a morphism $$(f,M):(A,\alpha)\to (A',\alpha')$$ consists of a continous function
$f:B\to B'$  
together with a Morita equivalence $_AM_{f^*A'}$
over $B$. If $(f,M):(A,\alpha)\to (A',\alpha')$ and
$(f',M'):(A',\alpha')\to (A'',\alpha'')$
are two morphisms, then the composition
is $(f\circ f', M\otimes_{f^*A'}f^*M')$ which is only associative up to isomorphism. 
So define a 2-morphism between two morphisms $(g,L)$ and $(g,L')$ to 
be an isomorphism $L\cong L'$ over $B$,
and we have just defined the 2-category {of} 
{\em C$^*$-dynamical systems over spaces}.
Apparently, NC pairs and polarisable NC pairs form 2-subcategories.

\begin{thm}
\label{THEMAINTHEOREM}
The duality functor ${\_\rtimes G}$ 
defined on all C$^*$-dynamical systems over spaces 
with group $G$
restricts to a duality of polarisable NC pairs:
$$
\xymatrix{
{\begin{array}{c}
\textrm{C$^*$-Dynamical Systems}\\
\textrm{{over} Spaces with Group $G$}
\end{array}}
\ar[r]^-\sim
&
{\begin{array}{c}
\text{C$^*$-Dynamical Systems}\\
\text{{over} Spaces with Group $\widehat G$}
\end{array}}
\\
\ar@{}[u]|{{\resizebox{0.4cm}{!}{$\cup$} }}
\textrm{NC Pairs}
&
\ar@{}[u]|{{\resizebox{0.4cm}{!}{$\cup$} }}
\textrm{NC Dual Pairs}
\\
\ar@{}[u]|{{\resizebox{0.4cm}{!}{$\cup$} }}
{\begin{array}{c}
\textrm{Polarisable}\\
\textrm{NC Pairs}
\end{array}}
\ar[r]^-\sim&
\ar@{}[u]|{{\resizebox{0.4cm}{!}{$\cup$} }}
{\begin{array}{c}
\textrm{Polarisable}\\
\textrm{NC Dual Pairs}
\end{array}}
}.
$$

Moreover, if $(A,\alpha)$  is a polarisable
NC pair  with polarisable atlas $U_i,\hat\omega_i,\omega_i$,
then 
$U_i, \omega_i\rtimes\hat\omega_i,
\hat\omega_i\overline{\rtimes}\omega_i$ 
(defined point-wise as in Definition \ref{DefiOfTheDualMoped})
is a polarisable atlas for the dual. 

\end{thm}

\begin{proof}

Let $(A,\alpha)$ be a polarisable NC pair.
We have to show first that $A\rtimes_\alpha G$ is 
locally $\omega$-trivial and that the dual $\widehat G$-action 
on it is locally transverse in the sense of Definition \ref{DefiNCPair} (1).
The proof is organised in 8 steps. 
\\\\
{\bf Step 0:}
Fix a polarisable atlas $U_i,\hat\omega_i,\omega_i$ 
with polarisations $\varphi_i$, i.e.
$\hat\omega_i\wedge\omega_i \sim
(\varphi_i)^*\wedge$.
\\\\
{\bf Step 1:}
We have an isomorphism
\begin{eqnarray*}
C_0(\overline{U_i})\rtimes_{\hat\omega_i\wedge\omega_i} (\widehat G\times G)
&\stackrel{\cong}{\longrightarrow} 
&
C_0(\overline{U_i})\rtimes_{\wedge } (\widehat G\times G)\\
f&\mapsto& c_i\cdot f\circ \varphi_i^{-1},
\end{eqnarray*}
where $c_i:\widehat G\times G\to C(\overline U_i,\UU)$ 
relates $\hat\omega_i\wedge\omega_i \sim
(\varphi_i)^*\wedge$.
(Note: it is not difficult to show that $\varphi_i$ is measure-preserving. 
If it weren't, there would occur a positive factor in this isomorphism 
to rescale the Haar measure.) 
The dual action  of $G\times \widehat G$ is transformed
into the dual action pre-composed with $\hat\varphi_i^{-1}$.
\\\\
{\bf Step 2:} 
We have to make a version of 
the Takai duality isomorphism explicit:
\begin{eqnarray*}
T: C_0(\overline{U_i})\rtimes_{\wedge } (\widehat G\times G)
\stackrel{\cong}{\longrightarrow}
\K(L^2(\widehat G))\otimes C_0(\overline{U_i})
\end{eqnarray*}
which is given for $f\in L^1(\widehat G\times G,C_0(\overline U_i))$
by 
$$
\tilde f(\chi,\chi'):= 
 \int_G f(\chi-\chi',g)\lk g,\chi'\rk\ dg \in C_0(\overline U_i),
$$
where 
$\tilde f$ is the integral kernel for $Tf$, i.e. for 
$\xi\in L^2(\widehat G)$ and $u\in \overline U_i$ we have
$$
(Tf)|_u(\xi)(\chi)=\int_{\widehat G}\tilde f(\chi,\chi')|_u\ \xi(\chi')\ d\chi'.
$$
In fact, it is a lengthy but straight forward calculation
that $T$ defines a $^*$-homomor\-phism which is injective,
because $T$ is composed  by injective transformations
such as the Fourier transformation. The image of
$L^1(\widehat G\times G,C_0(\overline U_i))$ 
is dense in $C_0(\overline{U_i})\otimes \K(L^2(\widehat G))$
and thus $T$ is an isomorphism.
\\\\
{\bf Step 3:}
The dual action of $G\times\widehat G$ on 
$C_0(\overline U_i)\rtimes_{\wedge } (\widehat G\times G)$
transforms under $T$ into\footnote{
We denote  multiplication operators on $L^2$
by the same symbol as their defining functions, e.g. 
$\lk g_0,\_\rk: \xi\mapsto \lk g_0,\_\rk\xi$, for $\xi\in L^2(\widehat G)$. 
} 
\begin{eqnarray}
\label{EqTheNewFreedomAct}
\hat\alpha_i(g,\chi):= T \circ \lk(g,\chi),(\_,\_)\rk\circ T^{-1}=
 {\rm Ad}\Big(\lk g,\_\rk R(\chi)\Big)\otimes \id_{C_0(\overline U_i)},
\end{eqnarray}
where $R:\widehat G\to \U(L^2(\widehat G))$ is the 
right regular representation.
$\hat\alpha_i$ is Morita equivalent 
to the action 
$\hat\alpha_i\otimes \id_{\K(L^2( G))}$
on $ \K(L^2(\widehat G))\otimes C_0(\overline U_i)
\otimes \K(L^2( G))$, and this action is exterior equivalent to 
$\id_{\K(L^2(\widehat G))}\otimes \hat\alpha'_i$
with
$$
\hat\alpha'_i(g,\chi):= 
\id_{C_0(\overline U_i)}\otimes {\rm Ad}\Big( L(g)\lk \_,\chi\rk \Big),
$$
where $L:\widehat G\to \U(L^2(\widehat G))$ is the 
left regular representation.
The cocycle which implements the exterior equivalence
is given by 
$$
v_{(g,\chi)}:=\lk g,\_\rk R(\chi)\otimes\id_{C_0(\overline U_i)}\otimes \lk\_,\chi\rk^{-1} L(g)^{-1}.
$$
The important thing to notice here is that
 $\hat\alpha'_i$ is in the image of 
$\xi_{\overline U_i,G\times\widehat G}$
from  (\ref{EqDefiningXi}):
It is the image 
of  $\vee^{-1}\in Z^2(G\times\widehat G,C(\overline U_i,\UU))$
which reads explicitly 
\begin{eqnarray*}
(g,\chi)\vee^{-1} (h,\psi)=\lk h,\chi\rk^{-1}.
\end{eqnarray*}
As a consequence a crossed product with the
action $\hat\alpha_i$ (or $\hat\alpha'_i$) 
precomposed with $\hat\varphi_i^{-1}$
is Morita equivalent to 
a twisted group algebra 
with cocycle $(\hat\varphi_i)_*\vee^{-1}$.
\\\\
{\bf Step 4:}
Note that for the Heisenberg cocycle $\wedge$
on $\widehat G\times G$ we have 
$\lk h_\wedge(\chi,g),(\psi,h)\rk =
\lk(g,-\chi),(\psi,h)\rk$ which means
that 
$h_\wedge$ is the map $(\chi,g)\mapsto (g,-\chi)$
and so 
$$
\big((h_\wedge)_*\wedge\big) ((g,\chi),(h,\psi))=
\lk g,-\psi\rk
\sim 
(g,\chi)\vee(h,\psi).
$$
This relation holds in $Z^2(G\times \widehat G,\UU)$
but also in $Z^2(G\times \widehat G,C(\overline U_i,\UU))$,
where we consider the above cocycles as constant in the
fibres.
We can use this intermediate step to 
compute the cocycle of $\hat\alpha_i'$ 
from above: 
\begin{eqnarray*}
\Big( (\hat\varphi_i)_*\vee\Big)^{-1}
&\sim&
\Big( (\hat\varphi_i)_* (h_\wedge)_*\wedge\Big)^{-1}\\
&=&\Big( (\hat\varphi_i)_* (h_\wedge)_* (\varphi_i)_*(\varphi_i)^*\wedge\Big)^{-1}\\
&=&\Big( (\hat\varphi_i \circ h_\wedge\circ \varphi_i)_*\ 
\big((\varphi_i)^*\wedge\big)\Big)^{-1}\\
&=&\Big( \big( h_{(\varphi_i)^*\wedge}\big)_*\ 
\big((\varphi_i)^*\wedge\big)\Big)^{-1}\\
&\sim&\Big( \big( h_{\hat\omega_i\wedge\omega_i}\big)_*\ 
(\hat\omega_i\wedge\omega_i)\Big)^{-1}.
\end{eqnarray*}
So apart from an interchange of $G$ and $\widehat G$, 
this is exactly the formula we analysed in 
Lemma \ref{LemComputationOfWidehatOmega}.
\\\\
{\bf Step 5:}
By the continuity and openness of $\varphi_i$, the equality 
$$
\hat\varphi_i(u) \circ h_\wedge\circ\varphi_i(u)
= h_{\hat\omega_i\wedge\omega_i|_u} 
=\left(
\begin{matrix}
h_{\hat\omega_i|_u}& \id_G\\
-\id_{\widehat G}& h_{\omega_i|_u}
\end{matrix}
\right)
$$
shows that $(u,g)\mapsto (u,h_{\omega_i|_u}(g))$ is continuous and open, 
so
\begin{eqnarray*}
\phi_i:\overline U_i&\to& {\rm Aut}(G)\\
u&\mapsto& (g\mapsto g+h_{\hat\omega_i|_u}( h_{\omega_i|_u}(g)))
\end{eqnarray*}
is continuous for the Bracconier topology on ${\rm Aut}(G)$.
We already remarked in (\ref{EqBihomomor})
that the involved cocycles are all bihomomorphisms,
so  as explained in  Remark \ref{RemImportantNotice}
we have
 \begin{eqnarray*}
\Big( (\hat\varphi_i)_*\vee\Big)^{-1}
\sim \omega_i\overline\rtimes\hat\omega_i \cdot 
\hat\omega_i\overline\rtimes\omega_i
\cdot (\phi_i\times\id_{\widehat G})_*\vee^{-1}.
\end{eqnarray*}
This relation holds in $\in Z^2(G\times\widehat G,C(\overline U_i,\UU))$,
i.e. this is a local statement rather than just a point-wise statement.
\\\\ 
{\bf Step 6:} 
We can now locally compute the dual  
\begin{eqnarray*}
(A\rtimes_\alpha G)|_{\overline U_i} 
&\sim&
C_0(\overline U_i)\rtimes_{\hat\omega_i\wedge\omega_i} (N^\perp\times G)\\
&\sim&
C_0(\overline U_i)\rtimes_{\hat\omega_i\wedge\omega_i} (\widehat G\times G)\rtimes_{\rm dual} (N\times 0)\\
&\cong&
C_0(\overline U_i)\rtimes_{\varphi_i^*\wedge} (\widehat G\times G)\rtimes_{\rm dual} (N\times 0)\\
&\cong&
C_0(\overline U_i)\rtimes_{\wedge} (\widehat G\times G)
\rtimes_{\rm dual\circ \hat \varphi_i^{-1}} (N\times 0)\\
&\cong&
 \K(L^2(\widehat G))
\otimes C_0(\overline U_i)
\rtimes_{\alpha_i\circ \hat \varphi_i^{-1}} (N\times 0)\\
&\sim&
C_0(\overline U_i)
\rtimes_{(\hat \varphi_i)_*\vee^{-1}} (N\times 0)
\\
&\cong&
C_0(\overline U_i)
\rtimes_{\omega_i\overline\rtimes\hat\omega_i} (N\times 0)\\
&\cong&
C_0(\overline U_i)
\rtimes_{\omega_i\rtimes\hat\omega_i} (N\times 0),
\end{eqnarray*}
where the last isomorphism is pullback by
$\phi_i$. 
This implies the local $\omega$-triviality of
$A\rtimes_\alpha G$
and shows that $\omega\rtimes\hat\omega_i$ gives an atlas.
\\\\
{\bf Step 7:}
Using these Morita equivalences,
we can similarly compute the dual action of $\widehat G$ on  
$(A\rtimes_\alpha G)|_{\overline U_i}$
(note that all equivalences are $G$-equivariant):
\begin{eqnarray*}
(A\rtimes_\alpha G)|_{\overline U_i}\rtimes_{\rm dual} \widehat G 
&\sim&
C_0(\overline U_i)\rtimes_{\hat\omega_i\wedge\omega_i} (N^\perp\times G)
\rtimes_{\rm dual}\widehat G\\
&\sim&
C_0(\overline U_i)\rtimes_{\hat\omega_i\wedge\omega_i} (\widehat G\times G)\rtimes_{\rm dual} (N\times \widehat G)\\
&\sim&
C_0(\overline U_i)
\rtimes_{(\hat \varphi_i)_*\vee^{-1}} (N\times \widehat G)\\
&\cong&
C_0(\overline U_i)
\rtimes_{
\omega_i\overline\rtimes\hat\omega_i\cdot 
(\phi_i\times\id_{\widehat G})_*\vee^{-1}\cdot
\hat\omega_i\overline\rtimes\omega_i 
} 
(N\times \widehat G)\\
&\cong&
C_0(\overline U_i)
\rtimes_{
\omega_i\rtimes\hat\omega_i\cdot 
\vee^{-1}\cdot
\hat\omega_i\overline\rtimes\omega_i 
} 
(N\times \widehat G)\\
&\cong&
C_0(\overline U_i)
\rtimes_{
\omega_i\rtimes\hat\omega_i\cdot 
\vee\cdot
\hat\omega_i\overline\rtimes\omega_i 
} 
(N\times \widehat G),
\end{eqnarray*}
where the last 
step is pullback along 
the inversion $\ominus$ in $N$ 
together with the isomorphism 
induced by the cochain that
relates
$\ominus^*(\omega_i\rtimes \hat\omega_i)
\sim  
\omega_i\rtimes \hat\omega_i$.
This again requires 
the computation of Remark \ref{RemImportantNotice}.
As $\omega_i\rtimes\hat\omega_i$ and 
$\hat\omega_i\overline\rtimes\omega_i$
are point-wise transverse, this completes the proof.
\end{proof}

By Lemma \ref{LemTheDualMackeyObstructions},
we can associate to a given NC pair $(A,\alpha)$ 
an obstruction function 
\begin{eqnarray}\label{EqTheFctDefinedByAction}
\theta: B\to H^2(N,\UU)
\end{eqnarray}
which at each point of $B$ has the Mackey obstruction
of the local and restricted $N$-action as its value.
This function is continuous. In fact, it is locally 
continuous as 
it is given by point-wise evaluation\footnote{
Here evaluation is an isomorphism because $N$ is 
discrete, so the topology of almost everywhere point-wise
convergence on $Z^2(N,\UU)$ coincides with the compact open topology
for which we can apply the exponential law for 
locally compact Hausdorff spaces $X,Y,Z$:
$C(X\times Y,Z)\cong C(X,C(Y,Z))$}
$$
Z^2(N,C(\overline U_i,\UU))\cong C(\overline U_i,Z^2(N,\UU))
\to C(\overline U_i,H^2(N,\UU)).
$$
Similarly, to a NC dual pair we associate
an obstruction function
\begin{eqnarray}
\label{EqTheFctDefinedByDualAction}
\hat\theta: B\to H^2(N^\perp,\UU).
\end{eqnarray}
Using these functions we can 
talk about the commutative subtheory 
inside polarisable NC pairs:

\begin{defi}
A polarisable NC pair $(A,\alpha)$
 is called {\it (point-wise) commutative} if the obstruction function
$\hat\theta$ defined by its dual $(A\rtimes_\alpha G,\hat\alpha)$
vanishes.
\end{defi}

If we restrict to 
the groups 
$G=\R^n,N=\Z^n$ 
the next proposition shows that 
this notion of  commutativity
reproduces the familiar objects 
from classical T-duality.

\begin{prop}\label{PropPrincipalityIsImplied}
Let $G=\R^n,N=\Z^n$, and let 
$(A,\alpha)$ be a (polarisable) 
 NC pair.
Then  $(A,\alpha)$ is point-wise commutative
if and only if
there is a locally trivial principal $G/N$-bundle $E\to B$ and a locally trivial
bundle of compact operators $F\to E$ such that
\begin{eqnarray}
\label{EqCommCase}
A\cong \Gamma_0(E,F).
\end{eqnarray}
Moreover, $F\to E$ is trivialisable over (a neighbourhood of) 
all fibres $E|_b, b\in B$,
and there is a $G$-action on $F$ that covers the $G/N$
action on $E$ such that the isomorphism in (\ref{EqCommCase})
can be chosen to be $G$-equivariant.
\end{prop}

\begin{proof}
It is clear that if an isomorphism 
(\ref{EqCommCase}) exists, then $\hat\theta=0$.

Conversely, if $\hat\theta(b)=0$, then
by Lemma \ref{LemAboutThetaTransverse}
$\hat\theta(b)^\tri=\{0\}$ which determines the 
class of each fibre to be 0.
So we have locally 
$A|_{\overline U_i}\sim C_0(\overline U_i)\rtimes_{\hat\omega_i} N^\perp$, where
the cocycle $\hat\omega_i:Z^2(\widehat G,C(\overline U_i,\UU))$ is such that
its restrictions 
$$\hat\eta_i:=\hat\omega_i|_{N^\perp}
\in Z^2(N^\perp,C(\overline U_i,\UU))\cong C(\overline U_i,Z^2(N^\perp,\UU))
$$ 
are point-wise in the image 
of the
boundary operator $d:C^{1}(N^\perp,\UU) \to C^{2}(N^\perp,\UU)$:
$\hat\eta_i|_b\in B^2(N^\perp,\UU)$.

The kernel of $d$ is ${\rm Hom}(N^\perp,\UU)\cong G/N=\T^n$ which 
is a Lie group.
{So just as in the proof of \cite[Theorem 2.1]{Rosi}} 
we apply the Palais 
cross-section theorem \cite[4.1]{Pa} which implies 
that $d: C^{1}(N^\perp,\UU)\to B^2(N^\perp,\UU)$
is a locally trivial $G/N$-bundle.
Then, after passing to a possibly finer covering of $B$ which we again call $U_i$,
we have lifts $\hat \nu_i$ in
\begin{eqnarray*}
\xymatrix{
& C^1(N^\perp,\UU)\ar@{->>}[d]^d\\
\overline U_i\ar[r]_-{\hat\eta_i}\ar[ur]^-{\hat\nu_i}& B^2(N^\perp,\UU)  
}
\end{eqnarray*}
which implies the existence of  local $G$-equivariant isomorphisms
\begin{eqnarray*}
A|_{\overline U_i}\cong \K\otimes C_0(\overline U_i)\rtimes_{d\circ\hat\nu_i} N^\perp\cong
\K\otimes C_0(\overline U_i)\rtimes_{1} N^\perp\cong
C_0(U_i\times G/N,\K).
\end{eqnarray*}
Now, consider the transition 
on the overlap $U_{ij}:=U_i\cap U_j$ of two charts
\begin{eqnarray*}
\xymatrix{
C_0(\overline U_{j}\times G/N,\K)
\ar@{->>}[d]& 
C_0(\overline U_{i}\times G/N,\K)
\ar@{->>}[d]\\
C_0(\overline U_{ji}\times G/N,\K)\ar[r]^\cong_{\varphi_{ij}}
&
C_0(\overline U_{ij}\times G/N,\K). 
}
\end{eqnarray*}
It is $G$-equivariant, and it induces a $G$-equivariant automorphism 
on the spectrum $\gamma_{ji}: \overline U_{ij}\times G/N\cong \overline U_{ji}\times G/N$.
By equivariance $\gamma_{ji}$ is of the form
$\gamma_{ji}(u,z)=(u,g_{ji}(u)+z)$ for some 
$g_{ji}:\overline U_{ij}\to G/N$. Clearly, the $g_{ji}$ define a $\check {\rm C}$ech
cocycle and thus a principal $G/N$-bundle $E\to B$.
Precomposition of $\varphi_{ij}$ with $(\gamma_{ji}^{-1})_*$ yields then 
a spectrum fixing automorphism of $C_0(\overline U_{ij}\times G/N,\K)$, i.e.
a function $\zeta_{ji}:\overline U_{ij}\times G/N\to \PU$ such that
$\varphi_{ij}(f)(u,z)= \zeta_{ji}^{-1}(u,z)[f(u,g_{ji}(u)+z)]$
for $f\in C_0(\overline U_{ji}\times G/N,\K)\twoheadleftarrow 
C_0(\overline U_j\times G/N,\K)$.
The $\zeta_{ij}$ satisfy the twisted 
$\check {\rm C}$ech identity
\begin{eqnarray*}
\zeta_{kj}(u,g_{ji}(u)+z) \ \zeta_{ji}(u,z)=\zeta_{ki}(u,z)
\end{eqnarray*}
and hence define a bundle of compact operators $F\to E$.
It is then clear from the construction that $F$ is trivialisable over
a neighbourhood of each fibre $E|_b\subset E$, and 
$F$ carries a $G$-action which covers the principal $G/N$-action on $E$.

For any $a\in A$ the quotient  maps 
$A\twoheadrightarrow C_0(\overline U_i\times G/N,\K)$
define a compatible family of functions, i.e. 
they define a section of $F\to E$. Hence we get a $G$-equivariant map
\begin{eqnarray*}
A\to \Gamma(E,F).
\end{eqnarray*}
Because this map is $C_0(B)$-linear, one can use a partition 
of unity on $B$ as an approximate identity on $C_0(B)$ 
to show that this map takes values in the sections vanishing at infinity 
only and that this assignment is injective and surjective    
\begin{eqnarray*}
A\stackrel{\cong}{\to} \Gamma_0(E,F)\subset \Gamma(E,F).
\end{eqnarray*}
\end{proof}

\subsection{NC Bundles}\label{NCbundle}
If $(A,\alpha)$ is a NC pair we have seen that it
is not reasonable to ask for the cohomology classes of the fibres
of $A$ rather than to handle this issue with the ambiguity which is measured 
by the set $\hat\theta^\tri$ as defined in (\ref{EqDefiOfThetaPerp}).
So if we want to deal with bundles without actions the following 
definition is appropriate.  

\begin{defi}
A {\it NC bundle} $(A,\hat\theta)$ is a locally $\hat\omega$-trivial
algebra together with a continuous function $\hat\theta:B\to H^2(N^\perp,\UU)$
such that there exists an altas $U_i,\hat\omega_i$ of $A$ which has the 
property that for all $i$ we have 
$[\hat\omega_i|_b]|_{N}\in \hat\theta(b)^\tri$, for all $b\in U_i$.
A NC bundle $(A,\hat\theta)$ is called {\it commutative}
if $\hat\theta=0$.

(The notion of {\it (commutative) NC dual bundles} is defined analogously.)
\end{defi}

These are the objects of the category of NC bundles. The morphisms
$(A,\hat\theta)\to (A',\hat\theta')$
are pairs $(f,M)$ of a continuous function $f:B\to B'$ such that
$\hat\theta =\hat\theta'\circ f$ and a 
$C_0(B)$-linear Morita equivalence ${}_AM_{f^*A'}$.
Here we assume that $A'$ is a $C_0(B')$-algebra.

If $(A,\alpha)$ is any polarisable NC pair, we can use 
the obstructions function  $\hat\theta:B\to H^2(N^\perp,\UU)$
of its dual $(A\rtimes_\alpha G,\hat \alpha)$
to define an assignment on objects 
$(A,\alpha)\mapsto (A,\hat\theta)$
which extends to a functor
\begin{eqnarray*}
\xymatrix{
\Theta: {\begin{array}{c}
\textrm{Polarisable}\\
\textrm{NC Pairs}
\end{array}}
\ar[r]&
{\begin{array}{c}
\textrm{NC Bundles}
\end{array}}
}
\end{eqnarray*}
that gives the {\it underlying NC bundle}
of a polarisable NC pair.  
Similarly, there is a functor 
\begin{eqnarray*}
\xymatrix{
\widehat\Theta: {\begin{array}{c}
\textrm{Polarisable}\\
\textrm{NC Dual Pairs}
\end{array}}
\ar[r]&
{\begin{array}{c}
\textrm{NC Dual}\\
\textrm{Bundles}
\end{array}}
}
\end{eqnarray*}
that gives the underlying NC dual bundle of 
a polarisable NC dual pair.
An {\it extension} of a NC (dual) bundle $(A,\hat\theta)$
is a polarisable NC (dual) pair $(A,\alpha)$ such that
$\Theta(A,\alpha)=(A,\hat\theta)$, and a NC (dual) bundle
is called {\it dualisable} if it has an extension.
Now, if $(\hat A,\theta)$ is a dualisable NC dual bundle,
one might ask the question whether it has an extension 
that is the dual of a commutative polarisable NC  pair.
In other words, we ask whether 
 the functor
\begin{eqnarray*}
\xymatrix{
\Xi:{\begin{array}{c}
\textrm{Polarisable}\\
\textrm{NC Pairs}
\end{array}}
\ar[r]^-\sim&
{\begin{array}{c}
\textrm{Polarisable}\\
\textrm{NC Dual Pairs}
\end{array}}
\ar[r]^-{\widehat\Theta}
&
{\begin{array}{c}
\textrm{NC Dual}\\
\textrm{Bundles}
\end{array}}
}
\end{eqnarray*}
and its restriction 
\begin{eqnarray*}
\xymatrix{
\Xi|_{\textrm{com}}:{\begin{array}{c}
\textrm{Commutative}\\
\textrm{Polarisable}\\
\textrm{NC Pairs}
\end{array}}
\ar@{}[r]|-{\resizebox!{0.34cm}{$\subset$}}
&
{\begin{array}{c}
\textrm{Polarisable}\\
\textrm{NC Pairs}
\end{array}}
\ar[r]^-\sim&
{\begin{array}{c}
\textrm{Polarisable}\\
\textrm{NC Dual Pairs}
\end{array}}
\ar[r]^-{\widehat\Theta}
&
{\begin{array}{c}
\textrm{NC Dual}\\
\textrm{Bundles}
\end{array}}
}
\end{eqnarray*}
have the same (essential) image.
If the answer to this question were yes,
then the theory we developed wouldn't
be significantly richer than the 
classical commutative and semi-comutative 
theory of Mathai and Rosenberg, 
because all NC bundles could be understood
as crossed products of classical bundles (with all their possible $G$-actions).
However, this is not the case:

\begin{prop}
\label{PropTheImages}
Let $G:=\R^2$ and $N:=\Z^2$. Then the essential images 
of the functors $\Xi$ and $\Xi|_\text{com}$
do not coincide.
\end{prop}
\begin{proof}
In section \ref{SecTheTwistedHeis} below
we give an example of a NC bundle that 
cannot be obtained as the dual of a commutative 
polarisable NC pair.
\end{proof}

\section{Example: The Heisenberg Bundle and Its Relatives} 
\label{SecExamplesOverTheCircle}

In this section we always let $G:=\R^2$ and $N:=\Z^2$ and we use the canonical 
identifications $\widehat{\R^2}\cong\R^2$ and $(\Z^2)^\perp\cong \Z^2$.
We will typically denote elements in $\T=\R/\Z$ or $\T^2=\R^2/\Z^2$ by $x,y,z$ 
or by $\dot{g}$ if $g$ is in $\R$ or $\R^2$.

By $S^1$ we denote the unit interval $[0,1]$ with glued end-points $0\sim 1$.
Sometimes it is convenient to have have one of the endpoints, say 1,  
thickened to a whole (non-empty) interval $\fat:=[1_-,1_+]$, i.e. we then consider the 
space
$S^{\fat}:=([0,1]\sqcup \fat)/_{\!\sim}$,
where $1\sim 1_-$ and $0\sim 1_+$ which, of course, is homeomorphic to $S^1$.
For $t\in \R$ we denote by $\omega_t$ the 2-cocycle on $\R^2$ given by
$
\omega_t(g,h):=\exp(2\pi i t g_2 h_1).
$
We will use the notation $\lambda^t$, $t\in\R$, for the actions
$\R^2\to \Aut(\K)$ which are given by
$
{\lambda^t}_g:={\rm Ad}(L_{-t}(g))
$
for the left regular $\omega_{-t}$-representation
$$ 
(L_{-t}(g)(\xi))(h):=\omega_{-t}(g,h-g) \xi(h-g),\quad \xi\in L^2(\R^2).
$$
The actions $\lambda^t$ have Mackey obstruction $+t\in\R\cong H^2(\R^2,\UU)$.

\subsection{The Heisenberg Bundle}

Consider the function $\omega_{S^1}:S^1\to Z^2(N,\UU)$ given by
$
\omega_{S^1}(\dot s):=\omega_s|_{N\times N}
$
which is well-defined as $\omega_n|_{N\times N}=1$, for $n\in\Z$.
The Heisenberg bundle is the C*-algebra over $S^1$ given by
$$
\hat A_0:=\K\otimes (C(S^1)\rtimes_{\omega_{S^1}} N).
$$
Its fibres are the stable non-commutative tori 
$\hat A|_{\dot s} \cong\K\otimes(\C\rtimes_{\omega_s|_{N\times N}}N)$.
We equip it with the canonical $\widehat G$-action $\hat\alpha_0:=\id\otimes\inf$.
The following proposition is immediate.
\begin{prop}
$(\hat A_0,\hat\alpha_0)$ is a (polarisable) NC dual pair.
\end{prop}
The obstruction function $\theta_0:S^1\to H^2(\Z^2,\UU)$ 
defined by  $(\hat A_0,\hat\alpha_0)$
 is just the 
 composition of the canonical identifications $S^1\cong\T\cong H^2(\Z^2,\UU)$.

Let us construct a (commutative) NC pair which will turn out to be 
the dual of the Heisenberg bundle.
Consider the trivial principal $\T^2$-bundle $E_0:=S^1\times\T^2\to S^1$.
Its total space $E_0$ is a compact orientable three manifold so its
3rd ($\check{\rm C}$ech) cohmology
is $\check H^3(E_0,\Z)\cong \Z$.
Let $F_1\to E$ be a locally trivial $\K$-bundle representing
the canonical generator of 3rd  cohomology:
\begin{eqnarray*}
\check H^3(E_0,\Z)&\stackrel{\cong}{\longrightarrow}& \Z.\\
{} [F_1]&\mapsto&{1}
\end{eqnarray*}
Let $\Gamma(E_0,F_1)$ be the section C*-algebra of $F_1\to E$.
We construct an action on this algebra 
which covers the principal $\T^2$-action on $E_0={\rm Prim}(\Gamma(E_0,F_1))$.
To do so, observe first that one can describe the
bundle $F_1$ slightly different.
Up to isomorphism any 
$\K$-bundle over $E_0$ can be obtained
from a function $\T^2\to \Aut(\K)$ which is used to glue 
the two boundary parts of the trivial $\K$-bundle   
\begin{eqnarray}
\label{EqTheTrivialKbundleOverThe}
[0,1]\times\T^2\times\K\to [0,1]\times \T^2
\end{eqnarray}
to each other. In particular, the bundle $F_1$ 
is obtained by using a classifying 
map $\T^2\to \mathcal B\UU=\PU(\HH)=\Aut(\K)$
of the canonical $\UU$-bundle over $\T^2$,
i.e. a function whose class in 2nd  
cohomology $\check H^2(\T^2,\Z)\cong\Z$ corresponds to 1. 

Let us identify such a function $\T^2\to\Aut(\K)$ by the following construction.
Choose $\K=\K(L^2(\R^2))$
and define a $C([0,1])$-linear 
action $\beta$ on  the sections $C([0,1]\times\T^2,\K)$ of 
(\ref{EqTheTrivialKbundleOverThe})
by
$\beta_g(f)(t, x)={\lambda^t}_g(f(t, x-\dot g))$.
For $t=1$ the cocycle involved in $\lambda^1$,
$
\omega_{-1}:(g,h)\mapsto \exp(2\pi i (-1) g_2 h_1),
$
becomes trivial when restricted to $N=\Z^2$. 
To continue we need a lemma:
\begin{lem}
\label{LemCommutes}
The canonical isomorphism  
$H^2(\R^2,L^\infty(\R^2/\Z^2,\UU)\to H^2(\Z^2,\UU)$
of (\cite[Thm. 6]{Moore3}) makes the diagramm
\begin{eqnarray*}
\xymatrix{
&H^2(\R^2,\UU)\ar[d]\ar[dr]&\\
&H^2(\Z^2,\UU)&H^2(\R^2,L^\infty(\R^2/\Z^2,\UU))\ar[l]_-\cong
}
\end{eqnarray*}
commute, where the vertical map is restriction and 
the diagonal map is induced by the inclusion of coefficients.
\end{lem}

\begin{proof}
See Appendix \ref{AppCommutes}.
\end{proof}

Because of this lemma
there is a Borel function\footnote{We always consider $L^\infty(X,\UU)$
as a subspace of $\U(L^2(X))$ by associating to an $L^\infty$-function the
multiplication operator that multiplies  $L^2$-functions point-wise 
with the given $L^\infty$-function. 
}
$c:\R^2\to L^\infty(\T^2,\UU)\subset \U(L^2(\T^2))$ such that
$\omega_{-1}=dc$, i.e.
\begin{eqnarray}
\label{EqTheBoundaryOf1}
\omega_{-1}(g,h)=c(h)(z-\dot g)\ c(g+h)(z)^{-1}\ c(g)(z),
\end{eqnarray}
for almost all $z\in\T^2$. Note at this point that 
(\ref{EqTheBoundaryOf1}) implies that
the restriction $c|_{\Z^2}:\Z^2\to L^\infty(\T^2,\UU)$
is a homomorphism, i.e. there is a measure-one set
$S\subset \T^2$ such that $n\mapsto c(n)(z)$ is in
$\widehat \Z^2$, for all $z\in S$.
The function $c$ determines a function 
$\tilde c\in L^\infty(\R^2\times\T^2,\UU)\subset\U(L^2(\R^2\times\T^2))$
which for each $g$ satisfies
\begin{eqnarray}
\label{EqTheBoundaryTransp}
\omega_{-1}(g,h)=\tilde c(h,z-\dot g)\ \tilde c(g+h,z)^{-1}\ c(g)(z),
\end{eqnarray}
for almost all $(h,z)\in \R^2\times\T^2$.
(By choosing a Borel representative of $\tilde c$) we consider for $z\in \T^2$ 
the unitary multiplication 
operator $\tilde c(\_,z)\in \U(L^2(\R^2))$ which sends a funtion $\xi \in L^2(\R^2)$
to $h\mapsto \tilde c(h,z)\xi(h)$. These operators define a Borel function  
$z\mapsto\eta_0(z):={\rm Ad}(\tilde c(\_,z))\in \PU(L^2(\R^2)).$

\begin{lem}
\label{LemEinEtaZumVerkleben}
\begin{enumerate}
\item There exists a continuous function
\begin{eqnarray*}
\eta:\T^2&\to& \PU(L^2(\R^2))
\end{eqnarray*}
which agrees almost everywhere with $\eta_0$.
\item
The class of $\eta$ in $\check H^2(\T^2,\Z)\cong \Z$ is 
the generator $+1$, i.e. $\eta$ is a classifying map for the 
canonical line bundle over $\T^2$.
\item
The equality
$$
\eta(z)\ \lambda^{k+1}_g=\lambda^{k}_g\ \eta(z-\dot g) \in \PU(L^2(\R^2))
$$
holds, for all $k\in\Z$, $g\in\R^2, z\in \T^2$. 
\end{enumerate}
\end{lem}

\begin{proof}

1. Consider the countable dense subgroup $\mathbb Q^2\subset \R^2$. 
 By (\ref{EqTheBoundaryTransp}) there is for each $g\in\mathbb Q^2$ a set $S_g\subset\T^2$ of measure one such that 
\begin{eqnarray*}
\omega_{-1}(g,\_)=\tilde c(\_,z-\dot g_k)\ \tilde c(g_k+\_,z)^{-1}\ c(g_k)(z)\in U(L^2(\R^2))
\end{eqnarray*}
holds for all $z\in S_g$.
Choose some $z_0\in S\cap \bigcap_{g\in\mathbb Q^2} S_g$. Then
\begin{eqnarray*}
\omega_{-1}(g,\_)=\tilde c(\_,z_0-\dot g)\ \tilde c(g+\_,z_0)^{-1}\ c(g)(z_0)\in U(L^2(\R^2))
\end{eqnarray*}
holds for all $g\in\mathbb Q^2$.
The map 
\begin{eqnarray*}
g&\mapsto& {\rm Ad}(\tilde c(\_,z_0-\dot g))\\
&=&{\rm Ad}(\omega_{-1}(g,\_) \tilde c(g+\_,z_0))\\
&=&{\rm Ad}(\omega_{-1}(g,\_) L_0(-g) \tilde c(\_,z_0)L_0(g))
\end{eqnarray*}
is clearly continuous on $\mathbb Q^2$ (as one can see in line three) 
and clearly factors over the quotient $\mathbb Q^2/\Z^2$  (as one can see in line one).
Then define $\eta_{z_0}:\T^2\to \PU(L^2(\R^2))$  be its continuous extension 
(which is just given by the formula in line three, for all $g\in \R^2$),
and let $\eta(z):=\eta_{z_0}(z_0-\dot g)$.
By construction it's clear that $\eta$ satisfies $(1)$ of the lemma.

2. The $\check{\rm C}$ech classes of $\eta$ and of  $z\mapsto \eta(z-z_0)$ agree.
We compute the $\check{\rm C}$ech class of the latter.
Choose continuous, local sections $\sigma_k:V_k\to\R^2$
of the quotient map $\R^2\to\T^2$ such that
the open domains $V_k\subset \T^2$ cover $\T^2$.
Then by line two from above
$$
\eta_k(z):= \omega_{-1}(-\sigma_k(z),\_)\ \tilde c(-\sigma_k(z)+\_,z_0)\in\U(L^2(G)) 
$$
are continuous 
unitary local lifts of $z\mapsto\eta(z-z_0)$, so its class in 
$\check H^1(\T^2,\UU)\cong \check H^2(\T^2,\Z)$ 
is given by the cocycle
\begin{eqnarray*}
\eta_{kl}(z)&:=&\eta_k(z)\eta_l(z)^{-1}\\
&=&\omega_{-1}(-\sigma_k(z),\_)\ \tilde c(-\sigma_k(z)+\_,z_0)\
\tilde c(-\sigma_l(z)+\_,z_0)^{-1}\omega_{-1}(-\sigma_l(z),\_)^{-1}\\
&=&\omega_{-1}(-\sigma_k(z)+\sigma_l(z),\_)\ \tilde c(-\sigma_k(z)+\_,z_0)\\
&&\hspace{4.4cm}\cdot \tilde c((-\sigma_k(z)+\_)+(\sigma_k(z)-\sigma_l(z)),z_0)^{-1}\\
&=&\omega_{-1}(-\sigma_k(z)+\sigma_l(z),\_)\ 
\omega_{-1}(\sigma_k(z)-\sigma_l(z),\sigma_k(z)+\_)\\ 
&&\hspace{7cm}\cdot c(\sigma_k(z)-\sigma_l(z))(z_0)^{-1}\\
&=&
\omega_{-1}(\sigma_k(z)-\sigma_l(z),\sigma_k(z))\ 
c((\sigma_k(z)-\sigma_l(z)))(z_0)^{-1} \in \UU.
\end{eqnarray*}
We claim that $c((\sigma_k(z)-\sigma_l(z)))(z_0)^{-1}$
is a coboundary term. In fact, $\sigma_k(z)-\sigma_l(z)$ is in $\Z^2$, and so choose
by surjectivity of $\widehat{\R^2}\to\widehat{\Z^2}$ an
extension $\chi\in \widehat{\R^2}$ 
of the character $n\mapsto c(n)(z_0)$ which  gives
$\chi(\sigma_k(z)) \chi(\sigma_l(z))^{-1}=c((\sigma_k(z)-\sigma_l(z)))(z_0)$.
This indeed is a coboundary term and does not effect the class of $\eta$. 
The expressions
$z=(z_1,z_2)\mapsto\omega_{-1}(\sigma_l(z)-\sigma_k(z),\sigma_k(z))=
\lk \sigma_k(z)_2-\sigma_l(z)_2,z_1\rk$
are transition functions for the 
canonical line bundle on $\T\times \T$
\cite[Sec. 2.6]{Sch1}, i.e.
they give the class of the canonical line bundle.

3. Just multiply both sides of (\ref{EqTheBoundaryOf1}) 
by the left regular representation $L_{-k}(g)$
and apply ${\rm Ad}:\U(L^2(G))\to \PU(L^2(G))$ to both sides.

\end{proof}

Part 3. of this lemma just says for $k=0$ that 
\begin{eqnarray*}
\xymatrix{
C(1\times \T^2,\K)\ar[r]^{\beta|_1(g)}\ar[d]^{\eta_*}&C(1\times\T^2,\K)\ar[d]^{\eta_*}\\
C(0\times \T^2,\K)\ar[r]^{\beta|_0(g)}&C(0\times\T^2,\K)
}
\end{eqnarray*}
commutes, where $\beta|_0,\beta|_1$ are the actions  
on the fibres over $t=0,1$ given by the fibre-wise action $\beta$, 
and $\eta_*(f)(1,z):=\eta(z)(f(0,z))$.
So by construction we have shown the 
following proposition.
\begin{prop}
\label{PropTheDefiDualOfHeisi}
Let $(A_0,\alpha_0)$ be the pullback in the category of  
C*-dynamical systems of the diagram 
\begin{eqnarray*}
\xymatrix{
(A_0,\alpha_0)\ar@{..>}[dd]\ar@{..>}[rr]&&
(C([0,1]\times\T^2,\K),\beta)
\ar[dd]^{i_0^*\times i_1^*}\\
&&\\
(C(1\times\T^2,\K),\beta|_1)\ar[rr]^-{\eta_*\times\id}&&
(C(0\times\T^2,\K)\oplus C(1\times\T^2,\K),\beta|_0\times\beta|_1)
}
\end{eqnarray*}
then there is a canonical isomorphism $A_0\cong\Gamma(E_0,F_1)$.
\end{prop}

It is clear that $(A_0,\alpha_0)$
is a polarisable NC pair with a trivial obstruction function 
$\hat \theta_0=0:S^1\to H^2(\Z,\UU)$, so its dual is commutative:

\begin{prop}
\label{PropDualityForHeisenberg}
$(A_0,\alpha_0)$ and $(\hat A_0,\hat\alpha_0)$
are dual to each other, i.e $A_0\rtimes_{\alpha_0} \R^2$ with its dual 
$\widehat{\R^2}$-action is Morita equivalent to $(\hat A_0,\hat\alpha_0)$.
\end{prop}

\begin{proof}
Consider the function $\omega_{[0,1]}:[0,1]\to Z^2(\Z^2,\UU)$ given by
$
\omega_{[0,1]}(t):=\omega_t|_{\Z^2\times \Z^2},
$
then the Heisenberg bundle $\hat A_0$ with its induced action
$\hat\alpha_0$ is the pullback in  
\begin{eqnarray}
\label{DiagHeisiIsPullback}
\xymatrix{
(\hat A_0,\hat \alpha_0)\ar@{..>}[dd]\ar@{..>}[rr]&&
(\K\otimes (C([0,1])\rtimes_{\omega_{[0,1]}}\Z^2),\id\otimes\inf) 
\ar[dd]^{i_0^*\times i_1^*}\\
&&\\
(C(1\times\T^2,\K),\inf)\ar[rr]^-{\rm diag}&&
(C(0\times\T^2,\K)\oplus C(1\times\T^2,\K),\inf\otimes\inf).
}
\end{eqnarray}
In the diagram of Proposition \ref{PropTheDefiDualOfHeisi} 
we can take crossed product with $\R^2$. 
Because taking crossed products is a continuous functor,
this leads to another pullback diagram.
But this new diagram is stably isomorphic to  
diagram (\ref{DiagHeisiIsPullback}).
\end{proof}

\begin{rem}
The duality stated in Proposition \ref{PropDualityForHeisenberg} has already been
observed in \cite[Section 4]{MR2}.
However, they follow a  different approach which is less explicit 
then the one presented here, and we can use the intermediate 
steps of our approach  for the construction of the twisted Heisenberg bundle
which we do next. 
\end{rem}

\subsection{The Twisted Heisenberg Bundle}
\label{SecTheTwistedHeis}

Let us use the pullback description of the two  NC (dual) pairs from above 
to construct a chimaera out of the two.
Consider the algebra $\K\otimes (C([0,1])\rtimes_{\omega_{[0,1]}}\Z^2)$,
where  $\omega_{[0,1]}$ is as in 
(\ref{DiagHeisiIsPullback}).
This algebra carries the fibre-wise action
$\gamma$ that is given in each fibre by
$\gamma|_t:= {\lambda^{\frac{2}{1+t}}} \otimes \inf$. 
As $t$ moves from 0 to 1, the Mackey obstuction of 
$\lambda^{\frac{2}{1+t}}$ moves from 2 to 1. 
For $k=1$ part 3. of Lemma \ref{LemEinEtaZumVerkleben} 
implies that 
$\eta(z){\lambda^{2}}_g={\lambda^{1}}_g\eta(z-\dot g)$,
i.e.
\begin{eqnarray*}
\xymatrix{
C(1\times \T^2,\K)\ar[r]^{\gamma|_1(g)}\ar[d]^{\eta^*}&C(1\times\T^2,\K)\ar[d]^{\eta^*}\\
C(0\times \T^2,\K)\ar[r]^{\gamma|_0(g)}&C(0\times\T^2,\K)
}
\end{eqnarray*}
commutes, for $\eta^*(f)(1,z):=\eta(z)^{-1}(f(0,z))$.
So naively, what we now could consider is the pullback in
\begin{eqnarray*}
\xymatrix{
(\hat A_1,\hat \alpha_1)\ar@{..>}[dd]\ar@{..>}[rr]&&
(\K\otimes (C([0,1])\rtimes_{\omega_{[0,1]}}\Z^2 ),\gamma)
\ar[dd]^{i_0^*\times i_1^*}\\
&&\\
(C(1\times\T^2,\K),\gamma|_1)\ar[rr]^-{\eta^*\times \id}&&
(C(0\times\T^2,\K)\oplus C(1\times\T^2,\K),\gamma|_0\times \gamma|_1).
}
\end{eqnarray*}
Indeed $\hat A_1$ is a bundle over $S^1$ with fibers 
$\hat A_1|_{\dot s}\cong\K\otimes (\C\rtimes_{\omega_s} N)$. 
However, it fails to be 
$\omega$-trivial around $\dot 0=\dot 1\in S^1$.
We can get around this by thickening the gluing-point to the 
interval $\fat$:
We define the {\sl twisted Heisenberg bundle} $\hat A_{\fat}$ 
together with its action $\hat\alpha_{\fat}$ to be the pullback in 
\begin{eqnarray*}
\xymatrix{
(\hat A_{\fat},\hat \alpha_{\fat})\ar@{..>}[dd]\ar@{..>}[rr]&&
(\K\otimes (C([0,1])\rtimes_{\omega_{[0,1]}}\Z^2 ),\gamma)
\ar[dd]^{i_0^*\times i_1^*}\\
&&\\
(C(\fat\times\T^2,\K),\gamma|_{\fat})\ar[rr]^-{(\eta^*\circ i_{1_+}^*)\times i_{1_-}^*}&&
(C(0\times\T^2,\K)\oplus C(1\times\T^2,\K),\gamma|_0\times \gamma|_1),
}
\end{eqnarray*}
where the action $\gamma|_{\fat}$ is just $\gamma|_1$ in each fibre.

\begin{prop}
\label{PropHeisenbergsCousin}
The twisted Heisenberg bundle
$(\hat A_{\fat},\hat\alpha_{\fat})$ 
is a (polarisable) NC dual pair over
$S^{\fat}$.
\end{prop}

\begin{proof}
Let $1_\pm \in\fat$ be the middle. We consider the open cover of $S^\fat$ 
given by the two open sets $U:=((0,1]\sqcup[1_-,1_\pm))/_{\!\sim}$ and
$V:=([0,\frac{1}{2})\sqcup(1_-,1_+])/_{\!\sim}$. Then
\begin{eqnarray*}
\hat A_\fat|_{\overline U}&\sim&
\Big\{(f,g)\in C([0,1])\rtimes_{\omega_{[0,1]}}\Z^2\times C([1_-,1_\pm]\times\T)\ \Big|\ f|_1=g|_{1_-}\Big\}\\
&\sim&
C(\overline U)\rtimes_{\omega_{\overline U}} \Z^2,
\end{eqnarray*}
where we let $\omega_{\overline U}(s):=\omega_{1+s}$, for $s\in[0,1]$ and
$\omega_{\overline U}(s):=\omega_{2}$, for $s\in[1_-,1_\pm]$.
If we let $\hat\omega_{\overline U}(s):=\omega_{\frac{2}{1+s}}$, for $s\in[0,1]$ and
$\hat\omega_{\overline U}(s):=\omega_{1}$, for $s\in[1_-,1_\pm]$, then
$\omega_{\overline U}$ and $\hat\omega_{\overline U}$ are point-wise transverse
(cp. section \ref{Sec2DimNCT}).
Moreover, by construction of the action $\hat\alpha_\fat$ we have
\begin{eqnarray*}
(\hat A_\fat\rtimes_{\hat\alpha_\fat} \widehat\R^2)|_{\overline U}&\cong&
\K\otimes (C(\overline U)\rtimes_{\omega_{\overline U}\vee\hat\omega_{\overline U}} (\Z^2\times \widehat \R^2)).
\end{eqnarray*}
We have shown that $U$ is a chart for $(\hat A_\fat,\hat\alpha_\fat)$.
To show that $V$ is a chart we need one more step in between.
In fact, we use that we can extend the isomorphism 
$\eta^*$ (fibre-wise) to an isomorphism
\begin{eqnarray}
\label{EqTheEtaExtension}
\eta^*:C(\fat\times \T^2, \K)\to 
C(\fatn\times\T^2, \K), 
\end{eqnarray}
where $\fatn:=[0_-,0_+]:=\fat-1$ is just the intervall $\fat$
shifted by $1$.
This isomorphims turns the action
$\gamma|_\fat$ into
$ \gamma|_{\fatn}$,
which in each fibre $s\in\fatn$ just is
$\gamma|_{\fatn}|_s=\gamma|_0$.
Then we find
\begin{eqnarray*}
\xymatrix{
\hat A_\fat|_{\overline V}\ar[d]_-\cong\\
\Big\{(f,g)\in \K\otimes \big(C([0,\frac{1}{2}])\rtimes_{\omega_{[0,1]}}\Z^2)
\times C(\fat\times\T,\K)\big)\ \Big|\ f|_0=\eta^*g|_{1_+}\Big\}
\ar[d]_\cong^-{(f,g)\mapsto(f,\eta^*g)}
\\
\Big\{(f,h)\in \K\otimes \big(C([0,\frac{1}{2}])\rtimes_{\omega_{[0,1]}}\Z^2)
\times C(\fatn\times\T,\K)\big)\ \Big|\ f|_0=h|_{0_+}\Big\}
\ar[d]_-\cong\\
\K\otimes (C(\overline V)\rtimes_{\omega_{\overline V}} \Z^2)
}
\end{eqnarray*}
wherein $\omega_{\overline V}(s):=\omega_{1+s}$, for $s\in[0,\frac{1}{2}]$ and
$\omega_{\overline V}(s):=\omega_{1}$, for $s\in\fat$.
We have $\theta_\fat|_{\overline V}=[\omega_{\overline V}|_{\Z^2\times \Z^2}]$,
and we let $\hat\omega_{\overline V}(s):=\omega_{\frac{2}{1+s}}$, for $s\in[0,\frac{1}{2}]$ and
$\hat\omega_{\overline V}(s):=\omega_{2}$, for $s\in\fat$.
By construction
$\omega_{\overline V}$ and $\hat\omega_{\overline V}$ are point-wise transverse,
and  by construction of the action $\hat\alpha_\fat$ we have
\begin{eqnarray*}
(\hat A_\fat\rtimes_{\hat\alpha_\fat} \widehat\R^2)|_{\overline V}&\cong&
\K\otimes (C(\overline V)\rtimes_{\omega_{\overline V}\vee\hat\omega_{\overline V}} (\Z^2\times \widehat \R^2)).
\end{eqnarray*}

\end{proof}

\begin{rem}
\begin{enumerate}
\item 
The possibility of finding the extension
(\ref{EqTheEtaExtension})
is the crucial step in proving the local triviality
of $\hat A_\fat$.
This extension exists on a bundle of
commutative tori, but  could not have been found
if we had not thickened $1$ to $\fat$.
\item
The essential point for point-wise 
transversality is that the point-wise cocycles
$\omega_{[0,1]}(s):\Z^2\times \Z^2\to\UU$
have extensions
$\omega_{1+s}:\R^2\times \R^2\to\UU$
whose cohomology classes vary in the interval
 $[1,2]$ which (together with $[-2,-1]$) 
is the only interval of integer length  which is preserved 
by the map $s\mapsto\frac{2}{s}$.
The relevance of the second interval 
$[-2,-1]$ becomes clear below: It is the 
interval from which the (pre-)dual of 
$(\hat A_\fat,\hat\alpha_\fat)$
takes its Mackey obstructions.
\end{enumerate}
\end{rem}

Let us identify the pre-dual of the twisted 
Heisenberg bundle.
Knowing its local structure from the proof of
Proposition \ref{PropHeisenbergsCousin}
and knowing the calculation in section \ref{Sec2DimNCT} 
it is rather clear how the pre-dual looks like.
Let $(A_{\fat},\alpha_{\fat})$ to be the pullback in 
\begin{eqnarray*}
\xymatrix{
(A_{\fat}, \alpha_{\fat})\ar@{..>}[dd]\ar@{..>}[rr]&&
(\K\otimes (C([0,1])\rtimes_{\hat\omega_{[0,1]}}\Z^2 ),\delta)
\ar[dd]^{i_0^*\times i_1^*}\\
&&\\
(C(\fat\times\T^2,\K),\delta|_{\fat})\ar[rr]^-{(\eta^*\circ i_{1_+}^*)\times i_{1_-}^*}&&
(C(0\times\T^2,\K)\oplus C(1\times\T^2,\K),\delta|_0\times \delta|_1),
}
\end{eqnarray*}
where the action $\delta$ is point-wise given by 
$\delta|_t:= \lambda^{-(1+t)}\otimes \inf$,
$t\in[0,1]$,
the action $\delta|_{\fat}$ is just $\delta|_1$ in each fibre,
and the cocycle $\hat\omega_{[0,1]}$ is defined by
$\hat\omega_{[0,1]}(t):=\omega_{\frac{-2}{1+t}}\Big|_{N^\perp\times N^\perp}$, 
$t\in [0,1]$ and $N=\Z^2$.
Then one can prove that $(A_\fat,\alpha_1)$ 
is a (polarisable) NC pair over
$S^\fat$ just along the lines of Proposition \ref{PropHeisenbergsCousin}.
Moreover, the computation of its dual is a point-wise repetition
of section \ref{Sec2DimNCT}:
\begin{prop}
The NC pair $(A_\fat,\alpha_1)$ is dual to 
$(\hat A_\fat,\hat\alpha_1)$.
\end{prop}

Let $(\hat A_\fat,\theta_1):=\widehat\Theta (\hat A_\fat,\alpha_1)$
be the underlying NC bundle of the twisted Heisenberg bundle 
as explained in section \ref{NCbundle}.
Note that $\theta_1:S^\fat\to H^2(\Z^2,\UU)\cong\T$ has winding
number 1, and this is the key to proof what we have claimed in 
Proposition \ref{PropTheImages}.
Let us denote by
$F_n\to E$
a locally trivial $\K$-bundle on $E:=S^\fat\times \T^2$
such that $[F_n]\in \check H^3(E,\Z)\cong\Z$
corresponds to $n\in \Z$.
Then (up to isomorphism) all commutative NC pairs over $S^\fat$ 
are given by
$\Gamma(E,F_n)$ with 
suitable $\R^2$-actions.

\begin{prop}
$(\hat A_\fat,\theta_1)$ is not in the essential 
image of $\Xi|_\text{com}$.
\end{prop}

\begin{proof}
The proof makes use of the K-theory bundle introduced in
\cite{enoo}. This is the group bundle over $S^\fat$ whose fibres 
are given by the K-theory of the fibres.
According to \cite[Sec. 5]{enoo}, 
the $K_0$-bundle of the Heisenberg
bundle $\hat A_0$ is the $\Z^2$-bundle 
over $S^\fat$ which is given by 
gluing the trivial $\Z^2$-bundle over the interval
with the matrix 
\begin{eqnarray*}
M:=\left(
\begin{matrix}
1 &1 \\
0&1
\end{matrix}
\right).
\end{eqnarray*}
The Heisenberg bundle is a crossed product of $\Gamma(E,F_1)$ 
by a fibre-preserving $\R^2$-action, so by Connes' Thom isomorphism their
fibres have the same K-theory, and thus the K-theory bundles are the same
(see \cite[Theorem 3.5 and Remark 4.4]{enoo}). 
This implies that $\eta_*$ induces the above map $M$ in K-theory.
The twisted Heisenberg bundle then has a trivial $K_0$-bundle,
because it has the additional gluing with $\eta^*$ (rather than with $\eta_*$).
(The bundle $\hat A_{-\fat}$ which is obtained from the 
Heisenberg bundle with twist  $\eta_*$ has $K_0$-bundle
glued by $M^2$.)
Therefore, if $\Gamma(E,F_n)$ is a crossed product of the twisted Heisenberg bundle 
it also must have a trivial K-theory bundle. But as $F_n$ is obtained 
by gluing with $\eta^n$, the $K_0$-bundle of 
$\Gamma(E,F_n)$ is given by gluing the trivial $\Z^2$-bundle 
with $M^n$ which is trivial if and only if $n=0$.
Now, consider $\Gamma(E,F_0)\cong\K\otimes C(E)$ 
with any (transverse) $\R^2$-action $\alpha$.
By s crossed product 
$\K\otimes C(E)\rtimes_\alpha\R^2$
is equivariantly Morita equivalent to 
$\K\otimes C(S^1\times \{0\})\rtimes_{\alpha|_{\Z^2}}\Z^2$ with 
the $\R^2$-action that is inflated from the 
dual action of $\T^2$ (see the discussion in section \ref{SecNCToverPoint1}).
However, $\K\otimes C(S^\fat)\rtimes_{\alpha}\Z^2$ with 
the dual $\T^2$-action 
is a NCP-torus bundle in the sense of \cite{enoo} which has a trivial 
K-theory bundle. 
Then \cite[Theorem 7.2]{enoo} implies that the 
Mackey obstruction function $S^\fat\to H^2(\Z^2,\UU)\cong\T$ 
of $\alpha|_{\Z^2}\in{\rm Aut}(\K\otimes C(S^\fat\times \{0\}))$ 
is null-homotopic.
Since every null-homotopic map $S^\fat\to \T\cong H^2(\Z^2,\UU)$ can be lifted
to a continuous map $S^\fat\to \R\cong H^2(\R^2,\UU)$, this implies 
that there exists an action $\mu:\R^2\to {\rm Aut}(\K\otimes C(S^\fat\times \{0\}))$ 
such that $\mu|_{\Z^2}$ is Morita equivalent to $\alpha|_{\Z^2}$. 
But this implies
(see section \ref{SecNCToverPoint1})   that  $\alpha$
is Morita equivalent to the $\R^2$-action  $\mu\otimes \inf$ on $(\K\otimes C(S^\fat))\otimes C(\T^2)$. 
This means we have a global chart,
and so the Mackey obstruction function of $\alpha|_{\Z^2}$
coincides with the obstruction function $\theta_1$ defined by the
polarisable NC pair $(\K\otimes C(E),\alpha)$ which has winding 
number 1. This is a contradiction.
\end{proof}

\begin{appendix}

\section{Example \ref{ExaTwoDimAbiguity} }
\label{AppSiggisExample}
\noindent
Let $G=\R^2$ and $N=\Z^2$, and choose 
${\frac{1}{3}},{\frac{2}{3}}\in\R/\Z=\T\cong H^2(N,\UU)$. 
Denote by $\hat\mu_3$ an action of  $\widehat{G}=\R^2$ on $\K$ with
Mackey obstruction $3\in \R\cong H^2(\R^2,\UU)$.
Then there is a $\widehat G$-equivariant Morita equivalence
$$
\Big(\K\otimes \C\rtimes_{\frac{2}{3}}N, \hat\mu_{3}\otimes {\rm inf} \Big)
\sim
\Big(\K\otimes \C\rtimes_{\frac{1}{3}}N, {\rm id}\otimes {\rm inf} \Big).
$$

\begin{proof} {\bf Step 1:}
The symmetry groups $S\subset N$ of the two classes ${\frac{1}{3}},{\frac{2}{3}}$ agree, 
i.e. the maps $h_{\frac{1}{3}},h_{\frac{2}{3}}:N\to \widehat N$ have the same kernel 
$S=(3\Z)^2$. The annihilator of $S$ in $\widehat G=\R^2$ is $S^\perp=({\frac{1}{3}}\Z)^2$,
and the annihilator of $S$ in $\widehat N$ is $S^{\perp}_N=({\frac{1}{3}}\Z/\Z)^2$.
We are in the following situation:
$$
\xymatrix{
&&N^\perp\ar@{^(->}[d]&&\\
&&S^\perp\ar@{^(->}[r]\ar@{->>}[d]&\widehat G\ar@{->>}[d]\ar@{->>}[r]&\widehat G/S^\perp\ar@{=}[d]\\
\widehat{N/S}\ar@{=}[r]&S^\perp/N^\perp\ar@{=}[r]&S^\perp_N\ar@{^(->}[r]&\widehat N\ar@{->>}[d]\ar@{->>}[r]&\widehat N/S^\perp_N\\
&&&\widehat S\ar@{=}[ru]&.
}
$$
The quotient map $N=\Z^2\to(\Z/3\Z)^2=N/S$ induces a map in cohomology,
and let us denote by $\frac{1}{3}$ also the pre-image of $\frac{1}{3}$ in 
\begin{eqnarray*}
\xymatrix{
\frac{1}{3}\in\ar@<3.6ex>@{|->}[d]
\hspace{-1.4cm}
&\frac{1}{3}\Z/\Z\ar@{^(->}[d]^{\text{inclusion}}\ar[r]^-\cong
&H^2(N/S,\UU)\ar[d]^{\text{induced map}}\\
\quad\frac{1}{3}\hspace{0.1cm}\in
\hspace{-1.2cm}
&\T\ar[r]^-\cong
&H^2(N,\UU).
}
\end{eqnarray*}
The action of $\widehat G$ on 
$\widehat S\cong{\rm Prim}(\C\rtimes_\frac{1}{3} N)
\cong{\rm Prim}(\C\rtimes_\frac{2}{3} N)$
has stabiliser $S^\perp$ which acts via the quotient
$S^\perp\to S^\perp_N=\widehat{N/S}$ and 
the dual actions on $\C\rtimes_\frac{1}{3} N/S$
and $\C\rtimes_\frac{2}{3} N/S$ on these two algebras.
Let us denote these actions by $\sigma$ and $\tau$, respectively.
Then there are $\widehat G$-equivariant isomorphisms
\begin{align*}
\Big(\K\otimes (\C\rtimes_\frac{2}{3} N)&,\hat\mu_3\otimes \inf\Big)\\
&\cong
\Big(\K\otimes\Ind^{\widehat G}_{S^\perp}(\C\rtimes_\frac{2}{3} N/S,\tau),
\hat\mu_3\otimes \Ind(\tau)\Big)
\\
&\cong
\Big(\Ind^{\widehat G}_{S^\perp}(\K\otimes(\C\rtimes_\frac{2}{3} N/S),\ 
\hat\mu_3|_{S^\perp}\otimes \tau),
\Ind(\hat\mu_3|_{S^\perp}\otimes \tau)\Big)
\end{align*}
and
\begin{eqnarray*}
\Big(\K\otimes (\C\rtimes_\frac{1}{3} N),\id \otimes \inf\Big)&\cong&
\Big(\K\otimes\Ind^{\widehat G}_{S^\perp}(\C\rtimes_\frac{1}{3} N/S,\sigma),
\id\otimes \Ind(\sigma)\Big)
\\
&\cong&
\Big(\Ind^{\widehat G}_{S^\perp}(\K\otimes(\C\rtimes_\frac{1}{3} N/S),\ 
\id\otimes \sigma),
\Ind(\id\otimes \sigma)\Big).
\end{eqnarray*}
\\
{\bf Step 2:} 
The product in $\C\rtimes_\frac{1}{3} N/S$ is given 
by 
$$
(f*f')(n)=\sum_{m\in N/S}f(m)f'(n-m)\exp\left(2\pi i\frac{1}{3} m_2(n_1-m_1)\right).
$$
It is a straight forward computation to see that the mapping which assigns
to a function
$f:N/S\to \C$ the matrix $\tilde f$ with entries
\begin{eqnarray*}
\tilde f(a,b):= \sum_{c=0,1,2}f(c,b-a)\exp\left(2\pi i\frac{1}{3}ca\right),
\quad a,b=0,1,2,
\end{eqnarray*}
defines an isomorphism
$
\ \tilde{}:\C\rtimes_\frac{1}{3} N/S \cong M_3(\C)= \K(L^2(\Z/3\Z)).
$
The action $\sigma$ transforms under this isomorphism 
to a conjugation action $\tilde \sigma:S^\perp\to \Aut(\K(L^2(\Z/3\Z)))$, 
i.e.
$\tilde\sigma_s={\rm Ad}(V(s))$, where the unitary 
$V(s)\in \U(L^2(\Z/3\Z))$
is given by
$$
(V(s)\psi)(a)=\exp\left(2\pi i\frac{1}{3}c(a+d)\right)^{-1}\ \psi(a+d), 
$$
for $s=(\frac{1}{3}c,\frac{1}{3}d)\in  S^\perp=(\frac{1}{3}\Z)^2,
\psi\in L^2(\Z/3\Z)$.
One can then start calculating 
$$
(\partial V)(s,s')=V(s')V(s+s')^{-1}V(s)= \exp\left(2\pi i\frac{1}{3}dc' \right)
= \exp\left(2\pi i\frac{2}{3}dc' \right)^{-1}.
$$
So after identifying $\frac{1}{3}\cdot\_:\Z^2\cong (\frac{1}{3}\Z)^2=S^\perp$, 
we find  the Mackey obstruction of ${\tilde\sigma}$ satisfying
$$
\xymatrix{
H^2(S^\perp,\UU)\ar[rr]_-{(\frac{1}{3}\cdot\_)^*}^-\cong&& H^2(\Z^2,\UU)\ar[r]^-\cong&\T\\
{\rm Ma}(\tilde\sigma)\ar@{}[u]|{\rotatebox{90}{\resizebox{0.4cm}{!}{$\in$} }}
\ar@{|->}[rrr]
&&&\frac{2}{3} .
\ar@{}[u]|{\rotatebox{90}{\resizebox{0.4cm}{!}{$\in$} }}
}
$$
\\\\
{\bf Step 3:}
The Mackey obstruction of $\hat\mu_{3}|_{S^\perp}$ is given 
by the cocycle
$$
\omega(s,s')=\exp\left(2\pi i\ 3\ s_2 s'_1\right)=\exp\left(2\pi i \frac{1}{3}dc'\right),
$$
{for } $s=(s_1,s_2)=(\frac{1}{3}c,\frac{1}{3}d),
s'=(s'_1,s'_2)=(\frac{1}{3}c',\frac{1}{3}d')\in S^\perp.
$
So here we find
$$
\xymatrix{
H^2(S^\perp,\UU)\ar[rr]_-{(\frac{1}{3}\cdot\_)^*}^-\cong&& H^2(\Z^2,\UU)\ar[r]^-\cong&\T\\
{\rm Ma}(\hat\mu|_{S^\perp})\ar@{}[u]|{\rotatebox{90}{\resizebox{0.4cm}{!}{$\in$} }}
\ar@{|->}[rrr]
&&&\frac{1}{3}  .
\ar@{}[u]|{\rotatebox{90}{\resizebox{0.4cm}{!}{$\in$} }}
}
$$
{\bf Step 4:}
A similar isomorphism as found in Step 2 can be used to 
identify $\C\rtimes_\frac{2}{3} N/S$  also with 
$\K(L^2(\Z/3\Z))$. In fact,
the product in $\C\rtimes_\frac{2}{3} N/S$ agrees 
with the product in  $\C\rtimes_\frac{2}{3} N/S$
up to a sign:
\begin{eqnarray*}
(f*f')(n)&=&\sum_{m\in N/S}f(m)f'(n-m)\exp\left(2\pi i\frac{2}{3} m_2(n_1-m_1)\right)\\
&=&\sum_{m\in N/S}f(m)f'(n-m)\exp\left(2\pi i\frac{1}{3} m_2(n_1-m_1)\right)^{-1}.
\end{eqnarray*}
Following this sign in the construction made in step 2
one finally finds that the Mackey obstruction of $\tau$
is exactly the inverse of $\sigma$:
$$
\xymatrix{
H^2(S^\perp,\UU)\ar[rr]_-{(\frac{1}{3}\cdot\_)^*}^-\cong&& H^2(\Z^2,\UU)\ar[r]^-\cong&\T\\
{\rm Ma}(\tau)\ar@{}[u]|{\rotatebox{90}{\resizebox{0.4cm}{!}{$\in$} }}
\ar@{|->}[rrr]
&&&\frac{1}{3}.
\ar@{}[u]|{\rotatebox{90}{\resizebox{0.4cm}{!}{$\in$} }}
}
$$
{\bf Step 5:}
Summing up all the Mackey obstructions,
we see that the two C*-dynamical systems
$$
(\K\otimes \C\rtimes_\frac{2}{3}N/S, S^\perp, \hat\mu_3|_{S^\perp}\otimes\tau)
\text{ and }
(\K\otimes \C\rtimes_\frac{1}{3}N/S, S^\perp, \id\otimes\sigma)$$
are equivalent.
This implies that also their induced systems are equivalent
which proves the claim.
\end{proof}

\section{Lemma \ref{LemCommutes}}
\label{AppCommutes} 
\noindent
The canonical isomorphism  
$H^2(G,L^\infty(G/N,\UU)\to H^2(N,\UU)$
of (\cite[Thm. 6]{Moore3})
makes the diagram 
\begin{eqnarray}
\label{TheLastCommDia}
\xymatrix{
&H^2(G,\UU)\ar[d]\ar[dr]&\\
&H^2(N,\UU)&H^2(G,L^\infty(G/N,\UU))\ar[l]_-\cong
}
\end{eqnarray}
commute, where the vertical map is restriction and 
the diagonal map is induced by the inclusion of coefficients.

\begin{proof}
Let us introduce some notation used in (\cite{Moore3}).
Let $I(X):=\{N\to X\}$ for any $X$. It has 
the structure of an $N$-module by left 
translation: $n\cdot f:=f(\_-n)\in I(N)$.
Let us denote by $A$ the quotient of
$I(\UU)$ by the constants $\UU$,
i.e. we have an $N$-equivariant short exact sequence
\begin{eqnarray}
\label{AShortExactOne}
1\to \UU\to I(\UU)\to A\to 1,
\end{eqnarray}
where $\UU$ has the trivial $N$-structure 
and $A$ has the quotient structure.
We obtain an $N$-equivariant embedding
$i:A\to I(A)$  by $(ia)(n):=(-n)\cdot a$, 
for $a\in A, n\in N$. The quotient of $I(A)$ by $i(A)$
is denoted by $U(A)$, i.e. we have another 
$N$-equivariant sequence
\begin{eqnarray*}
1\to A\to I(A)\to U(A)\to 1.
\end{eqnarray*}
By definition (cp. the axioms of group cohomology
as a derived functor (\cite[Sec. 4]{Moore3})) 
we have an exact sequence
\begin{eqnarray*}
1\to A^N\to I(A)^N\to U(A)^N\to H^1(N,A)\to 0,
\end{eqnarray*}
where the exponent denotes taking invariants.
Let $I^G_N(Y)$ denote the induced $G$-module
for any $Y$-module, i.e. equivalence classes of functions 
$f:G\to Y$, such that for all $n\in N$ 
$(f(g-n)= n\cdot (f(g))$ holds for almost all $g\in G$.
Two functions are identified if they agree almost everywhere.
$I^G_N(Y)$ is a $G$-module by left translation.
Note that $I^G_N(\UU)=L^\infty(G/N,\UU)$ as $G$-modules.
The functor $I^G_N(\_)$ from $N$-modules to $G$-modules 
is exact (\cite[Proposition 19]{Moore3}) and 
again by the axioms and (\cite[Proposition 18]{Moore3}) 
we have an exact sequence
\begin{eqnarray*}
1\to I^G_N(A)^G\to I^G_N(I(A))^G\to I^G_N(U(A))^G\to H^1(G,I^G_N(A))\to 0.
\end{eqnarray*}
The canonical isomorphism
$H^2(N,\UU)\to H^2(G,I^G_N(\UU))$ is defined by the diagram:
\begin{eqnarray*}
\xymatrix{
&1\ar[d]&1\ar[d]&\\
&A^N\ar[d]\ar[r]^-\cong&I^G_N(A)^G\ar[d]&\\
&I(A)^N\ar[d]\ar[r]^-\cong&I^G_N(I(A))^G\ar[d]&\\
&U(A)^N\ar[d]\ar[r]^-\cong&I^G_N(U(A))^G\ar[d]&\\
H^2(N,\UU)&\ar[l]_-\cong H^1(N,A)\ar[d]\ar@{..>}[r]^-\cong&H^1(G,I^G_N(A))\ar[d]
\ar[r]^-\cong&
H^2(G,I^G_N(\UU))\\
&0&0&
}.
\end{eqnarray*}
Here the top three isomorphisms are given by sending an 
invariant element $\nu$ to the constant function $\bar \nu:g\mapsto\nu$
(cp. \cite[Proposition 19]{Moore3}).
The dotted isomorphism is induced by the one above and the two 
isomorphisms to the left and right are the connecting morphisms in  
the long exact sequence induced by (\ref{AShortExactOne}).

Let us describe the dotted arrow. To do so, start with
some element $x\in H^1(N,A)$ and choose some $\nu\in U(A)^N$
which maps to $x$. Following the isomorphism to the 
right we get $\bar\nu\in I^G_N(U(A))^G\subset I^G_N(U(A))$.
Let us denote by $\tilde\nu\in I^G_N(I(A))$ an element 
that maps to $\bar\nu$. Then the image $y$ of $x$ under the
dotted arrow is (the class of) $d_G\tilde\nu:G\to I^G_N(A)$
defined by 
$$
G\ni g\to \tilde\nu(\_-g)\ \tilde\nu(\_)^{-1}\in 
I^G_N(i(A))\cong I^G_N(A).
$$

Let us now describe the element $x$ in terms of $\tilde\nu$:
We consider for each $n\in N$ a co-null set 
$S_n\subset G$ such that 
$\tilde\nu(g-n)=n\cdot \tilde\nu(g)$ holds for 
all $g\in S_n$.
Then let $S^1:=\bigcap_{n\in N}S_n$.
So  $\tilde\nu(g-n)=n\cdot \tilde\nu(g)$
holds for all $n\in N$ and all $g\in S^1$.
Moreover, there is another co-null set $S^2\subset G$ such that
$\tilde\nu(g)\textrm{ mod }A=\nu$ for all $g\in S^2$. 
Let $s\in S^1\cap S^2$, then 
$\tilde\nu(s)\in I(A)$ is a lift of $\nu\in U(A)^N$, 
and $x$ is represented by the class of 
$d_N(\tilde\nu(s)):N\to A$ defined by
$$
N\ni \mapsto n\cdot(\tilde\nu(s))\ \tilde\nu(s)^{-1}
=\tilde\nu(s-n)\ \tilde\nu(s)^{-1} \in i(A)\cong A.
$$
Note that that the operation restriction to $N$ and evaluation 
at $s$ transforms $d_G\tilde\nu$ into $d_N(\tilde\nu(s))$:
$(d_G\tilde\nu)(n)(s)=d_N(\tilde\nu(s))(n)$.

Let $\sigma:A\to I(\UU)$ be a Borel section of the 
quotient map. Then the element $[\omega_{N,\nu,s}]$
in $H^2(N,A)$ corresponding to $x=[d_n(\tilde\nu(s))]$
is given by 
$$
\omega_{N,\tilde\nu,s}(n,m) =  \sigma(d_N(\tilde\nu(s))(m))\ \sigma(d_N(\tilde\nu(s))(n+m))^{-1}\
\sigma(d_N(\tilde\nu(s))(n)).
$$
The element $[\omega_{G,\tilde\nu}]$ in $H^2(G,I^G_N(\UU))$ corrsponding 
to $y=[d_G\tilde\nu]$ is given by
$$
\omega_{G,\tilde\nu}(g,h) = \sigma_*((d_G\tilde\nu)(h)(\_-g))\ 
\sigma_*((d_G\tilde\nu)(g+h)(\_))^{-1}\
\sigma_*((d_G\tilde\nu)(g)(\_))
$$
where $\sigma_*:I^G_N(A)\to I^G_N(I(\UU))$ is just composition with
$\sigma$. 
Note again that that the operation restriction to $N$ and evaluation 
at $s$ transforms $\omega_{G,\tilde\nu}$ into $\omega_{N,\tilde\nu,s}$:
$\omega_{G,\tilde\nu}(n,m)(s)=\omega_{N,\tilde\nu,s}(n,m)$.

Now let us turn to the commutativity of 
(\ref{TheLastCommDia}).
Let $[\omega]\in H^2(G,\UU)$, and consider its 
image $y$ (also given by $\omega$) in 
$H^2(G,I^G_N(\UU))$. We can represent $y$ 
as above by finding some $\nu,\tilde\nu$ such that for all $g,h$
$$
\omega(g,h) \ (d_Gc)(g,h)= \omega_{G,\tilde\nu}(g,h)
$$
holds in $I^G_N(\UU)$, for some cochain $c:G\to I^G_N(\UU)$. 
So there is another co-null set $S^3$ such that 
$$
\omega(n,m) \ (dc)(n,m)(s)= \omega_{G,\tilde\nu}(n,m)(s)
$$
holds for all $n,m\in N$ and all $s\in S^3$.
If we choose $s_0\in S^1\cap S^2\cap S^3$, then 
by the above construction we have that the image of $[\omega]$
along the composition in diagram (\ref{TheLastCommDia}) is given by
the cocycle
\begin{eqnarray*}
\omega_{N,\tilde\nu,s_0}(n,m)&=&\omega_{G,\tilde\nu}(n,m)(s_0)\\
&=&\omega(n,m)\ c(m)(s_0)\ c(n+m)(s_0)^{-1}\ c(n)(s_0).
\end{eqnarray*}
It follows that  $\omega_{N,\tilde\nu,s_0}\sim \omega|_{N\times N}$ by the
cochain $n\mapsto c(n)(s_0)$.
This proves the lemma.
\end{proof}

\end{appendix}

\end{document}